\documentclass[11pt]{amsbook}
\usepackage{amsmath,amsfonts,amsthm,amssymb,amscd,graphicx}
\usepackage[spanish]{babel}
\def\classification#1{\def\@class{#1}}
\classification{\null}

\DeclareFontFamily{OT1}{rsfs}{}
\DeclareFontShape{OT1}{rsfs}{n}{it}{<-> rsfs10}{}
\DeclareMathAlphabet{\mathscr}{OT1}{rsfs}{n}{it}

\DeclareMathOperator{\mcd}{mcd}

\DeclareMathOperator{\Var}{Var}

\DeclareMathOperator{\Prob}{Prob}

\newtheorem{thm}{Teorema}[chapter]

\theoremstyle{definition}

\newtheorem{xca}[thm]{Ejercicio}

\theoremstyle{remark}

\numberwithin{section}{chapter}
\numberwithin{equation}{chapter}
\numberwithin{equation}{section}
\makeindex
\begin{document}
\frontmatter

\title{Azar y aritm\'etica\\
{\small Un cap\'itulo de la teor\'ia probabil\'istica de n\'umeros}}
\author{Harald Andr\'es Helfgott}
\address{H. A. Helfgott, School of Mathematics, University of Bristol, Bristol, BS8 1TW, United Kingdom}
\email{harald.helfgott@gmail.com}
\date{\today}
\begin{abstract}
Sea $\omega(n)$ el n\'umero de divisores primos de un entero
$n$. Sea $n$ un entero tomado al azar entre $1$ y $N$. Qu\'e se puede decir del
valor que entonces tomar\'a $\omega(n)$? Cu\'al es su esperanza? Cu\'al es su
distribucion en el limite? Cu\'al es la probabilidad que $\omega(n)$ tome
valores que se alejen mucho de su esperanza?

Estudiamos estas preguntas a guisa de introducci\'on a la teor\'ia de n\'umeros
probabil\'istica. Trataremos varios t\'opicos centrales de la teor\'ia
de probabilidades sin suponer conocimientos previos en el \'area.
No asumiremos ni teor\'ia de la medida ni an\'alisis complejo.
En los ejercicios, entre otros t\'opicos, se desarrollar\'an las
bases de la teor\'ia de cribas como una aplicaci\'on de ideas probabil\'isticas.

\end{abstract}

\maketitle

%    Dedication.  If the dedication is longer than a line or two,
%    remove the centering instructions and the line break.
\cleardoublepage
\thispagestyle{empty}
\vspace*{13.5pc}
\begin{center}
Al miaj gepatroj
\end{center}

\setcounter{page}{3}

%    Include unnumbered chapters (preface, acknowledgments, etc.) here.
%\include{notacion}

\tableofcontents
\begin{center} {\bf Prefacio} \end{center}

Este es un estudio de los factores primos de un n\'umero tomado al azar. 
El objeto principal es servir de introducci\'on a la teor\'ia probabil\'istica
de n\'umeros y, al mismo tiempo, a varios temas centrales en la teor\'ia de
las probabilidades en general:
la varianza, el l\'imite central, las grandes desviaciones, la entrop\'ia.

La historia de la teor\'ia de n\'umeros probabil\'istica comienza
 con Hardy\index{Hardy, G. H.} y
Ramanujan\index{Ramanujan, S.} \cite{HR}, quienes fueron los primeros en 
analizar el tema central de este libro:
 la distribuci\'on del n\'umero $\omega(n)$ de divisores primos
de un n\'umero entero aleatorio $n$. En el curso de la generaci\'on siguiente -- notablemente con el teorema de Erd\"os\index{Erd\"os, P.}
y Kac\index{Kac, M.} \cite{EK}, el cual estudiaremos
en la secci\'on \ref{sec:limcen}
-- se fueron asimilando conceptos y t\'ecnicas
 de la teor\'ia de  probabilidades en general al estudio incipiente del tema.
El \'area ha seguido desarroll\'andose hasta nuestros d\'ias, gracias tanto a especialistas en teor\'ia de n\'umeros como a probabilistas. 

No asumiremos ning\'un conocimiento de an\'alisis complejo ni de teor\'ia de 
la medida.
 Al final de cada secci\'on, se encontrar\'a una serie de notas y problemas; esencialmente se trata de ejercicios guiados o 
esbozos de pruebas a seguir y completar con l\'apiz y papel. Entre otros t\'opicos, las notas de fin de secci\'on desarrollan las bases de la teor\'ia de 
cribas, tanto como una aplicaci\'on de conceptos probabil\'isticos, como para
uso en el texto principal. Mi objetivo
 ha sido dar las pruebas que me parecen ser las 
m\'as naturales, antes que las m\'as conocidas.

El texto presente est\'a basado en las notas de clase de un curso que 
dict\'e en Julio y Agosto de 2007 bajo los auspicios del IMCA (Instituto 
de Matem\'atica y Ciencias Afines) en Lima, Per\'u.
Agradezco tanto al IMCA
como a la Universidad Mayor de San Marcos por su hospitalidad.

\begin{center} {\bf Notaci\'on} \end{center}

Sean $d$ y $n$ n\'umeros enteros. 
Escribimos $d|n$ cuando queremos decir que $d$ divide a $n$ exactamente,
es decir, sin dejar resto: $3|6$, $5|15$, $1|n$ para
todo $n$. Escribimos $d\nmid n$ cuando $d$ no divide a $n$, es decir,
cuando la division de $n$ por $d$ deja resto:
$4\nmid 6$, $7\nmid 15$, $(n+1)\nmid n$ para todo $n$.

La letra $p$ siempre designar\'a a un n\'umero primo. La funci\'on
$\Lambda(n)$ (funci\'on de von Mangoldt\index{funci\'on de von Mangoldt}) 
se define como sigue:
\[\Lambda(n) = \begin{cases} \log p &\text{si $n = p^{\alpha}$ para
alg\'un primo $p$ y alg\'un entero $\alpha>0$}\\
0&\text{si no es as\'i.}\end{cases}\]

Denotamos por $\lfloor x\rfloor$\index{$\lfloor x\rfloor$}
el m\'aximo entero $n$ que no sea mayor que $x$. Por ejemplo,
$\lfloor 2.75\rfloor = 2$, $\lfloor 7\rfloor = 7$, $\lfloor \pi\rfloor = 3$.

Cuando decimos ``logaritmo'' o escribimos $\log x$, tenemos siempre
en mente al logaritmo en base $e$, a menos que otra base se especifique
explicitamente (``logaritmo en base $2$'', por ejemplo).
Al contrario de los escritores franceses, utilizaremos la
notaci\'on $\log_2 x$ para
el logaritmo base $2$ de $x$, y $\log \log x$ para el logaritmo (base $e$)
del logaritmo (base $e$) de $x$.

Utilizaremos la notaci\'on $O$, $o$ de Landau
dadas dos funciones $f$, $g$,
(a) se escribe $f(x) = O(g(x))$ cuando existen constantes $c_1, c_2>0$ 
tales que $|f(x)/g(x)|<c_1$ para todo $x>c_2$; (b) se escribe 
$f(x) = o(g(x))$ cuando $\lim_{x\to \infty} f(x)/g(x) = 0$. Est\'a claro que
$f(x) = o(g(x))$ implica $f(x) = O(g(x))$, pero no viceversa.
 Ejemplos: $x^2 = O(x^3)$,
$x^2 = o(x^3)$, $f(x) = O(f(x))$ (para todo $f$), 
$\sin x = O(1)$, $x \sin x = O(x)$, $\sum_{n\leq x} 1/n = O(\log n)$,
$\sum_{n\leq x} (-1)^n/n = O(1)$, $\prod_{n\leq x} (1-1/n) = o(1)$.
En particular, $f(x) = O(1)$ quiere decir que $f$ est\'a acotada por una
constante, y $f(x) = o(1)$ quiere decir que $f$ tiende a cero cuando
$x$ va al infinito. Escribimos $O_{c}(1)$, $o_{\delta,z}(1)$ si dichas constantes
dependen de $c$ o $\delta$ y $z$, por ejemplo.

La expresi\'on ``$f(x)\ll g(x)$'' es un sin\'onimo de ``$f(x) = O(g(x))$'';\index{$\ll$} 
la expresi\'on ``$f(x)\gg g(x)$'' es un sin\'onimo de ``$g(x) = O(f(x))$''.\index{$\gg$}

Escribimos $f(x) \sim g(x)$ cuando queremos decir que $f$ es asint\'otica
con respecto a $g$, i.e., \[\lim_{x\to \infty} f(x)/g(x) = 1.\] Si decimos que $f(x) = o(g(x))$
(o $f(x) \sim g(x)$) ``cuando $x\to 0$'', queremos decir que
$\lim_{x\to 0} f(x)/g(x) = 0$ (o, respectivamente, 
$\lim_{x\to 0} f(x)/g(x) = 1$).

Denotamos por $\Prob(E)$ la probabilidad del evento aleatorio $E$, por
$\mathbb{E}(X)$ la esperanza de la variable aleatoria $X$ y por
$\Var(X)$ la varianza de la variable $X$.\index{esperanza}
\index{varianza}
Ver el ap\'endice \ref{ap:proba}.

\mainmatter
%    Include main chapters here.
%una introduccion: modelos, metodos, aplicaciones?
%Erdos; Cramer
 \chapter{Los divisores primos}
%O: los divisores de un n\'umero aleatorio?
\section{La esperanza}\label{sec:esperanza}
{\bf Definici\'on de esperanza. Ejemplos.}
Recordemos que, si se tiene una variable aleatoria $X$ que toma los valores
$x_1, x_2, \dotsc, x_n$ con probabilidades $p_1$, $p_2$, \dots, $p_n$, donde
$p_1+p_2+\dotsc+p_n=1$, entonces la esperanza\index{esperanza} se define como la cantidad
\[\sum_{i=1}^n p_i x_i .\]
As\'i, por ejemplo, si $X$ es el valor que da un dado arrojado al aire,
\[X= \begin{cases} 1 & \text{con probabilidad $1/6$}\\
2 & \text{con probabilidad $1/6$}\\
3 & \text{con probabilidad $1/6$}\\
\dotsc & \dotsc\\
6 & \text{con probabilidad $1/6$}\end{cases}\]
y por lo tanto
\[\mathbb{E}(X) = \sum_{i=1}^n p_i x_i = \frac{1}{6} \cdot 1 + \frac{1}{6} 
\cdot 2 + \dotsc + \frac{1}{6} \cdot 6 = \frac{7}{2} .\]
Si $X$ es un dado trucado, muy bien podr\'ia tener la distribuci\'on
\[X = \begin{cases} 1 &\text{con probabilidad $1/18$}\\
2,3,\dotsc, 5 &\text{con probabilidad $1/9$ en cada caso}\\
6 &\text{con probabilidad $1/2$}
\end{cases}\]
y entonces \[\mathbb{E}(X) = \sum_i p_i x_i = \frac{1}{18} \cdot 1 +
\frac{1}{9} \cdot 2 + \dotsc + \frac{1}{9} \cdot 5 + \frac{1}{2} \cdot 6 =
 \frac{83}{18}.\]
{\bf Esperanza y sumas.}
Denotamos por $\mathbb{E}(X)$ la esperanza de una variable aleatoria $X$.
Sean $X_1, X_2,\dotsc, X_n$ variables aleatorias. Es f\'acil ver que
\begin{equation}\label{eq:sumo}
\mathbb{E}(X_1 + \dotsc + X_n) = \mathbb{E}(X_1) + \dotsc +
\mathbb{E}(X_n).\end{equation}

{\em Aplicaci\'on 1.}
Denotemos por $\tau(n)$\index{$\tau(n)$} el n\'umero de divisores de un entero $n$.\index{divisores} Cu\'anto
es $\tau(n)$, en promedio?

Para que nuestra pregunta tenga sentido, debemos decir como estamos
escogiendo $n$. Fijemos un entero $N$. Tomamos $n$ al azar entre
$1$ y $N$, con la distribuci\'on uniforme. Estamos preguntando cu\'al
es el valor de $\mathbb{E}(\tau(n))$.

Para todo entero $m$, definimos
\begin{equation}\label{eq:maz}
X_m = \begin{cases}0 &\text{si $m\nmid n$,}\\ 1 &\text{si
    $m|n$.}\end{cases}\end{equation}
(Formalmente, lo que tenemos es una variable aleatoria $X$ que toma
valores $n$ entre $1$ y $N$ con la distribuci\'on uniforme, y
varias variables aleatorias $X_m$ que dependen de $X$.) Ahora bien,
\[\sum X_m = \tau(n).\]
Est\'a claro que lo que queremos calcular es $\mathbb{E}(\sum X_m)$.

Ya sabemos (\ref{eq:sumo}) que
\[\mathbb{E}\left(\sum_m X_m\right) = \sum_m \mathbb{E}(X_m).\]
Ahora bien, cu\'al es el valor de $\mathbb{E}(X_m)$? Calculamos:
\begin{equation}\label{eq:thecount}
\mathbb{E}(X_m) = \Prob(X_m = 1) = \frac{1}{N} 
\mathop{\sum_{n\leq N}}_{m|n} 1 = \frac{1}{N} 
\left\lfloor \frac{N}{m}\right\rfloor
= \frac{1}{N} \left(\frac{N}{m} + O(1)\right)  = \frac{1}{m} + O\left(
\frac{1}{N}\right).\end{equation}
(Aqu\'i, como de ahora en adelante, $O(1)$ quiere decir ``una cantidad
$x$ tal que $|x|\leq C$ para alguna constante $C$'' y $O(1/N)$ quiere decir 
``una cantidad $x$ tal que $|x|\leq C/N$''.
La ecuaci\'on $\lfloor N/m \rfloor = N/m + O(1)$ nos est\'a diciendo
simplemente que el valor absoluto de
 la diferencia entre $\lfloor N/m \rfloor$ y $N/m$ es siempre menor que
una constante -- en verdad, menor que $1$.)\index{$O(f), o(f)$}
Por lo tanto
\[\mathbb{E}\left(\sum_m X_m\right) = \sum_m \mathbb{E}(X_m) =
\sum_{m\leq N} \frac{1}{m} + \sum_{m\leq N} O\left(\frac{1}{N}\right) =
\log N + O(1).\]
concluimos que \begin{equation}\label{eq:gotz}
\mathbb{E}(\tau(n)) = \log N + O(1).\end{equation}

{\em Aplicaci\'on 2.} Sea $\omega(n)$\index{$\omega(n)$} el n\'umero de divisores primos
de $n$. Cu\'anto es $\omega(n)$, en promedio?\index{divisores!divisores primos}

Calculemos:
\[\mathbb{E}(\omega(n)) = \mathbb{E}\left(\sum_{p\leq N} X_p\right)
= \sum_{p\leq N} \mathbb{E}(X_p) = \sum_{p\leq N} \frac{1}{p}
+ \sum_{p\leq N} O\left(\frac{1}{N}\right)
.\]
Ahora bien
\begin{equation}\label{eq:ort}
\sum_{p\leq N} \frac{1}{p} = \log \log N + O(1) \end{equation}
(teorema de Chebyshev-Mertens, 1875; ver los ejercicios).\index{teorema!de
Chebyshev-Mertens}
Por lo tanto
\begin{equation}\label{eq:ost}\mathbb{E}(\omega(n)) = \log \log N + O(1).
\end{equation}

{\em Aplicaci\'on 3.} 
Cu\'antos factores primos de un tama\~no dado tiene un n\'umero tomado
al azar?

Precisemos el rango.
Sean dados $\delta_0$, $\delta_1$ tales que $0\leq \delta_0<\delta_1<1$.
Tomemos $n$ entre $1$ y $N$ bajo la distribuci\'on uniforme. Queremos
saber la esperanza $\mathbb{E}(Y)$ del n\'umero $Y$
de factores primos de $n$ entre $N^{\delta_0}$ y
$N^{\delta_1}$.

Calculemos:
\[\begin{aligned}
\mathbb{E}(Y) &= \mathbb{E}\left(\sum_{N^{\delta_0}< p \leq N^{\delta_1}} X_p
\right) = \sum_{N^{\delta_0}<p\leq N^{\delta_1}} \left(\frac{1}{p} + O\left(
\frac{1}{N}\right)\right)\\ &=
\log \delta_1 - \log \delta_0 + o(1) .\end{aligned}\]

Podemos plantearnos una segunda pregunta:
 cu\'al es la probabilidad que $n$ tenga por lo menos
un factor primo entre $N^{\delta_0}$ y $N^{\delta_1}$? 
Un n\'umero $n\leq N$ no puede tener m\'as de $1/\delta_0$
factores primos mayores que $N^{\delta_0}$. Por lo tanto,
\[\Prob(Y>0) \leq \frac{1}{1/\delta_0} \mathbb{E}(Y) \leq
\delta_0 \cdot (\log \delta_1 - \log \delta_0 + o(1)) .\]

Esto es s\'olo una cota superior.
 Veremos m\'as tarde como estimar $\Prob(Y>0)$ de manera m\'as precisa.

{\bf Desigualdad de Markov.}\index{desigualdad!de Markov}
Sea $X$ una variable aleatoria que toma
siempre valores no negativos. Sea $t\geq \mathbb{E}(X)$. Entonces
\begin{equation}\label{eq:mark}
\Prob(X\geq t) \leq \frac{\mathbb{E}(X)}{t}
\text{\;\;\;\;\;\;\;\;\;\;\;\;(desigualdad de Markov).}\end{equation}
Esto tiene sentido: si, en promedio, cae $1$ cm de lluvia al d\'ia, la
probabilidad que caigan m\'as de $10$ cm no puede ser m\'as de $0.1$.
(Por otra parte, dado que cae $1$ cm de lluvia al d\'ia en promedio, 
la probabilidad que caigan $0$ cm puede ser tan
cercana a $1$ como se quiera: muy bien podr\'ian haber cien a\~nos de
sequ\'ia y un d\'ia de diluvio. Esto nos muestra que una desigualdad
tan general como la de Markov puede valer solo para la cola superior de la
distribuci\'on, no para la cola inferior.) 

La prueba es sencilla: por la definici\'on de la esperanza, tenemos
\[\mathbb{E}(X) \geq 0 \cdot \Prob(X<t) + t\cdot \Prob(X\geq t)
 = t \cdot \Prob(X \geq t),\]
y por lo tanto
\[\Prob(X\geq t) \leq  \frac{\mathbb{E}(X)}{t} .\]
\newline

{\em Aplicaciones.} Obtenemos de manera inmediata que
\begin{equation}\label{eq:gaudi}\begin{aligned}
\Prob(\tau(n)\geq t) \leq \frac{\log N + O(1)}{t},\\
\Prob(\omega(n)\geq t) \leq \frac{\log \log N + O(1)}{t}.
\end{aligned}\end{equation}
Podemos mejorar la segunda cota en (\ref{eq:gaudi})
de la manera siguiente.
Es f\'acil ver que $\tau(n)\geq 2^{\omega(n)}$. Luego
\begin{equation}\label{eq:ods}\Prob(\omega(n)\geq t) \leq \Prob(\tau(n) \geq 2^t) \leq \frac{\log N + O(1)}{2^t}
.\end{equation}
Qu\'e tan mejor es esto que la segunda l\'inea de (\ref{eq:gaudi})?
Consideremos $t = (1 + \epsilon) \log_2 \log N$. Entonces (\ref{eq:gaudi})
nos da $\Prob(\omega(n)\geq t) \leq \frac{\log 2 + o(1)}{1 + \epsilon}$, mientras
que (\ref{eq:ods}) nos da
\[\Prob(\omega(n)\geq t) \leq \frac{1 + o(1)}{(\log N)^{\epsilon}},\] lo cual es
una cota mucho m\'as fuerte (es decir, baja).

Por otra parte, si $t$ est\'a entre $\log \log N$ y $\log_2 \log N$,
la desigualdad (\ref{eq:ods}) no nos da nada. Esto se debe al hecho que,
si bien $\tau(n)=2^{\omega(n)}$ gran parte del tiempo,
$\mathbb{E}(\tau(n)) = \log N$, mientras que $\mathbb{E}(\omega(n)) = \log \log N$;
en otras palabras, $\mathbb{E}(\tau(n)) \ne 2^{\mathbb{E}(\omega(n))}$. Las colas
superiores de la distribuci\'on de $\omega(n)$ cobran gran efecto cuando
$\omega(n)$ se pone en el exponente, al punto que afectan considerablemente 
la esperanza de $\mathbb{E}(\tau(n))$ (o la de $\mathbb{E}(2^{\omega(n)})$, la cual
es del mismo orden de magnitud\footnote{
Es decir, es de tama\~no comparable, 
qu\'itese o p\'ongase un factor constante.}). 

Tendremos la oportunidad de estimar las distribuciones de $\omega(n)$ y
$\tau(n)$ con mayor precisi\'on m\'as tarde.
\\

\begin{center}
{\bf Notas y problemas}
\end{center}
\begin{enumerate}
\item\label{it:raul} {\em Sumas por partes.}\index{sumas por partes}
La siguiente t\'ecnica es \'util a menudo; la necesitaremos inmediatamente
en la prueba de Chebyshev-Mertens y una y otra vez en el futuro.
Digamos que tenemos que calcular
\[\sum_{n=1}^N h(n),\]
donde $h(n) = (f(n+1) - f(n))\cdot g(n)$.
Entonces
\[\begin{aligned}
\sum_{n=1}^N h(n)
&= \sum_{n=1}^N (f(n+1) - f(n))\cdot g(n) = 
\sum_{n=1}^N f(n+1) g(n) - \sum_{n=1}^N f(n) g(n)\\
&= \sum_{n=2}^{N+1} f(n) g(n-1) - \sum_{n=1}^N f(n) g(n)\\
&= (f(N+1) g(N) - f(1) g(1)) - \sum_{n=2}^N f(n) (g(n) - g(n-1)) .
\end{aligned}\]
Esta t\'ecnica ({\em sumar por partes}) 
 es \'util cuando la suma $\sum_{n=1}^N f(n) (g(n) - g(n-1))$ es, por
alg\'un motivo, m\'as facil de calcular que $\sum_{n=1}^N 
(f(n+1) - f(n)) \cdot g(n)$, o ya ha
sido evaluada.

Se puede ver que el proceso es an\'alogo a la integraci\'on por partes.
(Uno puede, incluso, ver a la sumaci\'on por partes como un caso especial
de la integraci\'on por partes, mediante el uso de una integral de Lebesgue.)

\item\label{it:chem} Probaremos el teorema de Chebyshev-Mertens\index{teorema!de Chebyshev-Mertens} 
(ecuaci\'on (\ref{eq:ort})).\index{n\'umero de primos hasta $N$}
\begin{enumerate}
\item
Todo n\'umero entero positivo puede ser
expresado como un producto de primos de manera 
\'unica\footnote{Este hecho es a veces llamado el {\em teorema
fundamental de la aritm\'etica}.\index{teorema!fundamental de la aritm\'etica}
Puede parecer extra\~no que un
enunciado tan familiar tenga un nombre tan grandilocuente; empero,
el hecho que un enunciado nos sea sumamente natural no quiere decir que
no deba ser probado, o que sea cierto. Hay an\'alogos del conjunto
de enteros $\mathbb{Z}$, los as\'i llamados {\em anillos de enteros}
de los {\em campos algebraicos}; en la gran mayor\'ia de 
ellos, el teorema fundamental de la aritm\'etica deja de ser cierto. 
(Si bien todo elemento a\'un se factoriza en elementos irreducibles, 
ya no lo hace de manera \'unica.) Tenemos, por ejemplo, las dos
factorizaciones $6 = 2\cdot 3 = (1 + \sqrt{-5}) (1 - \sqrt{-5})$
en el anillo de enteros del campo algebraico $\mathbb{Q}(\sqrt{-5})$.}. 
En otras palabras, para todo entero positivo $n$,
\[n = \prod_p p^{v_p(n)},\]
donde $v_p(n)$ es el m\'aximo entero no negativo $k$ tal que $p^k|n$.

Tome logaritmos a ambos lados y muestre que
\begin{equation}\label{eq:orgo}\log n = \sum_{d|n} \Lambda(d),\end{equation}
donde $\Lambda(d) = \log p$ si $d$ es una potencia $p^{\alpha}$ de un
\index{$\Lambda(n)$}
primo $p$, y $\Lambda(d) = 0$ si no es as\'i ({\em funci\'on de von
Mangoldt}).\index{funci\'on de von Mangoldt}

\item Sea $X_d$ como antes, es decir, la variable aleatoria que toma el
valor $1$ cuando $d|n$ y el valor $0$ cuando $d\nmid n$. Sea
$Y = \sum_{d|n} \Lambda(d) X_d$. Entonces, por (\ref{eq:orgo}), $Y$ siempre
toma el valor $\log n$. Concluya que
\begin{equation}\label{eq:asa}
\mathbb{E}(Y) = \log N + O(1).\end{equation}

\item Al mismo tiempo, tenemos que
\[\mathbb{E}(Y) = \sum_{d\leq N} \Lambda(d) \mathbb{E}(X_d)  =
\sum_{d\leq N} \Lambda(d) \cdot \frac{1}{N} \left\lfloor \frac{N}{d}\right\rfloor,\]
as\'i que
\begin{equation}\label{eq:chlam}\sum_{d\leq N} \Lambda(d) \cdot \frac{1}{N} \left\lfloor \frac{N}{d}\right\rfloor
= \log N + O(1).\end{equation}
Como $\frac{1}{N} \lfloor \frac{N}{d}\rfloor = \frac{1}{d} - O\left(
\frac{1}{N}\right)$,  estamos a un paso de obtener una estimaci\'on de
$\sum_{d\leq N} \frac{\Lambda(d)}{d}$:
\begin{equation}\label{eq:colto}\begin{aligned}
\sum_{d\leq N} \frac{\Lambda(d)}{d} 
&= \log N + O(1) + \sum_{d\leq N} 
\Lambda(d) \cdot O\left(\frac{1}{N}\right)\\
&= \log N + O(1) + \frac{1}{N} \cdot O\left(\sum_{d\leq N} 
\Lambda(d) \right).\end{aligned}\end{equation}
S\'olo nos falta acotar
$\sum_{d\leq N} \Lambda(d)$.

\item\label{it:istry}
Por (\ref{eq:chlam}), tenemos
\[\begin{aligned}
\sum_{d\leq N} \Lambda(d) \cdot \frac{1}{N} \left\lfloor \frac{N}{d}\right\rfloor
= \log N + O(1),\\
\sum_{d\leq N} \Lambda(d) \cdot \frac{1}{N/2} \left\lfloor \frac{N/2}{d}\right\rfloor
= \log \frac{N}{2} + O(1)
\end{aligned}\]
y por lo tanto
\[\begin{aligned}
\sum_{\frac{N}{2} \leq d \leq N} \Lambda(d) &\leq \sum_{d\leq N}
\Lambda(d) \cdot \left(\left\lfloor \frac{N}{d} \right\rfloor - 2 \left\lfloor 
\frac{N}{2 d}\right\rfloor\right)\\ &= N (\log N + O(1)) - 
N \left(\log \frac{N}{2} + O(1) \right) = O(N)
\end{aligned}\]
para todo $N$. Por lo tanto,
\begin{equation}\label{eq:shotwo}
\begin{aligned}
\sum_{d\leq N} \Lambda(d) &= \sum_{\frac{N}{2} <d \leq N} \Lambda(d) + 
\sum_{\frac{N}{4} < d \leq \frac{N}{2}} \Lambda(d) + 
\sum_{\frac{N}{8} < d \leq \frac{N}{4}} \Lambda(d) + \dotsc\\
&= O(N) + O(N/2) + O(N/4) + \dotsc = O(N).
\end{aligned}\end{equation}
(Aqu\'i lo que hemos hecho es dividir una suma en {\em intervalos di\'adicos},
es decir, intervalos de la forma $M < d \leq 2 M$; este es un procedimiento
muy com\'un en el an\'alisis.)\index{intervalos di\'adicos}

\item
De (\ref{eq:colto}) y (\ref{eq:shotwo}), deducimos que
\begin{equation}\label{eq:osto}
\sum_{d\leq N} \frac{\Lambda(d)}{d}  = \log N + O(1).
\end{equation}

\item Si bien (\ref{eq:osto}) ya es un resultado \'util, lo que queremos
en verdad es estimar la suma $\sum_{p\leq N} \frac{1}{p}$. Ahora bien,
la contribuci\'on de los enteros $d$ de la forma $d = p^{\alpha}$, $\alpha\geq
2$, a la suma (\ref{eq:osto}) es negligible, o, m\'as precisamente, $O(1)$.
(Por qu\'e? Porque $\sum_n \frac{\log n}{n^2}$ es convergente.) 
Tenemos entonces que
\begin{equation}\label{eq:sumpri}
\sum_{p\leq N} \frac{\log p}{p} = \log N + O(1).\end{equation}

Para liberarnos del factor $\log p$, podemos hacer una suma por partes
(ver la nota \ref{it:raul} m\'as arriba). Utilice tal t\'ecnica (c\'omo?) para
concluir que
\[\sum_{p\leq N} \frac{1}{p} = \sum_{2\leq n\leq N} 
\frac{\log n + O(1)}{n (\log n)^2} + O(1).\]
Aproximando la suma mediante una integral, muestre que
\[\sum_{2\leq n\leq N} \frac{\log n + O(1)}{n (\log n)^2} = \log \log N + O(1)\]
y por lo tanto
\begin{equation}\label{eq:cruj}\sum_{p\leq N} \frac{1}{p} = \log \log N + O(1)
\text{\;\;\;\;\;\;
({\em teorema de Chebyshev-Mertens}).}\end{equation}
\index{teorema!de Chebyshev-Mertens}
\item Denotemos por $\pi(N)$ el n\'umero de primos entre $1$ y $N$.
(Aqu\'i $\pi$ es la letra griega que corresponde a $p$, la cual es la primera
letra de la palabra {\em primo}; no hay otra connexi\'on con el n\'umero 
$\pi=3.14159\dotsc$.) 
Utilice la t\'ecnica de la suma por partes para deducir de (\ref{eq:shotwo})
que \begin{equation}\label{eq:kubisch}\pi(N) \ll \frac{N}{\log N}.
\end{equation}

\'Este es un resultado de Chebyshev.\index{Chebyshev, P.}
\item En verdad, siguiendo un procedimiento similar al que acabamos
de poner en pr\'actica, Chebyshev prob\'o resultados m\'as fuertes y m\'as 
precisos; en particular, mostr\'o que
\begin{equation}\label{eq:justu}
(\log 2) \frac{N}{\log N} (1+o(1)) \leq \pi(N) \leq
(\log 4) \frac{N}{\log N} (1+o(1))
.\end{equation}
N\'otese que Chebyshev dio
 una cota inferior, no s\'olo una cota inferior
como (\ref{eq:kubisch}). Probar
(\ref{eq:justu})
puede ser un problema interesante para el lector; alternativamente, se puede
consultar \cite[\S 2.2]{IK}, por ejemplo. Aqu\'i nos hemos querido
concentrar en derivar Chebyshev-Mertens\index{teorema!de Chebyshev-Mertens} (\ref{eq:cruj}) de la manera
m\'as breve posible.

 M\'as tarde, en 1896, Hadamard\index{Hadamard, J.} y 
de la Vall\'ee Poussin\index{de la Vall\'ee Poussin, C.} mostraron
 (independientemente el uno del otro) que
\begin{equation}\label{eq:dodoto}
\pi(N) \sim \frac{N}{\log N}\;\;\;\;\;\;\;
\text{(teorema de los n\'umeros primos)}\end{equation}\index{teorema!de los n\'umeros primos}
La mayor\'ia de las demostraciones
de (\ref{eq:dodoto}) requieren iniciar el estudio de la funci\'on
zeta de Riemann.\index{funci\'on zeta de Riemann} (Ver, e.g.,
\cite[\S 5.6]{IK}.) Existen tambi\'en pruebas ``elementales''\footnote{En el
sentido de no utilizar el an\'alisis complejo.}, generalmente complicadas.

Utilizaremos el teorema de los n\'umeros primos (\ref{eq:dodoto})
muy poco en estas notas, ya que Chebyshev-Mertens
nos ser\'a casi siempre suficiente. 
\end{enumerate}

%\item Las generalizaciones (en el sentido no matem'atico) 
%son peligrosas (en la matem\'atica y en otras partes). 
%Empero, podemos decir que, si gran parte de la teor\'ia probabil\'istica 
%de n\'umeros desciende de la desigualdad de Markov (aplicada
%a las variables $X$, $X^2$, $(X - \mathbb{E}(X))^2$ (Chebyshev), 
%$X^k$ y $e^X$, como veremos m\'as adelante), gran parte de la 
%combinatoria probabil\'istica depende de la siguiente aseveraci\'on, 
%a\'un m\'as 
%simple: para toda variable aleatoria $X$ y todo $\mathbb{R}$,
%\begin{equation}\label{eq:argo}
%\mathbb{E}(X) > t \Rightarrow \Prob(X>t)>0.\end{equation}

%En otras palabras, si la esperanza de la lluvia de ma\~nana m\'as es $3$ mm,
%se deduce que es posible que caigan m\'as de $3$ mm de lluvia.

%\begin{enumerate}
%\item 
%Muchas de las aplicaciones combinatoriales 
%de (\ref{eq:argo}) son aplicaciones del caso particular $t=0$. En otras
%palabras, si la esperanza de una variable es positiva, la variable... 
%junto con la observaci\'on que, si algo pasa con probabilidad positiva,
%a veces existe. (Tomar el ejemplo de Alon & Noga?)
%\item
%Lo de $Z/p$
%\end{enumerate}
%La aseveraci\'on (\ref{eq:argo}) es llamada por algunos {\em principio del
%palomar}, {\em principio de las casillas}\footnote{La traducci\'on m\'as
%idiom\'atica de ``pigeonhole principle''.} o {\em principio de Dirichlet}.
\end{enumerate}
\section{La varianza}\label{sec:varianza}\index{varianza}
%{\bf Varianza y sumas.} 
La {\em varianza} de una variable aleatoria est\'a dada por
\begin{equation}\label{eq:og}\Var(X) := \mathbb{E}((X-\mathbb{E}(X))^2) = 
\mathbb{E}(X^2) - \mathbb{E}(X)^2 .\end{equation}
Sean $X$, $Y$ dos variables independientes. Entonces
\begin{equation}\label{eq:agasto}
\mathbb{E}(XY) = \mathbb{E}(X) \mathbb{E}(Y),\end{equation}
y luego
\begin{equation}\label{eq:agar}
\Var(X + Y) = \mathbb{E}((X+Y)^2) - \mathbb{E}(X+Y)^2 = 
\Var(X) + \Var(Y) .\end{equation}
En general, si $X_1, X_2, \dotsc, X_n$ son variables independientes en pares
(es decir, si $X_i$, $X_j$ son independientes para $i,j\in \{1,2,\dotsc,n\}$
distintos cualesquiera),
\begin{equation}\label{eq:marcost}
\Var(X_1 + X_2 + \dotsc + X_n) = \Var(X_1) + \dotsc + \Var(X_n) .\end{equation}

%{\bf C\'omo acotar colas con la varianza.}
\begin{thm}{(Desigualdad de Chebyshev)}\label{thm:herper}\index{desigualdad!de Chebyshev}
Para toda variable aleatoria $X$ y todo $x>0$,
\[\Prob(|X - \mathbb{E}(X)|\geq x) \leq
\frac{\Var(X)}{x^2}.\]
En particular, si $X = X_1 + X_2 + \dotsc + X_n$, donde
$X_1, X_2, \dotsc, X_n$ son
variables independientes en pares,
\begin{equation}\label{eq:kavort}
\Prob(|X - \mathbb{E}(X)|\geq x) \leq
\frac{\Var(X)}{x^2} = \frac{\Var(X_1) + \dotsc + \Var(X_n)}{x^2} .\end{equation}
\end{thm}
\begin{proof}
Utilizamos la definici\'on (\ref{eq:og}) y la desigualdad de Markov\index{desigualdad!de Markov}
(\ref{eq:mark}):
\[\Prob(|X - \mathbb{E}(X)|\geq x) =
\Prob(|X - \mathbb{E}(X)|^2\geq x^2) 
\leq \frac{\mathbb{E}(|X - \mathbb{E}(X)|^2)}{x^2}
= \frac{\Var(X)}{x^2} .\]
Si $X = X_1 + \dotsc + X_n$, utilizamos (\ref{eq:marcost}) para evaluar
$\Var(X)$.
\end{proof}

{\em Aplicaci\'on.} Sabemos ya que, en promedio, un n\'umero $n\leq N$ tiene
$\log \log N + O(1)$ factores primos. Qu\'e tan comunes son los n\'umeros
que tienen muchos menos o muchos m\'as factores primos?

Sea $X_p$ como en (\ref{eq:maz}). 
 Ahora bien, $X_p^2 = X_p$, puesto que $0^2 = 0$ y
$1^2 = 1$. Por lo tanto, 
\[\Var(X_p) = \mathbb{E}(X_p^2) - \mathbb{E}(X_p)^2 =
\mathbb{E}(X_p) - \mathbb{E}(X_p)^2 =  
\mathbb{E}(X_p) - O(1/p^2). \]
La ecuaci\'on (\ref{eq:ost})
nos dice que $\sum_{p\leq N} \mathbb{E}(X_p) = \log \log N + O(1)$.
 As\'i,
\[\sum_{p\leq N} \Var(X_p) = \sum_{p\leq N} \mathbb{E}(X_p) - 
\sum_{p\leq N} \frac{1}{p^2} =
\sum_{p\leq N} \mathbb{E}(X_p) - O(1)
= \log \log N + O(1) .\]

Podemos concluir que, para $\omega(n) = X = \sum_{p\leq N} X_p$, 
la desigualdad de Chebyshev es v\'alida con\index{desigualdad!de Chebyshev}
$\mathbb{E}(X) = \log \log N + O(1)$? No nos apresuremos: las
variables $X_{p_1}$, $X_{p_2}$, $p_1\ne p_2$ no son exactamente
independientes. 
(Por ejemplo, si $p_1,p_2\geq \sqrt{N}$, entonces
$X_{p_1}$ y $X_{p_2}$ no pueden ser $1$ simult\'aneamente (por qu\'e?);
esto nunca pasar\'ia con variables verdaderamente independientes.)
Nos basta, empero, que la igualdad (\ref{eq:agasto}) sea v\'alida de manera
aproximada. Veamos: para $p_1\ne p_2$,
\[\begin{aligned}
\mathbb{E}(X_{p_1} X_{p_2}) &= \mathbb{E}(X_{p_1 p_2}) = \frac{1}{p_1 p_2}
+ O\left(\frac{1}{N}\right),\\
\mathbb{E}(X_{p_1}) \mathbb{E}(X_{p_2}) &= 
(1/p_1 + O(1/N)) (1/p_2 + O(1/N)) = \frac{1}{p_1 p_2} + 
O\left(\frac{1}{N}\right).\end{aligned}\]
Por lo tanto,
\[\Var(X_{p_1} + X_{p_2}) = \mathbb{E}((X_{p_1} + X_{p_2})^2) -
\mathbb{E}(X_{p_1} + X_{p_2})^2 = \Var(X_{p_1}) + \Var(X_{p_2}) + O(1/N),\]
y, de la misma manera,
\[\Var(\sum_{p\leq M} X_p) = \sum_{p\leq M} \Var(X_p) + O(M^2/N)\]
para todo $M$. Ahora bien, podemos escoger $M$ de tal manera que
el t\'ermino de error $O(M^2/N)$ sea peque\~no -- digamos, $M = N^{1/3}$.
Concluimos, por la desigualdad de Chebyshev, que
\begin{equation}\label{eq:sibkul}
\Prob(|X - \mathbb{E}(X)|\geq x) \leq \frac{\log \log N + O(1)}{x^2} 
\end{equation}
para $X = \sum_{p\leq N^{1/3}} X_p$.

Ahora bien, cu\'al es la diferencia entre $X$ y $\omega(n)$? Un n\'umero
$n\leq N$ no puede tener m\'as de dos divisores primos $>N^{1/3}$: m\'as
no caben. Por lo tanto, $|X - \omega(n)|$ nunca es m\'as de $2$.
Obtenemos que
\begin{equation}\label{eq:ogo}
\Prob(|\omega(n) - \log \log n|\geq x) \leq \frac{\log \log N + O(1)}{
(x+O(1))^2}
.\end{equation}
Dicho de otra manera,
\begin{equation}\label{eq:gortos}\Prob\left(|\omega(n) - \log \log N| \geq t \sqrt{\log \log N}\right) \leq
\frac{1 + O(1/\sqrt{\log \log N})}{t^2} \end{equation}
para todo $t\geq 1$. (La constante impl\'icita en $O(1/\sqrt{\log \log N})$
no depende de $t$.) Tanto el resultado (\ref{eq:gortos}) como la prueba
que hemos presentado se deben a Tur\'an \cite{T}; Hardy y Ramanujan hab\'ian
dado antes una prueba m\'as complicada de un resultado
ligeramente m\'as d\'ebil \cite{HR}.

{\em Ejemplos.} Escogemos $t = 10$, y obtenemos que
\[\Prob(|\omega(n) - \log \log N| \geq 10 \sqrt{\log \log N}) \leq
\frac{1}{100} + o(1);\]
escogemos $t = \epsilon \sqrt{\log \log N}$, y obtenemos que
\begin{equation}\label{eq:amalia}\begin{aligned}
\Prob(\omega(n) > (1 + \epsilon) \log \log N) &\leq \frac{1}{\epsilon^2
\log \log N} + o_{\epsilon}(1),\\
\Prob(\omega(n) < (1 - \epsilon) \log \log N) &\leq \frac{1}{\epsilon^2
\log \log N} + o_{\epsilon}(1)\end{aligned}\end{equation}
para todo $\epsilon>0$.

\begin{center}
{\bf Notas y problemas}
\end{center}
\begin{enumerate}\index{n\'umero de primos hasta $N$}
\item Nuevamente nos planteamos la pregunta: cu\'antos primos hay entre 1 y $N$? Tratemos de ver si podemos atacar el
  problema usando la varianza, antes que la esperanza. La idea central est\'a
clara: los primos son algo que se desv\'ian de la norma, y podemos usar
la desigualdad de Chebyshev (Teorema \ref{thm:herper}) \index{desigualdad!de Chebyshev}para obtener una cota
sobre la probabilidad de eventos que se desv\'ian de una norma. 
Podemos usar la desigualdad
(\ref{eq:ogo}) (la cual hemos probado utilizando la desigualdad de
Chebyshev)
 con $x = \log \log n$, y obtenemos
\[\Prob(\omega(n) = 1) \leq \frac{\log \log N + O(1)}{(\log \log N)^2 + O(\log \log
  N)} = \frac{1}{\log \log N + O(1)} .\]
Por lo tanto, hay a lo m\'as $\frac{N}{\log \log N + O(1)}$ primos entre $1$
y $N$. Esta es una cota sumamente d\'ebil:
la cota (\ref{eq:kubisch}) era mucho mejor.
Podemos utilizar la
desigualdad de Chebyshev de otra manera para obtener una cota m\'as fuerte?

Veremos que es as\'i, y luego veremos 
que la gran ventaja de lo que haremos sobre lo que hicimos en las
notas en la secci\'on anterior es que el que m\'etodo que seguiremos
ahora tambien sirve para obtener cotas sobre muchas cosas aparte del n\'umero
de primos entre $1$ y $N$.

Lo que estamos haciendo es evaluar la varianza de $X = \sum_p X_p$. Al
hacer tal cosa, utilizamos el hecho que las variables $X_p$, $p\leq N^{1/3}$
(digamos) son casi independientes en pares. Hay todav\'ia
un hecho m\'as general
 que no estamos usando: para $p_1<p_2<\dotsc<p_k$ cualesquiera
tales que $p_1 p_2\dotsb p_k$ es bastante menor que $N$, las variables
$X_{p_1}, X_{p_2},\dotsc , X_{p_k}$ son mutuamente independientes, o casi.
Qu\'e podemos hacer con esto?
\begin{enumerate}
\item Defina
\begin{equation}\label{eq:tortor}Z_p = \begin{cases} 
1/p &\text{si $p\nmid n$,}\\
-(1-1/p) & \text{si $p|n$,} \end{cases}\end{equation}
donde $n$ es un entero aleatorio entre $1$ y $N$. Verifique que 
$\mathbb{E}(Z_p) = O(1/N)$.
\item Para todo $d$ sin divisores cuadrados\footnote{Es decir,
$d$ no divisible por $4$, ni por $9$, ni por $16$, ni por $25$, \dots}, defina
\begin{equation}\label{eq:tartar}Z_d = \prod_{p|d} Z_p .\end{equation}
Verifique que $\mathbb{E}(Z_d) = O(\tau(d)/N)$.
Ya sabemos que $\tau(d)$ es peque\~no en promedio (ver (\ref{eq:gotz})).

Concluya que, si $d_1$, $d_2$ son distintos y carecen de divisores cuadrados,
\begin{equation}\label{eq:fornost}\begin{aligned}
\mathbb{E}(Z_{d_1}) \mathbb{E}(Z_{d_2})
 &= N^{-2} \cdot O(\tau(d_1)) O(\tau(d_2)),\\
\mathbb{E}(Z_{d_1} Z_{d_2}) &= N^{-1} \cdot O(\tau(d_1 d_2))
\leq N^{-1} \cdot O(\tau(d_1)) O(\tau(d_2)),\end{aligned}\end{equation}
y por lo tanto
\begin{equation}\label{eq:lamerm}
\mathbb{E}(Z_{d_1} Z_{d_2}) = \mathbb{E}(Z_{d_1}) \mathbb{E}(Z_{d_2}) + 
N^{-1} \cdot O(\tau(d_1)) O(\tau(d_2)) .\end{equation}
\item Defina $Z = \sum_{d\leq M}^* Z_d$, donde $M = N^{0.49}$.
(El asterisco $*$ en la suma $\sum_{d\leq M}^*$ quiere decir
que $d$ recorre
s\'olo a los enteros sin divisores cuadrados.)

Podemos ver que (\ref{eq:lamerm}) es una versi\'on aproximada de
(\ref{eq:agasto}); es razonable tratar de obtener una versi\'on
aproximada de (\ref{eq:marcost}) en consecuencia.
 Muestre que 
\begin{equation}\label{eq:scott} \begin{aligned}
\mathbb{E}(Z) &= O(N^{-1/2}),\\
\Var(Z) &= \sideset{}{^*}\sum_{d\leq M} \Var(Z_d) + O(N^{-0.01}) = 
\sideset{}{^*}\sum_{d\leq M} \mathbb{E}(Z_d^2) + O(N^{-0.01}) 
\end{aligned}\end{equation}
Muestre tambi\'en que $\mathbb{E}(Z_d^2) = \frac{\phi(d)}{d^2} + 
O\left(\frac{\tau(d)}{N}\right)$, donde $\phi(d) = d \cdot
\prod_{p|d} (1 - 1/p)$ ({\em funci\'on de Euler}).
Por consiguiente,
\begin{equation}\label{eq:anac}
\Var(Z) = \sideset{}{^*}\sum_{d\leq M} \frac{\phi(d)}{d^2} +
O\left(N^{-0.01}\right) \ll \log M
\leq \log N.\end{equation}

\item Por la desigualdad de Chebyshev,
\begin{equation}\label{eq:oates}
\Prob(|Z - \mathbb{E}(Z)| \geq x) \leq \frac{\Var(Z)}{x^2} .\end{equation}
Ahora bien, si el n\'umero $n$ es primo y mayor que $M$, entonces, para cada
$d$ sin divisores cuadrados, la variable
$Z_d$
tomar\'a el valor $\prod_{p|d} 1/p = 1/d$
(por la definici\'on de $Z_d$; ver (\ref{eq:tortor}) y (\ref{eq:tartar})).
Por lo tanto, si $n$ es primo y mayor que $M$, 
\[Z = \sideset{}{^*}\sum_{d\leq M} \frac{1}{d} \gg 
\frac{\sum_{d\leq M} \frac{1}{d}}{
\sum_m \frac{1}{m^2}} \gg
 \sum_{d\leq M} \frac{1}{d} \gg \log M 
\gg \log N,\]
donde utilizamos el hecho que la suma $\sum_m \frac{1}{m^2}$ converge.

Se infiere inmediatamente que
\[\Prob(\text{$n$ es primo y mayor que $M$}) 
\leq \Prob(|Z - \mathbb{E}(Z)|\geq x)\]
con $x = \sum_{d\leq M}^* 1/d - \mathbb{E}(Z) \gg \log N - O(N^{-1/2}) \gg
\log N$.
Por (\ref{eq:anac}) y (\ref{eq:oates}), concluimos que
\[
\Prob(\text{$n$ es primo y mayor que $M$}) \ll \frac{1}{\log N},\]
y por lo tanto
\begin{equation}\begin{aligned}
\Prob(\text{$n$ es primo}) &\ll \frac{1}{\log N} + \Prob(n\leq M)\\
&=\frac{1}{\log N} + \frac{M}{N} \ll \frac{1}{\log N}
\end{aligned}\end{equation}
para $n$ tomado al azar entre $1$ y $N$. En otras palabras, el n\'umero
de primos entre $1$ y $N$ es $\ll \frac{N}{\log N}$.
\end{enumerate}

\'Esta es esencialmente la misma cota que ya obtuvimos en
 \S \ref{sec:esperanza}, Problema \ref{it:istry}.  La ventaja del
m\'etodo presente reside en su suma flexibilidad: v\'ease el problema 
siguiente.

\item

Procediendo como lo hicimos en el problema anterior,
probaremos que
\begin{equation}\label{eq:tatan}
\Prob(\text{tanto $n$ como $n+2$ son primos}) \ll \frac{1}{(\log N)^2}
\end{equation}
para $n$ tomado al azar entre $1$ y $N$.

\begin{enumerate}
\item Comenzamos definiendo
\begin{equation}\label{eq:lauror}Z_p = \begin{cases} 
2/p &\text{si $p\nmid n$ y $p\nmid n+2$,}\\
-(1-2/p) & \text{si $p|n$ o $p|n+2$,} \end{cases}\end{equation}
donde $n$ es un entero aleatorio entre $1$ y $N$, y
\[Z_d = \prod_{p|d} Z_p.\]
De la misma manera que antes, se puede ver que
$\mathbb{E}(Z_{d_1} Z_{d_2})$ y $\mathbb{E}(Z_{d_1}) \mathbb{E}(Z_{d_2})$
son sumamente peque\~nos para $d_1, d_2$
distintos y sin divisores cuadrados.

\item Podr\'iamos definir, como antes, $Z = \sum_{d\leq M}^* Z_d$,
donde $M = N^{\frac{1}{2} - \epsilon}$. (El asterisco $*$ en la
suma $\sum_{d\leq M}$ denota que la suma recorre solo a los $d$ sin
 divisores cuadrados.) Esto podr\'ia dar resultados.
 Empero, tenemos el derecho de definir
$Z = \sum_{d\leq M}^* c_d Z_d$, para $c_d$'s arbitrarios; hacemos esto,
y afrontamos la tarea de encontrar los $c_d$ que nos den el
mejor resultado.

(Esta tarea de optimizaci\'on nos dar\'a una mejora cuantitativa, antes
que cualitativa; podr\'iamos obtener (\ref{eq:tatan}) sin este paso.
Por suerte, ciertos c\'alculos finales nos ser\'an m\'as simples
de esta manera que si escogieramos $c_d = 1$.)

Como antes, la idea es usar
\[\Prob(|Z|\geq x) \leq \frac{\mathbb{E}(Z^2)}{x^2},\]
donde $x$ es el valor que $Z$ toma cuando $n$ y $n+2$ son ambos primos.

Muestre que, cuando $n$ y $n+2$ son ambos primos,
\[Z = \sideset{}{^*}\sum_{d\leq M} 
c_d \frac{\tau(d)}{d}.\]
Muestre tambi\'en que
\begin{equation}\label{eq:stoste}\begin{aligned}
\mathbb{E}(Z^2) &= \sideset{}{^*}\sum_{d\leq M} 
 c_d^2 \mathbb{E}(Z_d^2) + 
O\left(N^{-1} \left(\sideset{}{^*}\sum_{d\leq M} 
|c_d| \tau(d)^3 \right)^2\right)\\
&= \sideset{}{^*}\sum_{d\leq M} 
 c_d^2 \frac{\tau(d)}{d} \prod_{p|d} \left(1 - \frac{2}{p}
\right) 
+ O\left(N^{-1} \left(\sideset{}{^*}\sum_{d\leq M} 
|c_d| \tau(d)^3 \right)^2\right).
\end{aligned}\end{equation}

\item Debemos, entonces, encontrar el m\'inimo de
\begin{equation}\label{eq:urpi}
\frac{\sideset{}{^*}{\sum_{d\leq M}} 
 c_d^2 \frac{\tau(d)}{d} \prod_{p|d} 
\left(1 - \frac{2}{p}\right)}{
\left(\sideset{}{^*}{\sum_{d\leq M}} 
 c_d \frac{\tau(d)}{d}\right)^2}\end{equation}

La preguntas es: como escogemos $c_d$ de tal manera que
(\ref{eq:urpi}) sea m\'inimo? O, m\'as bien: cu\'al es 
el m\'inimo valor tomado por (\ref{eq:urpi})?

Para $a_d$, $b_d$ cualesquiera,
\begin{equation}\label{eq:cutim}
\left(\sum_d a_d b_d\right)^2 \leq \sum_d a_d^2 \cdot \sum_d b_d^2
\text{\;\;\;\;\; (desigualdad de Cauchy)}
\end{equation}
con igualdad s\'olo cuando $\vec{a}$ y $\vec{b}$ son proporcionales,
i.e., cuando hay alg\'un $r$ tal que 
$a_d = r b_d$ para todo $d$ (o $a_d=0$ para todo $d$). (Prueba
de la desigualdad de Cauchy:
tenemos $\sum_{d<d'} (a_d b_d' - a_{d'} b_d)^2 \geq 0$ 
con igualdad
s\'i y s\'olo s\'i $a_d b_d' - a_{d'} b_d = 0$ para todo par $d$, $d'$,
lo cual a su vez ocurre si y s\'olo si $\vec{a}$ y $\vec{b}$ son proporcionales.
Expanda $\sum_{d<d'} (a_d b_d' - a_{d'} b_d)^2$, pase todos los t\'erminos
negativos 
al lado
derecho de la desigualdad
$\sum_{d<d'} (a_d b_d' - a_{d'} b_d)^2 \geq 0$, 
y sume $\sum a_d^2 b_d^2$ a cada lado.) 

La
desigualdad de Cauchy no es sino la familiar afirmaci\'on que el producto
de dos vectores es menor o igual que el producto de sus normas. 
En verdad, no necesitaremos la desigualdad de Cauchy, sino simplemente
el hecho (evidente) que (\ref{eq:cutim}) se vuelve una igualdad cuando
$a_d = b_d$ para todo $d$. La desigualdad de Cauchy s\'olo cumple el rol de
asegurarnos que estamos procediendo de la mejor manera posible (en este paso).

La expresi\'on (\ref{eq:urpi}) es igual a $\frac{\sum_{d\leq M}^* a_d^2}{
\left(\sum_{d\leq M}^* a_d b_d\right)^2}$ 
con \[a_d = c_d \sqrt{\frac{\tau(d)}{d}
\prod_{p|d} \left(1 - \frac{2}{p}\right)}
,\;\;\;\; b_d = 
\sqrt{\frac{\tau(d)}{d}} \left(\prod_{p|d} \left(1 - \frac{2}{p}\right)
\right)^{-1/2} .\]
Por la desigualdad de Cauchy, el m\'inimo de
$\frac{\sum_{d\leq M}^* a_d^2}{
\left(\sum_{d\leq M}^* a_d b_d\right)^2}$ es $\frac{1}{\sum_{d\leq M}^* b_d^2}$,
es decir,
\[\frac{1}{\sum_{d\leq M}^* \frac{\tau(d)}{d} \prod_{p|d} 
   \left(1 - \frac{2}{p}\right)^{-1}}.
\]
Este m\'inimo es alcanzado cuando $a_d = b_d$, i.e., cuando
$c_d = \prod_{p|d} \left(1 - \frac{2}{p}\right)^{-1}$.

Tenemos, entonces  -- utilizando 
(\ref{eq:stoste}) --
\[\begin{aligned}
\Prob(\text{$n$ y $n+2$ son primos}) &\leq \frac{\mathbb{E}(Z^2)}{x^2}\\
&\leq \frac{1}{\sum_d^* \frac{\tau(d)}{d} \prod_{p|d} 
   \left(1 - \frac{2}{p}\right)^{-1}} +  
O\left(\frac{M (\log N)^A}{N}\right) .
\end{aligned}\]
El t\'ermino $O\left(\frac{M (\log N)^A}{N}\right)$ es negligible
($\ll N^{-1/2}$). Ahora bien,
\[\begin{aligned}
\sideset{}{^*}\sum_{d\leq M} \frac{\tau(d)}{d} \prod_{p|d} \left(1 -
  \frac{2}{p}\right)^{-1}
&= \sideset{}{^*}\sum_{d\leq M} \frac{\tau(d)}{d} \prod_{p|d}
\left(1 + \frac{2}{p} + \frac{2^2}{p^2} + \frac{2^3}{p^3} + \dotsc \right)\\
&\geq \sum_{d\leq M} \frac{\tau(d)}{d}
\geq \left(\sum_{d\leq M^{1/2}} \frac{1}{d}\right)^2 \\ &\sim (\log M^{1/2})^2
\gg (\log N)^2 . 
\end{aligned}\]
Por lo tanto,
\begin{equation}
\Prob(\text{$n$ y $n+2$ son primos}) \ll \frac{1}{(\log N)^2} .
\end{equation}
\end{enumerate}
Lo que acabamos de hacer puede verse como una versi\'on del
m\'etodo llamado {\em criba de Selberg} (1950).\index{criba!de Selberg}

La primera prueba del resultado (\ref{eq:tatan}) fue dada por V. Brun (1920).

\item\label{it:gotor} Se sabe que $\Prob(\text{$n$ es primo}) \sim \frac{1}{\log N}$;
tal aseveraci\'on no es sino el teorema de los n\'umeros primos
(Hadamard -- de la Vall\'ee-Poussin, 1896),\index{teorema!de los n\'umeros primos} para el cual no se conoce
una demostraci\'on simple. Se cree que
\begin{equation}\label{eq:hl}
\Prob(\text{tanto $n$ como $n+2$ son primos}) 
\sim \frac{c_2}{(\log N)^2}\end{equation}
donde 
\[c_2 = 2 \prod_{p\geq 3} \frac{p (p-2)}{(p-1)^2} \sim 1.32032\dotsc\]
Empero, esta conjetura sigue sin
probarse; no se sabe siquiera si es que hay un n\'umero infinito
de primos $n$ tales que $n+2$ sea tambi\'en primo ({\em conjetura
de los primos gemelos})\index{conjetura de los primos gemelos}. (El enunciado (\ref{eq:hl}) es parte
de la {\em conjetura de Hardy-Littlewood},\index{conjetura de Hardy-Littlewood}
 la cual tambi\'en especifica, por
ejemplo, cu\'al debe ser la probabilidad que $n$, $n+2$ y $n+6$ sean todos
primos.)
%sieves
%pnt, t\omegain prime
%no hope of using sieves to get that; upper bounds, lower bounds -- but
%generally no asymptotics, and quite hard to get lower bounds (thus
%giving lower bounds for the tail of a distribution -- but not just that --
% also, technical issues (translation))

\item\label{it:gorgut}\index{problema de la tabla de multiplicaci\'on}
 Veamos ahora una bonita aplicaci\'on de los resultados
de esta secci\'on; tanto el resultado como la prueba se deben a 
Erd\H{o}s \cite{E}.\index{Erd\H{o}s, P.} El resultado es el siguiente: s\'olo una proporci\'on
$o(1)$ de los enteros $\leq N^2$ pueden expresarse como un producto $a\cdot b$
con $a, b\leq N$. (Cuando decimos una proporci\'on $o(1)$ (o, coloquialmente,
``proporci\'on $0$'') de los enteros $\leq N^2$, queremos decir ``$o(N^2)$
enteros $\leq N^2$''.) Probemos este resultado.
\begin{enumerate}
\item Sea $\epsilon$ un n\'umero peque\~no (digamos $\epsilon = 1/10$).
Entonces, por (\ref{eq:amalia}), 
\[\begin{aligned}
\Prob(\omega(a) < (1 - \epsilon) \log \log N) &= O\left(\frac{1}{\log \log
    N}\right) = o_{\epsilon}(1),\\ 
\Prob(\omega(b) < (1 - \epsilon) \log \log N) &= O\left(\frac{1}{\log \log
    N}\right) = o_{\epsilon}(1),\end{aligned}\]
donde $a$ y $b$ son enteros entre $1$ y $N$ tomados al azar.
\item Muestre que $\mathbb{E}(\omega(\mcd(a,b)))$ es $O(1)$.
\item Concluya que
\[\Prob(\omega(a \cdot b) < (2 - 3 \epsilon) \log \log N) = 
o_{\epsilon}(1) .\]
\item Nuevamente por (\ref{eq:amalia}),
\[\Prob(\omega(n) > (1 + \epsilon) \log \log N^2) \leq
o_{\epsilon}(1)\]
para $n$ tomado al azar entre $1$ y $N$. Tenemos, entonces, que
hay a lo m\'as $o(N^2)$ pares de enteros $a,b\leq N$ tales que
$\omega(a \cdot b) < (2 - 3 \epsilon) \log \log N$, y a lo m\'as
$o(N^2)$ enteros $n\leq N^2$ tal que
$\omega(n) > (1 + \epsilon) \log \log N^2$. Ahora bien, si
$n = a \cdot b$ (y $\epsilon$ es peque\~no y $N$ es grande),
 se debe dar o lo uno o lo otro (puesto que
$\log \log N^2 = \log \log N + \log 2$). Concluimos que
hay a lo m\'as $o(N^2)$ enteros $n$ entre $1$ y $N^2$ tales que
\[n = a\cdot b\]
para alg\'un par de enteros $1\leq a,b\leq N$. Esto es lo que quer\'iamos
demostrar.
\end{enumerate}
Examinaremos este problema en m\'as detalle cuando sepamos c\'omo obtener un
t\'ermino de error y darle un significado.

\item {\em Los n\'umeros desmenuzables. Momentos variables.}\index{n\'umeros desmenuzables} \index{momentos!momentos variables}
Se dice que un n\'umero es {\em desmenuzable}
si s\'olo tiene factores primos peque\~nos. 
Cu\'an comunes son los
n\'umeros desmenuzables? Esto es: cu\'antos n\'umeros $n$ de tama\~no $N$
tienen s\'olo factores primos $\leq N^{1/u}$, donde $u\to \infty$ cuando
$N\to \infty$? (Por ejemplo: cu\'antos n\'umeros de tama\~no $N$
tienen s\'olo factores primos
$\leq N^{\frac{1}{\log \log N}}$?)

\begin{enumerate}
\item Consideremos $N < n\leq 2 N$, $z = N^{1/u}$ y 
$P(z) = \prod_{p\leq z} p$. Para simplificar las cosas, comenzaremos
considerando
s\'olo enteros $n$ sin factores cuadrados; en otras palabras, nuestra meta
inicial es acotar la probabilidad que un n\'umero $n$ tomado al azar
entre $N$ y $2 N$ no tenga factores cuadrados ni factores primos $>N^{1/u}$.

Diremos que $n$ es {\em $u$-desmenuzable} si no tiene factores primos
$>N^{1/u}$. 
Si $n$ no tiene factores cuadrados, $n$ ser\'a $u$-desmenuzable 
si y s\'olo si
\begin{equation}\label{eq:jarusz}\gcd(n,P(z)) > N .\end{equation}
Como $P(z) = \prod_{p\leq z} p$ y $n = \prod_{p|n} p$ son productos, y como 
preferimos trabajar con sumas, sacamos logaritmos en (\ref{eq:jarusz}):
\[\log \gcd(n,P(z)) > \log N.\]
Ahora bien, $\log \gcd(n,P(z)) = \sum_{p\leq z} (\log p) \cdot X_p$, donde
$X_p = \begin{cases} 1 &\text{si $p|n$,}\\ 0 &\text{si $p\nmid n$.}
 \end{cases}$ Por lo tanto, nuestra tarea es acotar
\begin{equation}\label{eq:totoro}\Prob(X> \log N),\end{equation}
donde $X = \sum_{p\leq z} (\log p) \cdot X_p$ y $n$ es tomado
al azar entre $N$ y $2 N$.
\item Podemos acotar $\Prob(X> \log N)$ mediante Markov, o mediante Chebyshev, que no es sino Markov aplicado
a $X^2$, o a $(X - \mathbb{E}(X))^2$; de la misma manera, podemos acotar
$\Prob(X>\log N)$ mediante la desigualdad de Markov aplicada a $X^k$, para
un $k$ positivo de nuestra elecci\'on.

Para cualquier variable aleatoria $X$ y cualquier n\'umero par $k>0$ (o
para cualquier n\'umero $k>0$ y 
cualquier variable aleatoria $X$ que tome s\'olo valores no negativos),
\begin{equation}\label{eq:aet}
\Prob(X>t) \leq \frac{\mathbb{E}(X^k)}{t^k} .
\end{equation}
Esto no es sino Markov aplicado a $X^k$. Las desigualdades de Markov
y Chebyshev son los casos $k=1$ y $k=2$ de (\ref{eq:aet}). A la
utilizaci\'on de (\ref{eq:aet}) para $k$ general
 se le llama {\em acotaci\'on por momentos}. (La expresi\'on 
$\mathbb{E}(X^k)$ es llamada {\em un momento}.)\index{momentos}

\item Usaremos (\ref{eq:aet}) para estimar (\ref{eq:totoro}). Escogeremos
$k$ al final; no ser\'a una constante, sino una funci\'on de $N$.
Veamos:
\[\Prob(X> \log N) \leq \frac{\mathbb{E}(X^k)}{(\log N)^k},
\]
donde $X = \sum_{p\leq z} (\log p)\cdot X_p$. Proseguimos:
\[\begin{aligned}
\mathbb{E}(X^k) &= \mathbb{E}\left(\sum_{p_1,\dotsc,p_k\leq z}
 (\log p_1) \dotsb (\log p_k) \cdot X_{p_1} \dotsb X_{p_k}\right)\\
 &= \sum_{p_1,\dotsc,p_k\leq z}
 (\log p_1) \dotsb (\log p_k) \cdot \mathbb{E}\left(
X_{p_1} \dotsb X_{p_k}\right)\\
&= \sum_{p_1,\dotsc,p_k\leq z}
 (\log p_1) \dotsb (\log p_k) \cdot \mathbb{E}\left(
X_{p_1'} \dotsb X_{p_l'}\right),\end{aligned}\]
donde $p_1'<\dotsc<p_l'$ son los primos distintos entre $p_1,\dotsc,p_k$.
(Por ejemplo, si $k=4$ y $p_1 = 3$, $p_2 = 2$, $p_3 = 7$, $p_4 = 2$, entonces
$l=3$ y $p_1' = 2$, $p_2' = 3$, $p_3' = 7$.) Sabemos que
\begin{equation}\label{eq:hariman}\mathbb{E}(X_{p_1}' \dotsc X_{p_l}') = \frac{1}{N}
\left(\left\lfloor \frac{2 N}{m} \right\rfloor 
- \left\lfloor \frac{N}{m} \right\rfloor \right) 
\leq \frac{1}{N} \left\lfloor \frac{2 N}{m} \right\rfloor \leq \frac{2}{m},
\end{equation}
donde $m = p_1' p_2' \dotsc p_l'$. (Esta puede parecer una cota muy mala,
pero, en este problema,
 es mejor para nuestra salud que $\leq \frac{1}{m} + \frac{1}{N}$.)

Concluimos que
\begin{equation}\label{eq:thsain}
\mathbb{E}(X^k) = \sum_{p_1,\dotsc,p_k\leq z} 
\frac{2 (\log p_1) (\log p_2) \dotsb (\log p_k)}{p_1' p_2' \dotsb p_l'}
\end{equation}

\item Para estimar la suma (\ref{eq:thsain}), tenemos que estimar cuantas
veces los primos distintos $p_1' < p_2' < \dotsb < p_l'$ aparecen disfrazados de
$p_1, p_2, \dotsc, p_k$. Hay $l^k$ maneras de colorear $k$ objetos con $l$
colores. As\'i,
\[\mathbb{E}(X^k) \leq 2 \sum_{l=1}^k l^k \sum_{p_1' < \dotsb < p_l' \leq z}
 (\log z)^{k-l} \frac{(\log p_1') \dotsb (\log p_l')}{p_1' \dotsb p_l'} .\]
Ahora bien,
\[\sum_{p_1'<\dotsc<p_l'\leq z} \frac{(\log p_1') \dotsb (\log p_l')}{p_1'
  \dotsb p_l'} \leq \frac{1}{l!} \left(\sum_{p\leq z} \frac{\log p}{p}\right)^l
\]

Sabemos que $\sum_{p\leq z} \frac{\log p}{p} = \log z + O(1)$ 
(ver (\ref{eq:sumpri})). Entonces 
\[\mathbb{E}(X^k) \leq 2\cdot \sum_{l=1}^k \frac{l^k}{l!} \cdot
(\log z)^k \cdot \left(1 + O\left(\frac{1}{\log z}\right)\right)^k .\]

\item Como $\lim_{n\to \infty} (1 + 1/n)^n = e$, tenemos que
 $\left(1 + O\left(\frac{1}{\log z}\right)\right)^k
= e^{O\left(\frac{k}{\log z}\right)}$. Queda por estimar
\[\sum_{l=1}^k \frac{l^k}{l!} \leq k \cdot \max_{1\leq l\leq k} \frac{l^k}{l!} .\]
Como $\log l! = l \log l - l + O(\log l)$ (f\'ormula de Stirling),
tenemos que 
\[\frac{l^k}{l!} = O(l) \cdot \frac{l^k}{l^l e^{-l}} = 
O(l) e^{l + k \log l - l \log l}.\]
Verifique que, para $k$ fijo, la funci\'on $l + k \log l - l \log l$ llega a su m\'aximo
cuando $l$ es la soluci\'on a $k = l \log l$. As\'i,
\[\mathbb{E}(X^k) \ll (\log z)^k e^{O(k/\log z)} \cdot O(k^2) e^{l + k \log l - k},\]
donde $l$ es la soluci\'on a $k = l \log l$.

\item Aplicamos (\ref{eq:aet}):
\begin{equation}\label{eq:astprim}
\Prob(X>\log N) \ll e^{O(k/\log z)} O(k^2) \cdot
\left(\frac{\log z}{\log N}\right)^k e^{l + k \log l - k} .
\end{equation}
Estimaremos $e^{O(k/\log z)} k^2$ al final; nuestra tarea ahora es
encontrar el valor de $k$ para el cual
\begin{equation}\label{eq:mahl}
\left(\frac{\log z}{\log N}\right)^k e^{l + k \log l - k}
= \frac{e^{l + k \log l - k}}{u^k} \end{equation}
 es m\'inimo (donde $l$ es la soluci\'on a $k = l \log l$).
Muestre que
\[\frac{d}{dk} (l + k \log l - k - k \log u) = \log l - \log u.\]
Por lo tanto, el m\'inimo se encuentra cuando $l = u$, es decir, cuando
$k = u \log u$. N\'otese que $k$ no es una constante; por ello\index{momentos!momentos variables}
hablamos de {\em momentos variables}. Quiz\'as resulte algo sorprendente
que el valor \'optimo de $k$ es $u \log u$, puesto que esto es mayor
que $u$ (y, as\'i, $z^k > N$; es por esto que escogimos la cota
(\ref{eq:hariman})).

\item Escogemos, entonces, $k = u \log u$. Obtenemos, por (\ref{eq:astprim}) y 
(\ref{eq:mahl}),
\begin{equation}\label{eq:kraka}\Prob(X>\log N) \ll 
        e^{O\left(u \frac{\log u}{\log z}\right)} 
O(u \log u)^2 \cdot e^{u - u \log u}.\end{equation}
Concluimos que, si $z \to \infty$ y $u = \frac{\log N}{\log z} 
\to \infty$ cuando $N\to \infty$,\index{n\'umeros desmenuzables}
\begin{equation}\label{eq:astrid}
\Prob(\text{$n$ es $u$-desmenuzable y carece de divisores cuadrados}) \ll u^{- u (1 + o(1))} ,\end{equation}
donde $n$ es tomado al azar entre $N$ y $2 N$.
%Podemos mejorar esto un poco. No exijamos ya que $z\to
%\infty$ y $u\to \infty$. Si $z < \frac{1}{2} \log N$ (digamos),
%un entero $n>N$ no puede carecer de divisores
%cuadrados y a\'un as\'i ser $u$-desmenuzable (por qu\'e?). Si $z\geq
%\frac{1}{2} \log N$, (\ref{eq:kraka}) nos da que la probabilidad
%que un entero $n$ sea $u$-desmenuzable y carezca de divisores cuadrados es
%\begin{equation}\label{eq:longsock}\ll u^{- u (1 + O(1/\log u))} .\end{equation}
%Esto es trivialmente cierto para $z< \frac{1}{2} \log N$, por lo
%que acabamos de decir, y para $u$ debajo de una constante; 
%por lo tanto, (\ref{eq:longsock}) es cierto para todo $u$ y todo $z$.

\item La desigualdad (\ref{eq:astrid}) es todo
lo que necesitaremos en nuestra aplicaci\'on m\'as importante. Empero,
es v\'alido preguntarse que pasa si se retira la restricci\'on que
$n$ carezca de divisores cuadrados.
%, o si el entero $n$ se toma al azar entre $1$ y $N$.

Usando (\ref{eq:astrid}), vemos que
\begin{equation}\label{eq:gogo}\begin{aligned}
\Prob&(\text{$n$ es $u$-desmenuzable})\\ &=
\Prob(\text{$n$ es desmenuzable,\; $k^2|n$ y $\frac{n}{k^2}$ es libre de factores
cuadrados})\\
&= \sum_{1\leq k\leq K} \left(\frac{1}{k^2} + O\left(\frac{1}{N}\right)\right) \left(\frac{\log N/k^2}{\log z}\right)^{
-\left(\frac{\log N/k^2}{\log z}\right) (1 + o(1))} \\ 
&+ O\left(\sum_{K<k\leq \sqrt{N}}
\left(\frac{1}{k^2} + \frac{1}{N}\right)\right)\end{aligned}\end{equation}
para todo $K$. Ahora bien, para todo $\epsilon>0$ y todo $k\leq N^{\epsilon}$,
\[\begin{aligned}
\left(\frac{\log N/k^2}{\log z}\right)^{-
\left(\frac{\log N/k^2}{\log z}\right) (1 + o(1))} &= 
\left(\frac{\log N}{\log z} (1 + O(\epsilon)) \right)^{-
\left(\frac{\log N}{\log z}\right) (1 + O(\epsilon)) (1 + o(1))}\\ &= 
\left(\frac{\log N}{\log z}\right)^{- 
\left(\frac{\log N}{\log z}\right) (1 + O(\epsilon))}
\end{aligned}\]
 suponiendo que $z\leq \frac{1}{2} N^{1-2 \epsilon}$, digamos. 
Por lo tanto, haciendo que $K=N^{\epsilon}$, obtenemos de (\ref{eq:gogo})
que \[\Prob(\text{$n$ es $u$-desmenuzable}) = u^{-u (1 + O(\epsilon))} + O\left(
\frac{1}{N^{\epsilon}}\right).\] Si (digamos) $u=(\log N)^{f(N)}$, donde
$f(N)\to \infty$ cuando $N\to \infty$, entonces $N^{-\epsilon} = 
O(u^{-u})$. Haciendo que $\epsilon \to 0$, obtenemos
\begin{equation}\label{eq:uqey}
\Prob(\text{$n$ es $u$-desmenuzable}) = u^{-u (1 + O(\epsilon))} + O(u^{-u})
= u^{-u (1 + o(1))}
\end{equation}
bajo la condici\'on que $z\geq (\log N)^{f(N)}$ para alguna funci\'on $f$
que satisfaga $f(N)\to \infty$ cuando $N\to \infty$.

Una condici\'on de ese tipo, i.e., una cota inferior $z\geq $ para
$z$, es en verdad necesaria para que (\ref{eq:uqey}) sea cierto. Por ejemplo,
pruebe que, si $z$ es una constante, $\Prob(\text{$n$ es $u$-desmenuzable})$
es mucho mayor que $u^{-u (1 + o(1))}$.
\end{enumerate}
%\item Acabamos de ver lo que puede ser la instancia m\'as simple y
%elegante de la siguiente estrategia general: para mostrar que una
%ecuaci\'on tiene pocas soluciones (enteras, racionales, etc.), encontramos
%dos distribuciones naturales que entran en conflicto; toda soluci\'on debe
%ser poco probable respecto a una distribuci\'on o la otra, y por lo tanto hay
%pocas soluciones. El problema, claro est\'a, reside en encontrar estas
%dos distribuciones ``naturales'' -- es decir, distribuciones tales que
%los eventos que caigan en sus colas correspondan a conjuntos de soluciones
%con pocos elementos. 
%moraleja: como hacer - poco com\'un
\end{enumerate}
\section{El l\'imite central}\label{sec:limcen}\index{l\'imite central}
Ya conocemos $\mathbb{E}(\omega(n))$ y $\Var(\omega(n))$ para
$n$ un n\'umero tomado al azar entre $1$ y $N$, donde $N$ es grande.
Quisieramos saber, de una vez por todas, cu\'al es la 
distribuci\'on de $\omega(n)$, en el l\'imite $N\to \infty$.

Como antes, comenzaremos recordando que $\omega(n)$ es una suma de
variables aleatorias, y enfocamos el problema de manera general.

La siguiente observaci\'on se remonta en alguna forma a de Moivre (1718):
si algo es la suma de muchas peque\~nas cosas
que nada o poco tienen que ver entre si, este algo tendr\'a una 
distribuci\'on en forma de campana. Antes de probar tal aseveraci\'on,
debemos ponerla en forma precisa.

\begin{figure}
\centering \includegraphics[height=2.5in]{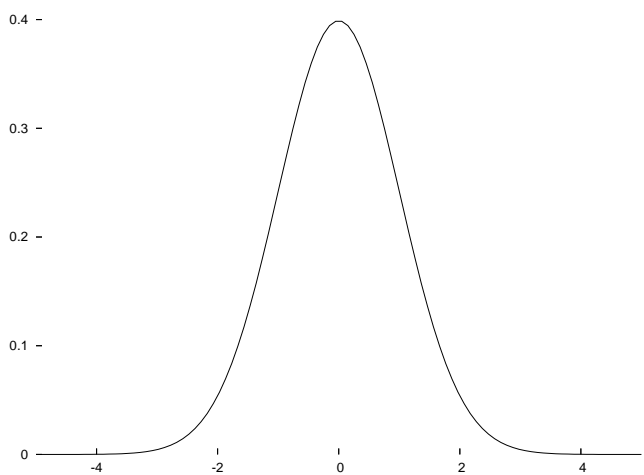}
\caption{La distribuci\'on normal $f(t) = \frac{1}{\sqrt{2 \pi}}
e^{-t^2/2}$. Este es el l\'imite central de la distribuci\'on de $\omega(n)$:
para $n$ tomado al azar entre $1$ y $N$, la probabilidad
$\Prob(\omega(n) \leq \log \log N + t \sqrt{\log \log N})$ tiende a
$\int_{-\infty}^t f(x) dx$.}\label{fig:normal}
\end{figure}

\begin{thm}[Teorema del l\'imite central]\index{l\'imite central!teorema del}
Sean $X_1, X_2, X_3,\dotsc$ variables aleatorias mutuamente independientes.
Asumamos que todas tienen la misma distribuci\'on; sea su esperanza
$E$ y su varianza $V$. 
Asumamos tambi\'en que $\mathbb{E}(X_j^k)$ es finita para todo $k\geq 0$.
Entonces 
$\frac{1}{\sqrt{n V}} \sum_{i=1}^n (X_i - E)$ tiende
en distribuci\'on a 
\begin{equation}\label{eq:comost}\frac{1}{\sqrt{2\pi}} e^{-t^2/2}\end{equation}
cuando $n\to \infty$.
\end{thm}
La distribuci\'on dada por la funci\'on de densidad 
(\ref{eq:comost}) es la afamada {\em distribuci\'on normal}
(ver figura \ref{fig:normal}).\index{distribuci\'on!normal}

Antes de comenzar la demostraci\'on del teorema, recordemos que la {\em transformada de Fourier}\index{transformada de Fourier}
$\hat{f}:\mathbb{R}\to \mathbb{C}$ de una funci\'on $f:\mathbb{R}\to
\mathbb{C}$ se define como sigue:
\[\hat{f}(t) := \int_{-\infty}^{\infty} e^{it x} f(x) dx.\]
La funci\'on 
$f(t) = \frac{1}{\sqrt{2 \pi}} e^{-t^2/2}$ es un {\em vector propio} de
la transformada de Fourier -- es decir, para ese $f$, 
la transformada $\hat{f}(t)$ resulta ser
un m\'ultiplo de $f(t)$: $\hat{f}(t) = \sqrt{2 \pi} f(t)$.
(La prueba est\'a en las notas al final de esta secci\'on.)
Esta propiedad de la funci\'on
$f(t) = \frac{1}{\sqrt{2 \pi}} e^{-t^2/2}$ ser\'a utilizada de manera crucial
hacia el fin de la prueba siguiente.
%Poco importa que se asuma que la esperanza sea $0$ y la varianza sea $1$;
%si no lo son, podemos multiplicar (\ref{eq:comost}) por la raiz cuadrada
%de la varianza y sumarle la esperanza. La suposici\'on t\'ecnica que
%$\mathbb{E}(e^{t X_i})$ sea finita sera satisfecha en todas nuestras
%aplicaciones. M\'as limitante es el hecho que asumamos que las
%variables aleatorias tengan la misma distribuci\'on y sean mutuamente
%independientes; veremos m\'as tarde como relajar estas restricciones.
\begin{proof}
Podemos asumir sin p\'erdida de generalidad que $E=0$ y $V=1$.

Dada una variable aleatoria $X$, definimos la {\em funci\'on caracter\'istica}
\index{funci\'on caracter\'istica}
$\widehat{X}:t\mapsto \mathbb{E}(e^{i t X})$. (Si $X$ es continua,
\begin{equation}\label{eq:arbolito}
\widehat{X}(t) = \mathbb{E}(e^{i t X}) = \int_{-\infty}^{\infty} e^{i t x} f(x) dx,\end{equation}
donde $f$ es la funci\'on de densidad de $X$;
si $X$ es discreta, 
\[\widehat{X}(t) = \mathbb{E}(e^{i t X}) = \sum \Prob(X = x) e^{i t x} .\]
En (\ref{eq:arbolito}) vemos de manera especialmente clara que 
$\widehat{X}$ no es sino una transformada de Fourier.)

Tenemos
\[\widehat{X}(0) = \mathbb{E}(e^{i\cdot 0 \cdot X}) = \mathbb{E}(1) = 1,\]
y, como
\[\widehat{X}'(t) = \mathbb{E}(i X e^{i t X}),\]
se obtiene que $\widehat{X}'(0) = \mathbb{E}(i X e^0) = i \mathbb{E}(X)$.
De la misma manera, $\widehat{X}''(0) = - \mathbb{E}(X^2)$,
$\widehat{X}'''(0) = - i \mathbb{E}(X^3)$, etc.

Sea ahora $X$ cualquiera de las variables $X_i$. Entonces
$\widehat{X}(0) = 1$, $\widehat{X}'(0) = E = 0$, $\widehat{X}''(0) = -V = - 1$,
$|\widehat{X}'''(t)| = |\mathbb{E}(- i X^3 e^{i t X})| \leq
|\mathbb{E}(|X|^3)|\leq 1 + |\mathbb{E}(X^4)| < \infty$. 
Ahora bien, si una funci\'on $f(t)$ es derivable $k+1$ veces alrededor
 del origen 
$t=0$, y $f^{(k+1)}(t)$ est\'a acotada por una constante $c$
 cuando $t$ est\'a cerca
del origen, entonces 
\begin{equation}\label{eq:yunza}
f(t) = f(0) + f'(0)\cdot t + \dotsb + \frac{f^{(k)}(0)}{k!} t^k +
O(c\cdot t^{k+1})
\text{\;\;\;\;\;\;\;\;\;\;\;\; (serie de Taylor truncada)}\end{equation}
cuando $t\to 0$. \begin{small} 
(As\'i como, cuando escribimos ``$O(f(n))$ cuando $n\to \infty$'' queremos
decir 
``entre $- C\cdot f(n)$ y $C \cdot f(n)$ cuando $n$ es mayor que una
constante'', escribimos, similarmente, ``$O(f(t))$ cuando $t\to \infty$''
cuando queremos decir ``entre $-C \cdot f(t)$ y $C \cdot f(t)$ cuando 
$|t|$ es menor que una constante'';  de la misma manera, as\'i como
``$o(f(n))$ cuando $n\to \infty$'' quiere decir ``entre $-g(n) \cdot f(n)$
y $g(n)\cdot f(n)$, donde $g$ es alguna funci\'on con 
$\lim_{n\to \infty} g(n) = 0$'', escribimos
``$o(f(t))$ cuando $t\to 0$'' cuando queremos decir
``entre $-g(t) \cdot f(t)$ y $g(t) \cdot f(t)$, donde 
$\lim_{t\to 0} g(t) = 0$''.)\end{small} Como $c\cdot t \to 0$ cuando $t\to 0$,
vemos por (\ref{eq:yunza}) que 
\[f(t) = f(0) + f'(0)\cdot t + \dotsb + \frac{f^{(k)}(0)}{k!} t^k + o(t^k).\]
Por lo tanto,
\[\widehat{X}(t) = 1 - \frac{t^2}{2} + o(t^2)\]
cuando $t\to 0$.

Para cualquier $r\ne 0$ y cualquier funci\'on $f$, la transformada de
Fourier $\widehat{\frac{1}{r} f}$ de $\frac{1}{r} f$ satisface
$\widehat{\frac{1}{r} f}(t) = \widehat{f}\left(\frac{t}{r}\right)$
(por qu\'e?). Por lo tanto, para $t$ fijo,
\begin{equation}\label{eq:ampar}
\widehat{\frac{1}{\sqrt{n}} X}(t) = \widehat{X}\left(\frac{t}{\sqrt{n}}
\right) = 1 - \frac{t^2}{2 n} + o\left(\frac{t^2}{2 n}\right)
\end{equation}
cuando $n\to \infty$. \begin{small} (F\'ijense en ``$t$ fijo'' y 
``$n\to \infty$''; hemos pasado a nuestro uso habitual de
$o(\cdot)$.)\end{small}

Cuando se tienen dos variables independientes $X$, $Y$, la variable
$X+Y$ tiene como distribuci\'on la convoluci\'on de las distribuciones
(por qu\'e?). 
Ahora bien, la transformada de Fourier $\widehat{f\ast g}$ de la
convoluci\'on $f\ast g$ de dos funciones es igual a $\widehat{f} \cdot
\widehat{g}$ (por qu\'e?).\index{transformada de Fourier}
En consecuencia, $\widehat{X+Y} = \widehat{X} \cdot \widehat{Y}$.
Repetiendo el proceso, obtenemos que la funci\'on caracter\'istica de
$X_1 + X_2 + \dotsc + X_n$ es $\widehat{X_1} \cdot \widehat{X_2}
\cdots \widehat{X_n}$.

Estamos considerando $S_n = \frac{1}{\sqrt{n}} \sum_{j=1}^n X_j$. 
Vemos ahora que, por (\ref{eq:ampar}), su funci\'on caracter\'istica debe ser
\begin{equation}\label{eq:rosta}
\widehat{S_n} = \left(\widehat{X}\left(\frac{t}{\sqrt{n}}\right)\right)^n =
\left(1 - \frac{t^2}{2 n} + o\left(\frac{t^2}{n}\right)\right)^n\end{equation}
cuando $n\to \infty$. 
%\mathop{\longrightarrow}^{n\to \infty}
% e^{-t^2/2} .
%\]

Para $\epsilon$ peque\~no, tenemos $(1 + \epsilon) = 
e^{\epsilon + O(\epsilon^2)}$. Por lo tanto, 
\[\left(1 - \frac{t^2}{2 n} + o\left(\frac{t^2}{2 n}\right)\right)^n
= \left(\left(1 + \frac{1}{n}\right)^{-\frac{t^2}{2} (1 + o(1))}\right)^n
= \left(\left(1 + \frac{1}{n}\right)^n\right)^{-\frac{t^2}{2} (1 + o(1))}\]
cuando $n\to \infty$. Ahora bien,
\[\left(1 + \frac{1}{n}\right)^n \to e\]
cuando $n\to \infty$. As\'i, concluimos que
\begin{equation}\label{eq:dobast}
\left(1 - \frac{t^2}{2 n} + o\left(\frac{t^2}{2 n}\right)\right)^n \to
e^{-t^2/2}
\end{equation}
cuando $n\to \infty$. Por (\ref{eq:rosta}) y (\ref{eq:dobast}), obtenemos
finalmente que
\[\widehat{S_n} \to e^{-t^2/2}
\]
cuando $n\to \infty$. \begin{small} (Puede que la velocidad de convergencia
dependa de $t$, pero esto no nos importa; el resultado que estamos por
utilizar es robusto en ese sentido.)\end{small}

%Tambi\'en tenemos $\widehat{\frac{1}{r} X}(t) = \widehat{X}(t/r)$.

La transformada $\widehat{W}$ de la normal $W =
\frac{1}{\sqrt{2 \pi}} e^{-t^2/2}$ es precisamente
$e^{-t^2/2}$ (ver las notas). Tenemos, entonces, que, para cada
$t$, $\widehat{S_n}(t)$
tiende a $\widehat{W}(t)$ cuando $n\to \infty$. 
%La transformada es invertible, ya que
%$\widehat{\widehat{X}} = 2 \pi X$ para cualquier $X$ ``sensato'';
Invocamos un resultado del an\'alisis ({\em teorema de convergencia de
  L\'evy};\index{teorema!de convergencia de L\'evy}
ver notas) y concluimos que $S_n$ tiende a $W$ cuando $n\to \infty$.
\end{proof}

La idea central de la siguiente prueba alternativa nos
ser\'a de utilidad cuando examinemos $\omega(n)$. El m\'etodo es 
llamado {\em m\'etodo de momentos}.\index{momentos!m\'etodo de momentos}
\begin{proof}[Esbozo de otra demostraci\'on.]\index{l\'imite central!teorema del}
Compararemos los {\em momentos} \[\mathbb{E}(S_n),\; \mathbb{E}(S_n^2),\;
\mathbb{E}(S_n^3), \dots\]\index{momentos}
 de la variable $S_n = \frac{1}{\sqrt{n}}
\sum_{j=1}^n X_j$ con los momentos $\mathbb{E}(S)$, $\mathbb{E}(S^2)$,
$\mathbb{E}(S^3)$, \dots
de una variable $S$ de distribuci\'on normal.

Por integraci\'on por partes, podemos ver que
$\mathbb{E}(S^k) = (k-1) (k-3) (k-5)\dotsb 3\cdot 1$ para $k$ par;
como $S$ es sim\'etrica con respecto al eje $y$, est\'a claro que
$\mathbb{E}(S^k) = 0$ para $k$ impar.
Podemos verificar que $\mathbb{E}(S_n^k) = 
(k-1) (k-3) (k-5)\dotsb 3\cdot 1 + o_k(1)$ para $k$ par,
y $\mathbb{E}(S_n^k) = o_k(1)$ para $k$ impar (ver el problema \ref{prob:mom}).

Como los momentos de $S_n$ convergen a los momentos de $S$ y
la distribuci\'on normal
satisface ciertas condiciones t\'ecnicas, podemos concluir que
$S_n \to S$ utilizando un resultado auxiliar est\'andar (nota
\ref{not:congmom}).
\end{proof}

%En general, si $S$ tiene una distribuci\'on ``sensata'' (en particular,
%la normal), y una sucesi\'on de variables $S_n$ es tal que
%$\lim_{n\to \infty} \mathbb{E}(S_n^k) = \mathbb{E}(S^k)$ para todo $k$,
%entonces las distribuci\'on de $S_n$ converge a la distribuci\'on de $S$
%cuando $n\to \infty$. 

{\em Condiciones del teorema del l\'imite central.}\index{l\'imite central!teorema del}
Hemos asumido tres cosas acerca de las variables
$X_j$: (a) que son mutuamente 
independientes, (b) que tienen la misma distribuci\'on,
(c) que, para cada $k$, $\mathbb{E}(X_j^k)$ est\'a acotada independientemente
de $j$ (lo cual es lo mismo que $\mathbb{E}(X_j^k)<\infty$, si (b)
se cumple).

Tanto (b) como (c) pueden relajarse; la {\em condici\'on de Lindeberg}\,
\index{condici\'on de Lindeberg}
funge por las dos (nota \ref{prob:lind}). Es m\'as dif\'icil
prescindir de (a); hay algunas herramientas est\'andar para tal tarea,
pero ninguna cubre todos los casos que aparecen en la pr\'actica.

Antes de ver como podemos arregl\'arnoslas sin (a), hagamos dos cosas:
primero, verifiquemos que la falta de (b) en el caso que m\'as nos
interesa es inocua; luego, veamos como, en muchos otros casos, la falta de (b)
(y de la condici\'on de Lindeberg) hace que la conclusi\'on sea
falsa -- es decir, que el l\'imite no sea normal.

Sean $X_2', X_3', X_5',\dotsc$ variables mutuamente independientes con
la siguiente distribuci\'on:
\begin{equation}\label{eq:ghae}X_p' = \begin{cases} 1 &\text{con probabilidad $1/p$}\\
0 &\text{con probabilidad $1-1/p$.} \end{cases}\end{equation}
(Escojemos los signos $X_2', X_3', X_5',\dotsc$ porque usaremos
$X_2, X_3, X_5,\dotsc$ m\'as tarde.)

Entonces $\mathbb{E}(X_p') = \frac{1}{p}$, 
$\mathbb{E}((X_p' - \mathbb{E}(X_p'))^2) = \Var(X_p') = \frac{1}{p} 
- \frac{1}{p^2}$, $\mathbb{E}(|X_p' - \mathbb{E}(X_p')|^3) \leq \frac{1}{p}$. 
Por consiguiente (ver (\ref{eq:yunza}) y las l\'ineas inmediatamente precedentes),
 la funci\'on caracter\'istica de $X_p' - 
\mathbb{E}(X_p')$ es
\[1 - \frac{t^2}{2} \left(\frac{1}{p} - \frac{1}{p^2}\right)
 + \frac{O(t^3)}{p}.\]
Por el mismo razonamiento, la funci\'on caracter\'istica de 
$\frac{1}{\sqrt{\log \log p}} (X_p' - \mathbb{E}(X_p'))$ es
\begin{equation}\label{eq:spirave}
1 - 
\frac{t^2}{2 \log \log n} \left(\frac{1}{p} - \frac{1}{p^2}\right)
+ \frac{O(t^3)}{p (\log \log n)^{3/2}}.
\end{equation}

Definimos $S_n' = \frac{1}{\sqrt{\log \log n}} \sum_{p\leq n} (X_p'
- \mathbb{E}(X_p'))$. Usando la regla $\widehat{X+Y} = \widehat{X} \cdot
\widehat{Y}$ y (\ref{eq:spirave}), vemos que
\[\begin{aligned} \widehat{S_n'} &=
 \prod_{p\leq n} \left(1 - 
\frac{t^2}{2 \log \log n} \left(\frac{1}{p} - \frac{1}{p^2}\right)
+ \frac{O(t^3)}{p (\log \log n)^{3/2}}\right)\\
&= \prod_{p\leq n} e^{- \frac{t^2}{2 \log \log n} \cdot \left(\frac{1}{p} - 
\frac{1}{p^2}\right) \cdot
(1 + O(t/\sqrt{\log \log n}))}\\
&= e^{- (1 + o_t(1))\cdot  \sum_{p\leq n} \frac{t^2}{2 \log \log n}
\left(\frac{1}{p} - \frac{1}{p^2}\right)} 
\end{aligned}\]
cuando $n\to \infty$.
Por el teorema de Chebyshev-Mertens (\ref{eq:cruj}),
\[\sum_{p\leq n} \frac{t^2}{2 \log \log n}
\left(\frac{1}{p} - \frac{1}{p^2}\right) = \frac{t^2}{2} (1 + O(1/\log
\log n))\]
cuando $n\to \infty$,
y en consecuencia
\[\widehat{S_n'} = e^{-t^2/2 \cdot (1 + o_t(1))} \to
e^{-t^2/2}\]
cuando $n\to \infty$. Por el teorema de convergencia de Levy,
concluimos que $S_n'$ converge en distribuci\'on a la 
normal (\ref{eq:comost}).

\begin{center}
* * *
\end{center}

Consideremos, en cambio, variables 
$X_2', X_3', X_5',\dotsc$ mutuamente independientes
tales que
\[X_p' = \begin{cases} 1 &\text{con probabilidad $\rho_p$}\\
0 &\text{con probabilidad $1 - \rho_p$,}\end{cases}\]
donde $\rho_2, \rho_3,\dotsc$ son tales que $\sum_p \rho_p$ converge.
Entonces $\sum_p \mathbb{E}(X_p') < \infty$ y
$\Var(\sum_p X_p') = \sum_p \Var(X_p') < \infty$. 
La funci\'on de distribuci\'on de 
$S_n' = \sum_{p\leq n} X_p'$ tender\'a a un l\'imite no normal 
$f(x)$. Ver el problema \ref{prob:astra}.

%Por el teorema de Chebyshev, se deduce que hay una constante $C$ tal que, para
%todo $n$, $\Prob(\sum_{p\leq n} X_p < C) > 0.9999$ (digamos) 
%y una constante $N$ tal que $\Prob(\sum_{p>N} X_p > 0.0001) <0.0001$
%(digamos). Est\'a
%claro, entonces, que 

%un factor de escala, como $\log \log n$.)
%; si tal cosa se introdujera, 
%obtendr\'iamos el valor $0$ con probabilidad $1$, en el l\'imite

%Ahora bien, $f(x) = 0$ para todo $x$ negativo (ya que $S_n$ adopta
%solo valores no negativos). Al mismo tiempo, $f(x)>0$ para todo $x$ positivo
%(ya que $S_n$ puede adoptar valores arbitrariamente grandes). Por
%consiguiente,
%$f$ no es sim\'etrica alrededor de $\lim_{n\to \infty} \mathbb{E}(S_n) 
%= \sum_p \rho_p$, o alrededor de ning\'un otro valor finito.
%Esto quiere decir que $f$ no es la distribuci\'on normal.

%El l\'imite $f$, en verdad, ser\'a en alg\'un sentido similar a la
%{\em distribuci\'on de Poisson}, aunque generalmente no sera igual a \'esta;
%se tratar\'a de alguna otra funci\'on de {\em cola pesada}, dependiente
%de cada $\rho_p$. 

\begin{center}
* * *
\end{center}

{\em El l\'imite central de $\omega(n)$.}\index{l\'imite central} Sean, ahora, $X_2, X_3, X_5,\dotsc$
variables dadas por
\begin{equation}\label{eq:mae}
X_p = \begin{cases} 1 &\text{si $p|n$,}\\ 0 &\text{si $p\nmid n$,}\end{cases}
\end{equation}
donde $n$ es un entero aleatorio entre $1$ y $N$.
Como ya sabemos, $\omega(n) = \sum_{p\leq n} X_p$. Queremos
probar que la distribuci\'on de $\omega(n)$ (o, m\'as bien dicho,
$\frac{1}{\sqrt{\log \log n}} (\omega(n) - \log \log n)$) tiende a la normal.

Cuando calculamos la varianza de $\omega(n)$, vimos que las variables
$X_p$ son casi independientes en pares: $X_{p_1}$ y $X_{p_2}$ son
aproximadamente independientes para $p_1, p_2 < N^{1/2 - \epsilon}$, $p_1\ne
p_2$ cualesquiera, y, en total, los t\'erminos de error son peque\~nos.
Empero, las variables $X_p$ est\'an muy lejos de ser mutuamente 
independientes. Qu\'e podemos hacer?

Podemos probar el teorema del l\'imite central para $\omega(n)$ por
el m\'etodo de momentos. Cuando calculamos el momento 
$\mathbb{E}(\omega(n)^k)$, s\'olo necesitamos el hecho que las variables
$X_p$ sean casi independientes ``de a $k$'': $k$ variables distintas 
cualesquiera
entre $X_2, X_3, X_5,\dotsc$ son aproximadamente independientes, por la
mismas razones que ya vimos para $k=2$. El t\'ermino de error depender\'a
de $k$, y por lo tanto la tasa de convergencia de
$\mathbb{E}(\omega(n)^k)$ a su l\'imite depender\'a de $k$; empero, al
m\'etodo de momentos esto no le importa (nota \ref{not:congmom}).

Pasemos esto en limpio. 
\begin{thm}[Erd\H{o}s-Kac \cite{EK}]\label{thm:erdkac}
\index{teorema!de Erd\H{o}s-Kac} \index{Erd\H{o}s, P.}
\index{Kac, M.} Sea $n$ un entero tomado
al azar entre $1$ y $N$ con la distribuci\'on uniforme. Entonces
 $Y_N = \frac{1}{\sqrt{\log \log N}}
(\omega(n) - \log \log N)$ tiende en distribuci\'on a
\[\frac{1}{\sqrt{2 \pi}} e^{-t^2/2}
\]
cuando $N\to \infty$.
\end{thm}
La demostraci\'on que veremos se debe a Billingsley;\index{Billingsley, P.} incluye
varias ideas de la prueba de Erd\H{o}s y Kac (1939) 
\index{teorema!de Erd\H{o}s-Kac} \index{Erd\H{o}s, P.} \index{Kac, M.}
y de la prueba de
Halberstam (1955). \index{Halberstam, H.}

Como antes, expresaremos $\omega(n)$ como una suma 
$\sum_{p\leq N} X_p$.
As\'i como, para utilizar la desigualdad de
Chebyshev\index{desigualdad!de Chebyshev}, tuvimos que truncar la suma
$\sum_{p\leq N} X_p$ (reemplaz\'andola por $\sum_{p\leq N^{1/3}} X_p$),
tendremos que truncarla ahora (a\'un m\'as, ya que la
reemplazaremos por $\sum_{p\leq g(N)} X_p$, donde $g(x)$ crece
m\'as lentamente que cualquier potencia de $x$). Primero mostraremos
que esta truncaci\'on no nos es da\~nina -- es decir, que el total
omitido $\sum_{g(N)<p\leq N} X_p$ es peque~no; luego procederemos
a determinar la distribuci\'on de la suma truncada
$\sum_{p\leq g(N)} X_p$ utilizando el m\'etodo de momentos.
\begin{proof}
Sea $g(x)$ una funci\'on tal que $g(x) = o_{\epsilon}(x^{\epsilon})$ para
todo $\epsilon>0$ y $\log \log x - \log \log g(x) = o(\sqrt{\log \log x})$;
podemos tomar, por ejemplo, $g(x) = x^{1/\log \log x}$.
%Entonces \[Y_N = \frac{1}{\sqrt{\log \log N}}
%(\sum_{p\leq N} X_p - \log \log N),\] 
%donde $X_p$
%es como en (\ref{eq:mae}).

Vemos inmediatamente que
\[
Y_N = (1 + o(1)) \cdot \frac{1}{\sqrt{\log \log g(N)}}
(\omega(n) - \log \log N) .\]
Ahora bien,
\[\begin{aligned}
\sum_{p\leq g(N)} 1/p &= \log \log g(N)\, +\, O(1) =
\log \log N + o(\sqrt{\log \log N}) \\ 
&= \log \log N + o(\sqrt{\log \log g(N)})
.\end{aligned}\]
Concluimos que
\begin{equation}\label{eq:ravel}
Y_N = (1 + o(1)) \cdot \frac{1}{\sqrt{\log \log g(N)}}
\left(\omega(n) - \sum_{p\leq g(N)} 1/p\right) + o(1).\end{equation}

Definamos ahora
\begin{equation}\label{eq:boulez}S_m = \frac{1}{\sqrt{\log \log m}} \cdot 
(\sum_{p\leq m} (X_p - 1/p)),\end{equation}
donde $X_p$ es como en (\ref{eq:mae}).
Est\'a claro que 
$\omega(n) = \sum_{p\leq N} X_p$, y por lo tanto
$\omega(n) - \sum_{p\leq g(N)} X_p$ es igual a
$\sum_{g(N) < p\leq N} X_p$. Calculamos la esperanza de esto \'ultimo:
\begin{equation}\label{eq:barbar}\begin{aligned}
 \mathbb{E}\left(\sum_{g(N) < p\leq N} X_p\right) &=
\sum_{g(N) < p \leq N} \frac{1}{p} 
= \log \log N - \log \log g(N) + O(1) \\ &= 
o\left(\sqrt{\log \log N}\right) =
o\left(\sqrt{\log \log g(N)}\right). \end{aligned}\end{equation}
Como $X_p$ toma s\'olo valores positivos, la desigualdad de Markov
nos permite deducir de (\ref{eq:barbar}) que
 \[\frac{1}{\sqrt{\log \log g(N)}} \left(\sum_{g(N) < p \leq N} X_p
\right) = o(1)\] con probabilidad $1 - o(1)$.
Concluimos (por (\ref{eq:boulez}) y (\ref{eq:ravel}))
que \begin{equation}\label{eq:moesg}
Y_N = (1 + o(1)) S_{g(N)} + o(1)\end{equation}
con probabilidad $1 - o(1)$ cuando $N\to \infty$. En consecuencia, 
si probamos que $S_{g(N)}$ tiende en distribuci\'on a la normal, habremos
probado que $Y_N$ tiende en distribuci\'on a la normal.

{\small Hasta ahora, nuestra labor ha sido s\'olo preparatoria: lo m\'as
que hemos hecho es truncar la suma $\sum_{p\leq N}$ y mostrar que 
el efecto de tal truncaci\'on es peque\~no.}
A continuaci\'on, nuestra tarea es averiguar cu\'ales son
 los momentos de $S_{g(N)}$. Sea $k\geq 0$.
Sean $X_p'$ como en (\ref{eq:ghae})
y $S_m' = \frac{1}{\sqrt{\log \log m}} \sum_{p\leq m} (X_p' -
1/p)$.
Para $k$ primos $p_1, p_2,\dotsc,p_k$ cualesquiera
(no necesariamente distintos), 
\[|\mathbb{E}(X_{p_1} X_{p_2}\dotsc X_{p_k}) -
 \mathbb{E}(X_{p_1}' X_{p_2}'\dotsc X_{p_k}')| =
\left|\frac{1}{N} \left\lfloor \frac{N}{d}\right\rfloor - \frac{1}{d} \right| 
\leq \frac{1}{N},\]
donde $d$ es el m\'inimo com\'un m\'ultiplo de $p_1, p_2,\dotsc,p_k$.
Por lo tanto,
\[\mathbb{E}(S_{g(N)}^k) = \mathbb{E}((S_{g(N)}')^k) +
O_k\left(g(N)^k \cdot \frac{1}{N}\right),\]
puesto que $g(N)^k$ es el n\'umero de t\'erminos que aparecen cuando se
expande $(S'_{g(N)})^k$.
Como $g(N) = o_{\epsilon}(N^{\epsilon})$, sabemos que
$O\left(g(N)^k\cdot \frac{1}{N}\right) = o_k(1)$.

Ya vimos (despu\'es de (\ref{eq:ghae})) que
la distribuci\'on de $S_m'$ tiende a la normal cuando $m\to \infty$; 
por lo tanto,
los momentos $\mathbb{E}((S_m')^k)$ de $S_m$ tienden a los 
momentos $\mathbb{E}(W^k)$ de la normal $W$. 
Tenemos, entonces, que
\[\lim_{N\to \infty} \mathbb{E}(S_{g(N)}^k) = \lim_{N\to \infty} 
\mathbb{E}((S_{g(N)}')^k) = \mathbb{E}(W^k)\]
para todo $k$. Concluimos que $S_{g(N)}^k$ 
converge en distribuci\'on a la normal, y, por lo tanto, $Y_N$ converge
en distribuci\'on a la normal.
\end{proof}

\begin{center}
{\bf Notas y problemas}
\end{center}
\begin{enumerate}
\item Debemos probar que la transformada de Fourier\index{transformada de Fourier}
\[\hat{f}(t) = \int_{-\infty}^{\infty} e^{i t x} f(x) dx\]
de $f(x) = \frac{1}{\sqrt{2 \pi}} e^{-x^2/2}$ es igual a $e^{-t^2/2}$. 
Una de las maneras m\'as simples es la siguiente. Tenemos que
\[\begin{aligned}
\hat{f}(t) &= \int_{-\infty}^{\infty} e^{i t x} \cdot \frac{1}{\sqrt{2 \pi}}
e^{-x^2/2} dx\\
&= \frac{1}{\sqrt{2\pi}} \int_{-\infty}^{\infty} e^{-\frac{1}{2} (x^2 - 2 i t x)}
dx.\end{aligned}\]
Completando cuadrados,
\[\begin{aligned}
\widehat{f}(t) &= \frac{1}{\sqrt{2\pi}} \cdot \int_{-\infty}^{\infty}
e^{(-i t)^2/2} e^{-(x-i t)^2/2} dx\\
&= e^{-t^2/2} \cdot \frac{1}{\sqrt{2\pi}} \int_{-\infty}^{\infty}
 e^{-(x-it)^2/2} dx.
\end{aligned}\]
Por el teorema de Cauchy en el an\'alisis complejo\footnote{El lector que no
quiera utilizar el an\'alisis complejo (el cual estamos evitando en general) puede saltarse este p\'arrafo y ver la prueba alternativa al final de la nota
presente.}, aplicado a la funci\'on $e^{-z^2/2}$ (anal\'itica en todo el
plano complejo), tenemos que
\begin{equation}\label{eq:gobot}
\int_{-\infty}^\infty e^{-(x-it)^2/2} dx = \int_{-\infty}^{\infty} e^{-x^2/2} dx
\end{equation}
As\'i tenemos que $\widehat{f}(t) = c \cdot e^{- t^2/2}$, donde
$c$ es la constante $c = \frac{1}{\sqrt{2\pi}} 
\int_{-\infty}^{\infty} e^{-x^2/2} dx$. S\'olo falta calcular $c$.
Haremos esto en una de las formas m\'as conocidas - 
evaluando la integral de $e^{-(x^2+y^2)}$ en el plano en dos
maneras distintas.

Por una parte,
\begin{equation}\label{eq:bull}\begin{aligned}
\int_{-\infty}^{\infty} \int_{-\infty}^{\infty} e^{-(x^2+y^2)/2} dx dy
&= \left(\int_{-\infty}^{\infty} e^{-x^2/2} \right) \cdot
\left(\int_{-\infty}^{\infty} e^{-y^2/2} \right)\\
&= \left(\int_{-\infty}^{\infty} e^{-x^2/2} \right)^2
\end{aligned}\end{equation}
Por otra parte, cambiando de coordenadas rectangulares $(x,y)$ a
coordenadas polares $(r,\theta)$, obtenemos
\[\int_{-\infty}^{\infty} \int_{-\infty}^{\infty} e^{-(x^2+y^2)} dx dy
= \int_0^{2\pi} \int_0^\infty e^{-r^2/2} r dr d\theta,
\]
puesto que $dx dy = r dr d\theta$. Un breve c\'omputo muestra que
\begin{equation}\label{eq:dosest}
\int_0^{2\pi} \int_0^\infty e^{-r^2/2} r dr d\theta
= 2\pi \int_0^{\infty} e^{-r^2/2} r dr = 2\pi \cdot 1 = 2\pi .
\end{equation}
Comparando (\ref{eq:bull}) y (\ref{eq:dosest}), vemos que
\[\int_{-\infty}^{\infty} e^{-x^2/2} dx = \sqrt{2 \pi}.\]
Por lo tanto, $c = 1$, y, as\'i, concluimos que $\widehat{f}(t) = 
e^{-t^2/2}$.

Como, en general, no estamos asumiendo ning\'un conocimiento del an\'alisis
complejo, es bueno indicar una prueba alternativa que no utilize
el teorema de Cauchy (el cual usamos en el paso (\ref{eq:gobot})).
Una manera conocida es la que sigue. La funci\'on 
$f(x) = \frac{1}{\sqrt{2 \pi}} e^{-x^2/2}$ puede ser descrita como la
soluci\'on (necesariamente \'unica) al problema de valor inicial
dado por las siguientes condiciones: $f'(x)= - x f(x)$, $f(0)=\frac{1}{
\sqrt{2\pi}}$.
Tenemos
\[\begin{aligned}
\widehat{f}'(t) &=
\int_{-\infty}^{\infty} \left(\frac{d}{dt} (e^{i t x} f(x))\right) dx= \int_{-\infty}^{\infty} i x e^{i t x} f(x) dx
= -i \int_{-\infty}^{\infty} (- x f(x)) e^{i t x} dx\\
&= -i \int_{-\infty}^{\infty} f'(x) e^{i t x} dx
\end{aligned}
\]
Hacemos una integraci\'on por partes, y obtenemos
\[\begin{aligned}
\widehat{f}'(t) =
i \int_{-\infty}^{\infty} f(x) \frac{d}{dx} (e^{i t x}) dx
= -t \int_{-\infty}^{\infty} f(x) e^{i t x} dx
= -t \widehat{f}(t).
\end{aligned}\]
Tambi\'en tenemos que \[\widehat{f}(0) = \frac{1}{\sqrt{2\pi}} 
\int_{-\infty}^{\infty} e^{-x^2/2} dx = 1,\]
como probamos con anterioridad. As\'i, la funci\'on $g(t) = \frac{1}{\sqrt{2\pi}}
\widehat{f}(t)$ satisface las condiciones $g'(t) = -t g(t)$, $g(0)=
\frac{1}{\sqrt{2\pi}}$. Como $f$ es la \'unica funci\'on que satisface tales
condiciones, concluimos que $g(t) = f(t)$. Por lo tanto,
\[\widehat{f}(t) = \sqrt{2\pi} g(t) = \sqrt{2\pi} f(t) = e^{-t^2/2},\]
como quer\'iamos demostrar.

\item El siguiente es un resultado muy \'util del an\'alisis 
``suave''; la prueba
involucra argumentos de compacticidad y convergencia, aparte de una
transformaci\'on oportuna.
\begin{thm}[Teorema de convergencia de P.\ L\'evy]\index{
teorema!de convergencia de L\'evy}\index{P. L\'evy}
Sean $X_1, X_2, X_3,\dotsc$ variables aleatorias
con funciones caracter\'isticas $\widehat{X_1}, \widehat{X_2},
\widehat{X_3},\dotsc$. Asumamos que,
para todo real $t$, la sucesi\'on $\widehat{X_1}(t), \widehat{X_2}(t),\dotsc$
tiene un l\'imite $f(t)$. Si $f$ es continua alrededor de $t=0$,
entonces $f$ es la funci\'on caracter\'istica $\widehat{X}$ de alguna
variable aleatoria $X$, y las variables $X_1, X_2, X_3,\dotsc$ convergen
a $X$ en distribuci\'on.
\end{thm}
La prueba puede encontrarse en \cite[Vol. 2, \S XV.3, Teorema 2]{Fe}.

\item\label{prob:mom} Sean $X_1, X_2, X_3,\dotsc$ variables independientes
con esperanza $0$, varianza $1$, y $\mathbb{E}(X_j^r)$ acotada para todo
$j$ y todo entero $r$ entre $0$ y $k$. Queremos estimar la
esperanza de
$(X_1 + X_2 + \dotsc + X_n)^k$. Comencemos expandiendo esta potencia
en sus t\'erminos: 
$(X_1 + X_2 + \dotsc + X_n)^k = 
X_1^k + X_1^{k-1} X_2 + \dotsb + X_1 X_2^2 X_4^{k-4} X_n + \dotsb $.

\begin{enumerate}
\item Muestre que los t\'erminos donde alg\'un $X_j$ aparece a la potencia $1$ 
tienen esperanza $0$. (Decimos que $X_j$ aparece a la potencia $\alpha$
si el t\'ermino es de la forma $\dotsc X_j^{\alpha} \dotsc$; por ejemplo,
$X_j$ aparece a la potencia $\alpha$ en $X_1 X_j^{\alpha} X_k^2$.)

\item Muestre que hay a lo m\'as $O_k(n^{(k-1)/2})$ t\'erminos 
(o $O_k(n^{k/2 - 1})$ si $k$ es par) donde
aparecen s\'olo potencias $\geq 2$ y por lo menos una potencia $\geq 3$.

\item Nos quedan los t\'erminos donde toda variable que aparece,
aparece a la potencia $2$. Tal cosa puede ocurrir s\'olo cuando
$k$ es par; escribamos $k = 2 l$. Muestre que cada t\'ermino
de la forma antedicha (es decir, cada t\'ermino que contenga s\'olo
cuadrados, e.g., $X_1^2 X_3^2 X_{10}^2$) ocurre
exactamente $\frac{(2 l)!}{2^l}$ veces.

\item Si $k$ es impar, concluimos que
$\mathbb{E}((X_1 + X_2 + \dotsc + X_n)^k) = O_{k,c}(n^{(k-1)/2})$.
Si $k$ es par, concluimos que
$\mathbb{E}((X_1 + X_2 + \dotsc + X_n)^k)$ es igual
al n\'umero de t\'erminos distintos que contengan solo cuadrados,
multiplicado por $\frac{(2 l)!}{2^l}$, m\'as
$O_{k,c}(n^{k/2 - 1})$. A continuaci\'on, estimaremos este n\'umero
de t\'erminos de manera indirecta.

\item Expandamos la expresi\'on $(x_1 + x_2 + \dotsc + x_n)^l$. 
Hay a lo m\'as $O_k(n^{l-1})$
 t\'erminos donde aparecen potencias $>1$. Muestre que cada t\'ermino
donde no aparecen potencias $>1$ ocurre $l!$ veces. 

\item Los t\'erminos de $(x_1 + x_2 + \dotsc + x_n)^l$
donde no aparecen potencias $>1$
est\'an en correspondencia uno a uno -- sin contar el n\'umero
de ocurrencias -- con los t\'erminos de
$(X_1 + X_2 + \dotsc + X_n)^{2 l}$ que contienen s\'olo cuadrados.

\item Definamos ahora $x_j = 1$ para todo $j$; entonces 
$(x_1 + x_2 + \dotsc + x_n)^l$ se vuelve $n^l$. 
Concluya que
hay $\frac{1}{l!} n^l + O_k(n^{l-1})$ t\'erminos 
distintos donde no aparecen potencias $>1$.

\item Obtenemos inmediatamente que
hay $\frac{1}{l!} n^l + O_k(n^{l-1})$ t\'erminos 
distintos en 
$(X_1 + X_2 + \dotsc + X_n)^k$ que contienen s\'olo cuadrados.
Concluya que \[\mathbb{E}((X_1 + X_2 + \dotsc + X_n)^k) =
((k-1) \cdot (k-3) \dotsb 3 \cdot 1) \cdot n^{k/2} + O_k(n^{k/2-1})\]
para $k$ par, y $\mathbb{E}((X_1 + X_2 + \dotsc + X_n)^k) = 0$ para $k$
impar.\\
\end{enumerate}

\item\label{not:congmom}
El m\'etodo de momentos es v\'alido gracias al resultado siguiente.
% los dos siguientes resultados.

%\begin{lem}
%Sea $X$ una variable aleatoria tal que $\frac{1}{k!}
%|\mathbb{E}(X^k)| < C^k$ para alg\'un $C>0$ y todo $k>0$.
%Entonces $X$ est\'a {\em determinada por sus momentos}:
%si una variable $Y$ satisface $\mathbb{E}(Y^k) = \mathbb{E}(X^k)$
%para todo $k>0$, entonces $Y$ tiene la misma distribuci\'on que $X$. 
%\end{lem}
%Decimos que $X$ e $Y$ tienen la misma distribuci\'on si
%$\Prob(X\leq t) = \Prob(Y\leq t)$ para todo $t$.
%\begin{proof}[Esbozo de una prueba]
%Como
%$\widehat{X}^{(n)}(0) = i^k \mathbb{E}(X^k)$, los momentos nos dan la serie
%de Taylor de $\widehat{X}$ alrededor del origen. La condici\'on
%$\frac{1}{k!} |\mathbb{E}(X^k)| < C^k$ nos asegura que la serie
%tiene radio de convergencia positivo. Con un poco m\'as de trabajo,
%obtenemos que la serie alrededor de todo punto $x=a$ tiene radio de 
%convergencia positivo. As\'i, podemos continuar
%la funci\'on $\widehat{X}$
%anal\'iticamente (de manera \'unica) a toda la recta real.
%Las series de Taylor alrededor de cada $x=a$
% est\'an determinadas por los momentos $\mathbb{E}(X^k)$;
%hemos probado, por ende, que $\widehat{X}$ est\'a determinada
%por los momentos $\mathbb{E}(X^k)$. 

%En otras palabras, $\widehat{X} = \widehat{Y}$. Invertimos la transformada
%de Fourier {\small (cuidado!)} y obtenemos que $X$ e $Y$ tienen la misma
%distribuci\'on.
%\end{proof}

\begin{thm}
Sean $X_1, X_2, X_3 \dotsc$ y $X$ variables aleatorias tales que
$\mathbb{E}(X_j^k)$ y $\mathbb{E}(X^k)$ son finitos para
$j,k\geq 0$ cualesquiera.
 Supongamos que los momentos de $X_j$
convergen a los momentos de $X$:
$\lim_{j\to \infty} \mathbb{E}(X_j^k) = \mathbb{E}(X^k)$.
Supongamos tambi\'en que $\mathbb{E}(X^k)<C^k$ para alg\'un $C>0$ y todo
$k>0$.
Entonces
$X_1, X_2, X_3,\dotsc$ convergen en distribuci\'on a $X$.
\end{thm}

La idea principal de la prueba es que los momentos de una distribuci\'on
$X$ determinan la serie de Taylor de $\widehat{X}(t)$ alrededor de $t=0$.
La condici\'on $\mathbb{E}(X^k)<C^k$ asegura que la serie de Taylor alrededor
de $t=0$ tenga radio de convergencia infinito; as\'i, la serie determina
$\widehat{X}$, y, por ende, determina $X$. Para ver una prueba completa,
consultar, por ejemplo, \cite[\S 30, Teoremas 30.1 y 30.2]{Bi}.

La condici\'on $\frac{1}{k!} |\mathbb{E}(X^k)| < C^k$ se cumple
para casi toda variable $X$ ``razonable''. He aqu\'i un ejemplo de un $X$
para el cual la condici\'on no se cumple, y, m\'as a\'un, la conclusi\'on
del teorema no es cierta: sea $X = e^Y$, donde $Y$ es una variable
de distribuci\'on normal.\index{distribuci\'on!normal}

\item\label{prob:lind}
Las condiciones del teorema del l\'imite central se pueden relajar de
varias formas. La siguiente es una de las formas m\'as comunes.
\begin{thm}[Teorema del l\'imite central -- Lindeberg]\
Sean $X_1, X_2, X_3,\dotsc$ variables aleatorias mutuamente independientes.
Sean
\[S_n = \sum_{j=1}^n (X_j - \mathbb{E}(X_j)),\;\;\;\;\;
s_n= \sqrt{\Var(S_n)} = \sqrt{\sum_{j=1}^n \Var(X_j)} .\]
Supongamos que, para todo $\epsilon>0$,\index{condici\'on de Lindeberg}
\begin{equation}\label{eq:lindeberg}\lim_{n\to \infty} \sum_{j=1}^n \frac{1}{s_n^2} 
 \int_{|t|\geq \epsilon s_n} t^2 f_j(t) dt < \infty,\;\;\;\;\;\;\;
\text{(condici\'on de Lindeberg)}\end{equation}
donde $f_j$ es la funci\'on de densidad de $X_j$. Entonces
$S_n/s_n$ tiende en distribuci\'on a
la normal $\frac{1}{\sqrt{2\pi}} e^{-t^2/2}$ cuando $n\to \infty$.
\end{thm}
Si $X_j$ es discreta, entonces, claro est\'a, la condici\'on
de Lindeberg se escribe
\[\lim_{n\to \infty} \sum_{j=1}^n \frac{1}{s_n^2} 
 \sum_{x:|x|\geq \epsilon s_n} x^2 \Prob(X_j = x) < \infty .\]
\begin{proof}[Esbozo de una prueba]
Se procede como en la primera demostraci\'on que dimos del teorema del
l\'imite central. La condici\'on de Lindeberg sirve para mostrar que
\[\lim_{n\to \infty} \sum_{j=1}^n 
\left|\widehat{X_j/s_n}(t) - \left(1 - \frac{1}{2} t^2 \Var(X_j)/s_n^2\right)\right| = 0\]
para todo $t$. Esto nos permite mostrar que
\[\begin{aligned}
\widehat{S_n/s_n} &= \prod_{j\leq n} 
\left(1 - \frac{1}{2} t^2 \Var(X_j)/s_n^2\right)^n + o(1)\\
&= \prod_{j\leq n} e^{- \frac{1}{2} t^2 \Var(X_j)/s_n^2} + o(1) = e^{-t^2/2} + o(1),
\end{aligned}\]
que es lo que deseamos.
\end{proof}
Para una prueba completa, ver \cite[XV.6, Teorema 1]{Fe}. La
condici\'on de Lindeberg (\ref{eq:lindeberg}) es b\'asicamente
necesaria (\cite[XV.6, Teorema 2]{Fe}).\index{condici\'on de Lindeberg}

\item\label{prob:astra} Sea 
\[C_p = \begin{cases} 1 &\text{si $p^2|n$}\\ 0 &\text{si $p^2\nmid n$,}
\end{cases}\]
donde $n$ es tomado al azar entre $1$ y $N$. Sea $C = \sum_{p\leq N} C_p$.
(Como, para $p>\sqrt{N}$, no hay entero $n\leq N$ tal que $p^2|n$,
 tenemos que $C = \sum_p C_p = \sum_{p\leq \sqrt{N}} C_p$.)  En otras palabras,
$C$ es la variable aleatoria que da el n\'umero de cuadrados de primos que 
dividen un entero tomado al azar entre $1$ y $N$. Estudiemos la distribuci\'on
de $C$.
\begin{enumerate} 
\item Muestre que
\[\mathbb{E}(C) = \sum_{p\leq \sqrt{N}} C_p = \sum_{p\leq \sqrt{N}}
\frac{1}{p^2} + O\left(N^{-1/2}\right) 
= \sum_p \frac{1}{p^2} +
O\left(N^{-1/2}\right).\]
%(La serie infinita $\sum_p \frac{1}{p^2}$ converge.)

\item Muestre que
\[\Prob(C = 0) = \prod_p \left(1 - \frac{1}{p^2}\right) + 
O\left(N^{-1/2}\right) .
\]
Por consiguiente, $\lim_{N\to \infty} \Prob(C=0)$ existe y es igual a
$\prod_p \left(1 - \frac{1}{p^2}\right)$. El evento $C = 0$ no es sino
el evento que $n$ carezca de divisores cuadrados
(i.e. $d^2\nmid n$ para todo entero $d>1$).

\item Muestre que, para todo $k\geq 0$,
\begin{equation}\label{eq:nacht}\begin{aligned}
\Prob(C=k) &= 
\mathop{\sideset{}{^*}\sum_{1\leq m\leq \sqrt{N}}}_{\omega(m) = k} 
 \frac{1}{m^2} \cdot \left(\prod_{p\nmid m} \left(1 - \frac{1}{p^2}\right) +
O\left((N/m)^{-1/2}\right) \right) \\ &=
\mathop{\sideset{}{^*}\sum_{m\geq 1}
}_{\omega(m) = k} \frac{1}{m^2}
\cdot \prod_{p\nmid m} \left(1 - \frac{1}{p^2}\right) +
O(N^{-1/2}) ,
\end{aligned}\end{equation}
donde $\sum^*$ denota una suma s\'olo sobre enteros sin divisores cuadrados.
Por lo tanto, $\lim_{N\to \infty} \Prob(C=k)$ existe y es igual a
\[\mathop{\sideset{}{^*}\sum_{m\geq 1}}_{\omega(m) = k} \frac{1}{m^2}
\cdot \prod_{p\nmid m} \left(1 - \frac{1}{p^2}\right).\]
Tenemos, entonces, que $C$ converge a una distribuci\'on discreta cuando
$N\to \infty$. Esta distribuci\'on no es la normal (ya que es discreta)
ni se le parece (no es sim\'etrica alrededor de $\mathbb{E}(C)$: 
la probabilidad de $C<0$ es cero, pero la probabilidad de 
$C>2 \mathbb{E}(C)$ tiende a un valor positivo). Lo crucial aqu\'i es que
$\lim_{N\to \infty} \mathbb{E}(C) < \infty$, es decir, el hecho que
$\sum_p \mathbb{E}(X_p) = \sum_{p} \frac{1}{p^2}$ converge. 
En cambio, cuando examin\'abamos
$X = \sum_p X_p$, ten\'iamos que, como $\sum_p \frac{1}{p}$ diverge,
la esperanza $\mathbb{E}(X)$ tend\'ia a $\infty$ cuando $N\to \infty$
y la distribuci\'on l\'imite era la normal.\index{distribuci\'on!normal}

\item\label{it:pois} Muestre que, para todo $k$,
\[\Prob(C = k) \ll \frac{\lambda^k}{k!}\]
donde $\lambda = \sum_p \frac{1}{p^2}$. La distribuci\'on
\begin{equation}\label{eq:gewitter}
\frac{\lambda^k}{k!} e^{-\lambda}\end{equation}
es la famosa {\em distribuci\'on de Poisson};\index{distribuci\'on!de Poisson}
 se trata del l\'imite
$n\to \infty$ de la distribuci\'on de $Y_n = Y_{n,1} + \dotsc + Y_{n,n}$, donde
$\{Y_{i,j}\}_{i,j\geq 1}$ son variables mutuamente independientes con
la distribuci\'on
\[
Y_{n,j} = \begin{cases}
1 &\text{con probabilidad $\lambda/n$}\\
0 &\text{con probabilidad $1 - \lambda/n$.}
\end{cases} 
\]
(Demuestre que $Y_n$ tiende, en efecto, a (\ref{eq:gewitter}).)

Hemos, entonces, acotado la distribuci\'on de $C$ por una distribuci\'on
de Poisson de esperanza $\lambda$. \'Esta no es una cota ``ajustada'':
pruebe que
\[\Prob(C = k) \ll_{\epsilon} \frac{\epsilon^k}{k!}\]
para todo $\epsilon>0$, y, en consecuencia,
\[\Prob(C = k) = o\left(\frac{\lambda^k}{k!}\right)
\]
cuando $k\to \infty$, e incluso
\[\Prob(C\geq k) = o\left(\frac{1}{k!}\right).\]

\item\label{it:ogor}
%Veamos otra variable con la misma distribuci\'on que $C$;
%le daremos una aplicaci\'on m\'as tarde. 
Sea
\[D_p = \begin{cases} 1 &\text{si $p|m$ y $p|n$}\\
0 &\text{de lo contrario,}\end{cases}\]
donde $(m,n)$ es un par de n\'umeros entre $1$ y $N$ tomados al
azar de acuerdo a la distribuci\'on uniforme. Sea $D = \sum_p D_p$.
Entonces el evento $D = 0$ no es sino el evento que $m$ y $n$ sean
primos entre si.

Muestre que \[\begin{aligned}
\Prob(D = 0) &= \prod_p \left(1 - \frac{1}{p^2}\right) + 
O(N^{-1}),\\
\Prob(D = k) &= 
\mathop{\sideset{}{^*}\sum_{m\geq 1}}_{\omega(m) = k} 
\frac{1}{m^2} \cdot
   \prod_{p\nmid m} \left(1 - \frac{1}{p^2}\right) + O(N^{-1}) .
\end{aligned}\]
Compare esto con (\ref{eq:nacht}). Como en (\ref{it:pois}), concluya que
\[\Prob(D\geq k) = o\left(\frac{1}{k!}\right) .\]
\end{enumerate}
%\item 
%Poisson; cola de $\rho_p$ 
%(dejar para su propia secci\'on)

%\item\label{nota:erdwin}
%Erd\"os-Wintner
%esperanza finita
\end{enumerate}
%los l\'imites no centrales?

\section{Grandes desviaciones: cotas superiores.
Valores cr\'iticos.}\label{sec:bambi}

Sean $X_1, X_2, \dotsc$ variables aleatorias mutuamente independientes
 con la distribuci\'on
\begin{equation}\label{eq:popeye} X_j = \begin{cases}
0 &\text{con probabilidad $1/2$}\\
1 &\text{con probabilidad $1/2$.}
\end{cases}\end{equation}
Sea $X = \sum_{j\leq n} X_j$. Sabemos que $\mathbb{E}(X) = \frac{1}{2} n$ y
que la distribuci\'on de $X$ ser\'a cercana a la normal alrededor de 
$\frac{1}{2} n$. Qu\'e pasa lejos de $\frac{1}{2} n$? 

En general, hablamos de peque\~nas desviaciones cuando la distancia 
entre el valor de $X$ y la esperanza $\mathbb{E}(X)$ es $O(\sqrt{\Var(X)})$,
y de {\em grandes desviaciones} cuando la distancia es comparable a\index{desviaciones!grandes}
$\Var(X)$ o $\mathbb{E}(X)$. Podemos preguntarnos, por ejemplo,
que tan a menudo se dan las grandes desviaciones
$\omega(n) < \frac{1}{2} \log \log n$ o $\omega(n) > 6 \log \log n$.

En el caso de las variables (\ref{eq:popeye}), podemos hacer los 
c\'alculos a mano. Tenemos $\Prob(X = m) = 2^{-n} \cdot \binom{n}{m}$
para todo $m$. Por lo tanto,
\[\Prob(Y > a n) = 2^{-n} \sum_{m\geq a n} \binom{n}{m}\]
para todo $a$, donde
\[\binom{n}{m} := \frac{n!}{(n-m)! m!} =
\frac{1 \cdot 2 \cdot 3 \cdot \dotsb n}{(1\cdot 2\cdot  \dotsb (n-m)) 
\cdot (1\cdot 2\cdot \dotsb m)} = 
\frac{n (n-1) \dotsb (n-m+1)}{1\cdot 2\cdot
 \dotsb m}\]
es el n\'umero de maneras de escoger $m$ cosas de entre $n$ cosas.
(Tenemos $n$ posibilidades para la primera cosa elegida, 
$n-1$ posibilidades para la segunda, \dots, $n-m+1$ posibilidades
para la $m$--\'esima, y no importa en qu\'e orden de los $m!$ \'ordenes
posibles hayamos elegido las $m$ cosas. Por lo tanto, hay
$\frac{n (n-1) \dotsb (n-m+1)}{1\cdot 2\cdot
 \dotsb m}$ maneras de escoger.)

Fijemos un $a\in \lbrack 1/2,1\rbrack$. Como
$m\mapsto \binom{n}{m}$ es decreciente para $m\geq \frac{1}{2} n$, tenemos
\begin{equation}\label{eq:ken1}2^{-n} \binom{n}{\lceil a n\rceil} \leq
\Prob(Y > a n) \leq (n+1) 2^{-n} \binom{n}{\lceil a n \rceil} .\end{equation}
Ahora bien,
\[
\log n! = \log 1 + \log 2 + \dotsb + \log n =
\int_1^n \log x\, dx + O(\log n) = n \log n - n + 1 + O(\log n),\]
as\'i que
\begin{equation}\label{eq:ken2}\begin{aligned}
\log \binom{n}{m} &= n \log n - ((n-m) \log (n-m) + m \log m) 
+ O(\log n)\\ &= - n \cdot \left(\left(1 - \frac{m}{n}\right) \cdot
 \log \left(1 - \frac{m}{n}\right) +
\frac{m}{n} \log \frac{m}{n}\right) + O(\log n). \end{aligned}\end{equation}
Por (\ref{eq:ken1}) y (\ref{eq:ken2}), para $a\geq \frac{1}{2}$,
\[\begin{aligned}
\log(\Prob(Y> a n)) &\sim \log \frac{2^{-n}}{(a^a (1-a)^{(1-a)})^n}
+ O(\log n)\\ &= n ( -a \log a - (1 - a) \log (1-a) - \log 2) + O(\log n) .
\end{aligned}\]
En otras palabras, $\Prob(Y> a n) = e^{(H + o(1)) n}$, donde
\begin{equation}\label{eq:entro}
H = - a \log \frac{a}{0.5} - (1-a) \log \frac{1-a}{0.5} .\end{equation}
La forma de (\ref{eq:entro}) puede resultar familiar para los f\'isicos
(o qu\'imicos). Elaboraremos esta observaci\'on para m\'as adelante;
ahora, pasemos a estudiar $\omega(n)$.
\begin{center}
* * *
\end{center}
Sean $X'_2, X'_3, X'_5,\dotsc$ variables mutuamente independientes con la 
distribuci\'on
\begin{equation}\label{eq:ovlid}
X'_p = \begin{cases}
1 &\text{con probabilidad $\frac{1}{p}$}\\
0 &\text{con probabilidad $ 1 - \frac{1}{p}$.}
\end{cases}
\end{equation}
Sea $X' = \sum_{p\leq N} X'_p$. La esperanza $\mathbb{E}(X')$ es
$\log \log N + O(1)$.
Nos preguntaremos cu\'anto es la probabilidad $\Prob(X' < a \log \log N)$,
$a<1$, o $\Prob(X'>a \log \log N)$, $a>1$.

Markov\index{desigualdad!de Markov} acota las cotas mediante $\mathbb{E}(X')$, Chebyshev\index{desigualdad!de Chebyshev} mediante
$\mathbb{E}(X'^2)$; muy bien podemos usar $\mathbb{E}(X'^k)$, o incluso
$\mathbb{E}(e^{X'})$. En efecto, para $a$ positivo y $\alpha\geq 0$,
\begin{equation}\label{eq:auto}
\Prob(X'> a \log \log N) \leq \frac{\mathbb{E}\left(e^{\alpha X'}\right)}{
e^{\alpha \cdot a \log \log N}},\end{equation}
y, para $a$ positivo y $\alpha\leq 0$,
\begin{equation}\label{eq:wash}
\Prob(X' < a \log \log N) \leq \frac{\mathbb{E}\left(e^{\alpha X'}\right)}{
e^{\alpha \cdot a \log \log N}}.\end{equation}
As\'i como la desigualdad de Chebyshev es simplemente la desigualdad
de Markov aplicada a la variable $X'^2$, las desigualdades (\ref{eq:auto})
y (\ref{eq:wash}) son simplemente la desigualdad de Markov
aplicada a la variable $e^{\alpha X'}$; tenemos la libertad de escoger
el par\'ametro $\alpha$ como m\'as nos convenga. (La aplicaci\'on
de la desigualdad de Markov a una variable de la forma $e^{\alpha X'}$
es a veces denominada {\em el m\'etodo de momentos exponenciales};\index{momentos exponenciales}
los momentos ``usuales'' son $\mathbb{E}(X')$, $\mathbb{E}(X'^2)$,
$\mathbb{E}(X'^3)$,\dots)

C\'omo evaluar $\mathbb{E}(e^{\alpha X'})$? Tenemos
\[e^{\alpha X'} = e^{\alpha \sum_{p\leq N} X'_p} = \prod_{p\leq N} e^{\alpha X'_p}
  .\]
Como $X'_2, X'_3, X'_5,\dotsc$ son mutuamente independientes,
las variables $e^{\alpha X'_2}, e^{\alpha X'_3}, e^{\alpha X'_5},\dotsc$
tambi\'en lo son, y por lo tanto
\[\mathbb{E}(e^{\alpha X'}) = \mathbb{E}\left(\prod_{p\leq N} e^{\alpha
    X'_p}\right) = \prod_{p\leq N} \mathbb{E}\left(e^{\alpha X'_p}\right) .
\]
Es f\'acil calcular que
\[\mathbb{E}\left(e^{\alpha X'_p}\right) = 1 + \frac{1}{p} (e^{\alpha} - 1) .\]
En consecuencia,
\begin{equation}\label{eq:manuc}
\mathbb{E}\left(e^{\alpha X'}\right) = \prod_{p\leq N}
(1 + \frac{1}{p} (e^{\alpha} - 1)) \ll \prod_{p\leq N} \left(1 + \frac{1}{p}
\right)^{e^{\alpha} - 1}.\end{equation}
(Esto es cierto a\'un si $\alpha < 0$, o incluso
$\alpha = - \infty$. Demuestre la desigualdad en (\ref{eq:manuc})
si \'esta no lo convence de inmediato; use el hecho que 
$\prod_{p\leq N} (1 + 1/p^2) \ll 1$, puesto que
$\prod_p (1 + 1/p^2)$
converge (por qu\'e?).) 
Ahora bien,
\begin{equation}\label{eq:casbah}
e^{\sum_{p\leq N} \frac{1}{p}} \ll \prod_{p\leq N} \left(1 +
  \frac{1}{p}\right) \ll e^{\sum_{p\leq N} \frac{1}{p}} ,
\end{equation}
y, como $\sum_{p\leq N} \frac{1}{p} = \log \log N + O(1)$ (Chebyshev-Mertens
(\ref{eq:cruj})), obtenemos que
\begin{equation}\label{eq:malmal}
\mathbb{E}\left(e^{\alpha X'}\right) \ll_{\alpha} \left(e^{\log \log N}\right)^{
e^{\alpha} - 1} = (\log N)^{e^{\alpha} - 1} .
\end{equation}

Podemos ahora sustituir (\ref{eq:malmal}) en (\ref{eq:auto}) y 
(\ref{eq:wash}). Nos queda escojer el valor \'optimo de $\alpha$. Para 
$a>1$ y $\alpha>0$,
\[\Prob(X' > a \log \log N) \ll_{\alpha} \frac{(\log N)^{e^{\alpha} - 1}}{e^{\alpha
\cdot a \log \log N}} = (\log N)^{e^{\alpha} - 1 - \alpha a} ;\]
para $a<1$ y $\alpha<0$,
\[\Prob(X' < a \log \log N) \ll_{\alpha} \frac{(\log N)^{e^{\alpha} - 1}}{e^{\alpha
\cdot a \log \log N}} = (\log N)^{e^{\alpha} - 1 - \alpha a} .\]
Debemos minimizar $e^{\alpha} - 1 - \alpha a$ para $a$ positivo y fijo,
cuid\'andonos que $\alpha$ sea negativo si $a>1$, y positivo si $a<1$.
Sacando derivadas, vemos que $e^{\alpha} - 1 - \alpha a$ es m\'inimo
para $a$ dado cuando $\alpha = \log a$. Entonces
\begin{equation}\label{eq:zarzu}
\begin{aligned}
\Prob(X'> a \log \log N) &\ll_a (\log N)^{- (a \log a + 1 - a)}
&\text{para $a>1$,}\\
\Prob(X'<a \log \log N) &\ll_a (\log N)^{- (a \log a + 1 - a)}
&\text{para $a<1$.}
\end{aligned}
\end{equation}
Hemos obtenido las cotas superiores
que dese\'abamos para las grandes desviaciones
de $X' = \sum_{p\leq N} X'_p$, donde $X'_p$ es como en (\ref{eq:ovlid}).
Las constantes impl\'icitas en (\ref{eq:zarzu}) dependen de $a$, pero
de manera continua;
por lo tanto, la misma constante puede servir para todo $a\leq A$,
$A$ dado. (La constantes impl\'icitas en (\ref{eq:zarzu}) dependen
de manera continua
de $e^{\alpha} - 1 = e^{\log a} - 1 = a - 1$, por lo cual
la dependencia en $a$ es continua a\'un cerca de $a = 0$.)

\begin{center}
* * *
\end{center}

Como de costumbre, pasamos a examinar las variables $X_2, X_3, X_5,\dotsc$
dadas por
\[X_p = \begin{cases} 1 &\text{si $p|n$,}\\ 0 &\text{si $p\nmid n$,}
\end{cases}\]
donde $n$ es tomado al azar entre $1$ y $N$ con la distribuci\'on uniforme.
Sea $X = \omega(n) = \sum_{p\leq N} X_p$. Las desigualdades
(\ref{eq:auto}) y (\ref{eq:wash}) son a\'un v\'alidas; partamos de ellas.

El problema es evaluar $\mathbb{E}(e^{\alpha X}) = 
\mathbb{E}(e^{\alpha \omega(n)})$. Si bien $e^{\alpha X} =
e^{\alpha \sum_{p\leq N} X_p} = \prod_{p\leq N} e^{\alpha X_p}$, 
no podemos
concluir que $\mathbb{E}(e^{\alpha X}) =
\mathbb{E}\left(e^{\alpha \sum_{p\leq N} X_p}\right) = \prod_{p\leq N} 
\mathbb{E}\left(e^{\alpha X_p}\right)$, ya que las variables $X_p$ no son
mutuamente independientes. Empero, podemos mostrar sin mucha dificultad que
\[\sum_{n\leq N} \frac{e^{\alpha \omega(n)}}{n} \leq \prod_{p\leq N}
\left(1 + e^{\alpha} \left(\frac{1}{p} + \frac{1}{p^2} + \dotsc\right)\right)
\ll \prod_{p\leq N} \left(1 + \frac{1}{p}\right)^{e^{\alpha}} \ll
(\log N)^{e^{\alpha}}
\] 
(por (\ref{eq:casbah})) y deducir con m\'as trabajo de esto que
\[\sum_{n\leq N} e^{\alpha \omega(n)} \ll N (\log N)^{e^{\alpha} - 1}
\]
(ver problema \ref{not:antartor}). Por lo tanto, $\mathbb{E}(e^{\alpha
\omega(n)}) \ll (\log N)^{e^{\alpha} - 1}$. Fijamos $\alpha = \log a$ como
antes, y concluimos (por (\ref{eq:auto}) y (\ref{eq:wash})) que
\begin{equation}\label{eq:vclib}\begin{aligned}
\Prob(\omega(n) > a \log \log N) &\ll_a (\log N)^{- (a \log a + 1 - a)}
&\text{\;\; si $a>1$,}\\
\Prob(\omega(n) <a \log \log N) &\ll_a (\log N)^{- (a \log a + 1 - a)}
&\text{\;\; si $a<1$.}
\end{aligned}\end{equation}
Las constantes impl\'icitas en (\ref{eq:vclib}) dependen de $a$ de manera
continua -- a\'un en la vecindad de $a = 0$.
\begin{center}
* * *
\end{center}

Tanto (\ref{eq:zarzu}) como (\ref{eq:vclib}) son cotas superiores. 
Encontraremos dentro de poco cotas inferiores muy cercanas a estas
cotas superiores; mientras tanto, content\'emonos con una aplicaci\'on de
(\ref{eq:vclib}).

Cu\'antos enteros $n$ entre $1$ y $N^2$ pueden expresarse como el producto
$n = a \cdot b$ de dos enteros $a$, $b$ entre $1$ y $N$? En otras palabras:
si escribimos una tabla de multiplicaci\'on de $1000$ por $1000$ (digamos),
habr\'a un mill\'on de n\'umeros en la tabla,
y todos estos n\'umeros estar\'an entre $1$ y un mill\'on; empero, como
hay muchas repeticiones, podemos preguntarnos: cu\'antos de los n\'umeros
entre $1$ y un mill\'on est\'an presentes en la tabla?

Este es el conocido {\em problema de la tabla de multiplicaci\'on}.\index{problema de la tabla de multiplicaci\'on}
Lo encontramos por primera vez en \S \ref{sec:varianza},
ejercicio \ref{it:gorgut}, donde probamos que el n\'umero de enteros
$n\leq N^2$ que aparecen en la tabla (i.e., $n = a\cdot b$
para alg\'un par $1\leq a,b\leq N$) es $o(N^2)$. Ahora 
estableceremos una cota superior bastante m\'as ajustada.

Comenzamos de la misma manera que antes: tomamos conciencia de que, si bien
la esperanza de $\omega(n)$ (para $n$ tomado al azar) es 
\[\mathbb{E}(\omega(n)) = \log \log N^2
+ O(1) = \log \log N + O(1),\]
la esperanza de $\omega(a\cdot b)$ (para $a$, $b$ tomados al azar) es
\[\mathbb{E}(\omega(a \cdot b)) = \mathbb{E}(\omega(a) + 
\omega(b) - \omega(\mcd(a,b))) = \mathbb{E}(\omega(a)) + \mathbb{E}(\omega(b))
- O(1) = 2 \log \log N + O(1) .\]
Por lo tanto, el n\'umero de divisores $\omega(n)$
de todo n\'umero $n = a \cdot b$ ($n\leq N^2$, $a,b\leq N$) debe
estar ya sea en la cola de la distribuci\'on de $\omega(n)$ o en la
cola de la distribuci\'on de $\omega(a \cdot b)$. Ahora que sabemos como
acotar las colas distantes (es decir, las grandes desviaciones) de 
$\omega(n) = X = \sum_p X_p$, podremos dar una buena cota superior
para la (poca) probabilidad de $\omega(n)$ dentro de una u otra distribuci\'on.

Estamos en una situaci\'on que ya es propia no s\'olo de las probabilidades,
sino de la estad\'istica. La situaci\'on es as\'i: hay una variable
aleatoria observable $U$; su distribuci\'on se desconoce, pero se tienen
dos conjeturas acerca de \'esta. Hablemos, entonces, de la distribuci\'on
$1$ y la distribuci\'on $2$. La esperanza $\mathbb{E}_1(U)$ de $U$ seg\'un
la distribuci\'on $1$ es (digamos) $u_1$, mientras que la
esperanza $\mathbb{E}_2(U)$ de $U$ seg\'un la distribuci\'on $2$ es 
$u_2>u_1$. La tarea es determinar cu\'al de las dos distribuciones es m\'as
veros\'imil.

Debemos ponernos de acuerdo en un {\em valor cr\'itico} $t$ entre
$u_1$ y $u_2$.\index{valor cr\'itico}
Si, despues de $n$ mediciones de $U$, vemos que el promedio
$X = \frac{1}{n} (U_1 + U_2 + \dotsc + U_n)$ es menor que $t$, decidiremos
la disputa en favor de la distribuci\'on $1$; si $X\geq t$, daremos la raz\'on
a $2$. La pregunta es: que {\em valor cr\'itico} $t$ debemos escoger?

Si la distribuci\'on $1$ es la verdadera distribuci\'on de $U$, la
probabilidad de decidir la disputa err\'oneamente a favor de la distribuci\'on
$2$ es $\Prob_1(X \geq t)$ (donde $\Prob_1(\dotsc)$ denota la probabilidad
de un evento si se asume que la distribuci\'on $1$ es la cierta).
 Si la 
distribuci\'on $2$ es la cierta, la probabilidad de errar a favor de la
distribuci\'on $1$ es $\Prob_2(X<t)$ (donde $\Prob_2(\dotsc)$ denota la probabilidad
de un evento si se asume que la distribuci\'on $2$ es la cierta). Resulta sensato, entonces, escoger
el valor de $t$ entre $u_1$ y $u_2$ para el cual
\[\Prob_1(X\geq t)+\Prob_2(X<t)\]
sea m\'inimo.

En el caso del problema de la tabla de multiplicaci\'on, estamos ante una
situaci\'on parecida. Fijaremos un $t$ entre $1$ y $2$. Para todo
$n = a\cdot b$, $a,b\leq N$, tendremos ya sea $\omega(n)<t \log \log N$ o
$\omega(n)\geq t \log \log N$. En el primer caso,
cualquier par $a,b\leq N$ tal que $n = a\cdot b$ tendr\'a que estar
en el conjunto $\{(a,b): 1\leq a,b\leq N,\; \omega(a\cdot b)<t \log \log N\}$,
y el n\'umero de elementos de este conjunto es $N^2$ multiplicado por
\[\Prob(\omega(a\cdot b) < t \log \log N),\]
donde $a$ y $b$ son tomados al azar entre $1$ y $N$ con la distribuci\'on
uniforme. En el segundo caso -- es decir, $\omega(n) \geq t \log \log N$ --
el n\'umero $n$ estar\'a en el conjunto $\{n:1\leq n\leq N,
\omega(n)\geq t\log \log N\}$,
cuyo n\'umero de elementos es
$N^2$ multiplicado por
\[\Prob(\omega(n) \geq t \log \log N),\]
donde $n$ es tomado al azar entre $1$ y $N^2$ con la distribuci\'on uniforme.
Por lo tanto, el conjunto de enteros $n\leq N^2$ tales que $n = a \cdot b$
para alg\'un par $a,b\leq N$ tiene a lo m\'as
\begin{equation}\label{eq:valor}N^2 \cdot (\Prob(\omega(n) \geq t \log \log N) +
\Prob(\omega(a\cdot b) < t \log \log N))\end{equation}
elementos. Queremos, entonces, escoger $t$ tal que (\ref{eq:valor}) sea
tan peque\~no como se pueda.

Para hacer esto, debemos primero estimar (\ref{eq:valor}).
Gracias a nuestra cota de grandes desviaciones (\ref{eq:vclib}),
\begin{equation}\label{eq:antij}\begin{aligned}
\Prob(\omega(n) \geq t \log \log N) &= \Prob(\omega(n)\geq t \log \log N^2 -
t \log 2)\\
&= \Prob(\omega(n) \geq t (1 - o(1)) \cdot \log \log N^2)\\
&\leq 
 (\log N)^{- I(t) (1 + o(1))} ,\end{aligned}\end{equation}
donde $I(t) = t \log t + 1 - t$.
Queda acotar $\Prob(\omega(a \cdot b) < t \log \log N)$. Est\'a claro que
\[\omega(a \cdot b) = \omega(a) + \omega(b) - \omega(\mcd(a,b)) .\]
Por \S \ref{sec:limcen}, ejercicio \ref{it:ogor},
 la probabilidad
que $\mcd(a,b)$ sea $> \epsilon \log \log N$ es $\ll (\log N)^{- 100}$ (o
cualquier otra potencia que se quiera).
(Hay varias maneras simples de evitar el uso del ejercicio \ref{it:ogor};
la ventaja de este \'ultimo es que nos permite continuar sin salir nunca
de un marco probabil\'istico.)
 Bastar\'a, entonces,
con estimar la probabilidad que $(\omega(a),\omega(b))$ quede en
un cuadrado $\mathbf{C}_{t_a,t_b}$ de lado $\epsilon \log \log N$ que contenga 
el punto $(t_a \log \log N, t_b \log \log N)$, 
donde $t_a + t_b\leq t + \epsilon$. El m\'aximo
de tal probabilidad, multiplicado por el n\'umero de cuadrados
a considerar (el cual es $O(1/\epsilon^2) = O_{\epsilon}(1)$), nos
dar\'a una cota superior para la probabilidad de $\omega(a \cdot b)
\leq t$.

Como $a$ y $b$ son variables independientes, la cota de grandes
desviaciones (\ref{eq:vclib}) nos dice que la probabilidad que
$(\omega(a),\omega(b))$ est\'e dentro de $\mathbf{C}_{t_a,t_b}$ es
\[\ll (\log N)^{- (I(t_a) + I(t_b)) + O(\epsilon)},\]
donde $I(t) = t \log t + 1 - t$.
Nos estamos preguntando, entonces, cu\'al es el m\'aximo de
$I(t_a) + I(t_b)$, bajo la condici\'on que $t_a + t_b \leq t$.
Sacamos la segunda derivada de $I(x)$ y vemos que siempre es positiva.
Por lo tanto, la gr\'afica de $x\mapsto I(x)$ se curva hacia arriba,
y, as\'i, $\frac{1}{2} (I(t_a) + I(t_b)) \geq I\left(
\frac{t_a + t_b}{2}\right)$. Como $t$ est\'a entre $1$ y $2$, y $I(x)$
es decreciente en el intervalo entre $1/2$ y $2/2 = 1$, vemos que
$I\left(\frac{t_a + t_b}{2}\right) \geq I(t/2 + \epsilon) 
= I(t/2) + O(\epsilon)$.
Concluimos que
\[(\log N)^{ - (I(t_a) + I(t_b)) + O(\epsilon)} \leq
(\log N)^{- 2 I(t/2) + O(\epsilon)} .\]
Por lo tanto, 
\begin{equation}\label{eq:pno}\Prob(\omega(a\cdot b) \leq t \log \log N) \leq
(\log N)^{- 2 I(t/2) + O(\epsilon)} .\end{equation}

Comparando (\ref{eq:antij}) y (\ref{eq:pno}), vemos que tenemos
que escoger $t\in \lbrack 1,2\rbrack$ de tal manera que
\[\min(I(t),2 I(t/2))\]
sea m\'aximo. Ahora bien, $I(t)$ es creciente y $2 I(t/2)$ es decreciente
en el intervalo $\lbrack 1,2\rbrack$. Por lo tanto, el m\'aximo de
$\min(I(t), 2 I(t/2))$ se alcanza cuando $I(t) = 2 I(t/2)$. Resolvemos
$I(t) = 2 I(t/2)$, y encontramos $t = \frac{1}{\log 2}$. Hemos
probado el siguiente resultado:

\begin{thm}[Erd\H{o}s \cite{E}]\index{Erd\H{o}s, P.}\index{problema de la tabla de multiplicaci\'on}
El n\'umero de enteros $n$ entre $1$ y $N^2$ que pueden ser expresados
como el producto $n = a \cdot b$ de dos enteros $a$, $b$ entre $1$ y $N$
es
\[O_{\epsilon}\left(N^2 \cdot (\log N)^{-I\left(\frac{1}{\log 2}\right) + \epsilon}\right)\]
para cualquier $\epsilon>0$, donde $I(t) = t \log t + 1 - t$.
\end{thm}
Num\'ericamente, $I\left(\frac{1}{\log 2}\right) = 0.08607\dotsc$.
\begin{center}
{\bf Notas y problemas}
\end{center}
\begin{enumerate}
\item\label{not:antartor}
 Debemos estimar $\sum_{n\leq N} e^{\alpha \omega(n)}$, donde $\alpha>0$.
Si se usan m\'etodos anal\'iticos -- que no trataremos para evitar el 
an\'alisis complejo -- esto es rutina (dentro del m\'etodo de
Selberg--Delange -- ver, e.g., \cite[II.5]{T}). Veremos como obtener
de manera elemental una cota superior ``del orden correcto'' 
(es decir, $\ll$ la asint\'otica verdadera).
\begin{enumerate}
\item\label{it:ino} Primero deduzca de (\ref{eq:ort}) que
\[\prod_{p\leq N} \left(1 + \frac{1}{p}\right) \ll \log N,\]
\begin{equation}\label{eq:canada}\prod_{p\leq N} \left(1 + \frac{1}{p} + \frac{1}{p^2} + \dotsc\right) 
\ll \log N .\end{equation}

\item\label{it:cencia} Acotemos $\sum_{n\leq N} e^{\alpha \omega(n)}/n$:
\begin{equation}\label{eq:put}\begin{aligned}
\sum_{n\leq N} \frac{e^{\alpha \omega(n)}}{n} &\leq
\prod_{p\leq N} \left(1 + e^{\alpha} \left(\frac{1}{p} + \frac{1}{p^2} +
\frac{1}{p^3} + \dotsc \right) \right)\\ &\leq
\prod_{p\leq N} \left(1 + \frac{1}{p} + \frac{1}{p^2} + \dotsc\right)^{
e^{\alpha}} \ll_{\alpha} (\log N)^{e^{\alpha}} .\end{aligned}\end{equation}
\item\label{it:irmin} Podr\'iamos concluir que
 $\sum_{n\leq N} e^{\alpha \omega(n)} \ll N \sum_{n\leq N} \frac{e^{\alpha
     \omega(n)}}{n} \ll N (\log N)^{e^{\alpha}}$, lo cual es correcto, pero
 esta cota se aleja de la realidad por un factor de $(\log N)$
(como veremos despu\'es). C\'omo 
obtendremos una mejor cota?\\

\begin{small}
(Como dijimos, es rutina obtener la ``suma parcial'' $\sum_{n\leq N} 
e^{\alpha \omega(n)}$ a trav\'es del an\'alisis (no elemental) de la funci\'on
$L(s) = \sum_{n=1}^{\infty} e^{\alpha \omega(n)} n^{-s}$. Estamos tratando de
obtener la suma parcial $\sum_{n\leq N} e^{\alpha \omega(n)}$ meramente
de la suma parcial $\sum_{n\leq N} e^{\alpha \omega(n)} n^{-1}$, y quiz\'as
de un par de propiedades de $e^{\alpha \omega(n)}$.)\\
\end{small}

La {\em convoluci\'on} $f\ast g$ de dos funciones\index{convoluci\'on}
 $f,g:\mathbb{Z}^+\to
\mathbb{C}$ es $(f\ast g)(n) = \sum_{d|n} f(d) g(n/d)$. A menudo, en la
teor\'ia anal\'itica de n\'umeros, es conveniente expresar una funci\'on
como la convoluci\'on de dos otras. Ya podemos obtener la suma parcial
de $h(n) = e^{\alpha \omega(n)} \ast {\bf 1}$, donde 
${\bf  1}:\mathbb{Z}^+\to \mathbb{C}$ es la funci\'on ${\bf 1}(n) =1$:
\begin{equation}\label{eq:mck}\begin{aligned}
\sum_{n\leq N} h(n) &= \sum_{n\leq N} \sum_{d|n} e^{\alpha \omega(d)} \cdot
1 = \sum_{d\leq N} e^{\alpha \omega(d)} \mathop{\sum_{m\leq N}}_{d|m} 1\\
&\leq \sum_{d\leq N} e^{\alpha \omega(d)} \frac{N}{d} = N \sum_{d\leq N}
\frac{e^{\alpha \omega(d)}}{d} \ll_{\alpha} N (\log N)^{e^{\alpha}} .\end{aligned}
\end{equation}
 Ahora bien, queremos conocer la suma parcial de $e^{\alpha \omega(n)}$, no la
de $e^{\alpha \omega(n)} \ast {\bf 1}$.

\item %En un mundo ideal, tendr\'iamos ${\bf 1} = {\bf 1} \ast {\bf 1}$;
%de esta manera podr\'iamos introducir convoluciones - y m\'as variables
%a manipular a nuestro agrado - cuando quisi\'eramos. Ahora bien, 
%${\bf 1} = {\bf 1} \ast {\bf 1}$ no es cierto, y tenemos que utilizar
%otras convoluciones que aproximen esta falsa identidad.

Utilicemos el hecho que $\log = {\bf 1} \ast \Lambda$. (Esto no es sino 
(\ref{eq:orgo}).) Tenemos que
\[\sum_{n\leq N} e^{\alpha \omega(n)} \log(n) = \sum_{n\leq N} \sum_{d|n} e^{\alpha
\omega(n)} \Lambda(n/d).\]
Ahora bien, $\omega(n)-1 \leq \omega(d)\leq \omega(n)$ para todo $n$ y todo
$d|n$ tales que $\Lambda(n/d)\ne 0$. As\'i,
\[\sum_{n\leq N} e^{\alpha \omega(n)} \log(n) \ll \sum_{n\leq N} \sum_{d|n} e^{\alpha
\omega(d)} \Lambda(n/d).\] 
De (\ref{eq:shotwo}) y (\ref{eq:mck}), se deduce que
\[\begin{aligned}
\sum_{n\leq N} \sum_{d|n} e^{\alpha
\omega(d)} \Lambda(n/d) &= \sum_{d\leq N} e^{\alpha \omega(d)} \sum_{m\leq N/d}
\Lambda(m)\\
&\ll \sum_{d\leq N}  e^{\alpha \omega(d)} \cdot N/d \ll_{\alpha}
 N (\log N)^{e^\alpha}.
\end{aligned}\]
Por lo tanto,
\[
\sum_{n\leq N} e^{\alpha \omega(n)} \log(n) \ll_{\alpha} N (\log N)^{e^{\alpha}}
.\]
Usando la t\'ecnica de la suma por partes (\S \ref{sec:esperanza}, 
nota (\ref{it:raul})), deduzca de esto que
\begin{equation}\label{eq:palomo}
\sum_{n\leq N} e^{\alpha \omega(n)} \ll_{\alpha} N (\log N)^{e^{\alpha} - 1} .
\end{equation}
La constante impl\'icita en $\ll$ depende de
$\alpha$, tanto aqu\'i como antes. La dependencia es continua
en $\alpha$ (o, si se prefiere, en $\beta = e^{\alpha}$; la dependencia
es continua a\'un en la vecindad de $\beta = 0$).
\end{enumerate}
\end{enumerate}
%Cu\'al es, por ejemplo,
%la probabilidad que $X> a n$ para $a>\frac{1}{2}$ fijo?
%y medianas

%tambi\'en: Sanov (entrop\'ia!)

%hypothesis testing!

\section{Grandes desviaciones: cotas inferiores.
Entrop\'ia.}\index{desviaciones!grandes}\index{entrop\'ia}

Consideremos dos compartimientos separados por una membrana porosa. 
Llen\'emoslos con un gas:
%\Begin{figure}
%\centering 
\begin{center}
\includegraphics[height=1in]{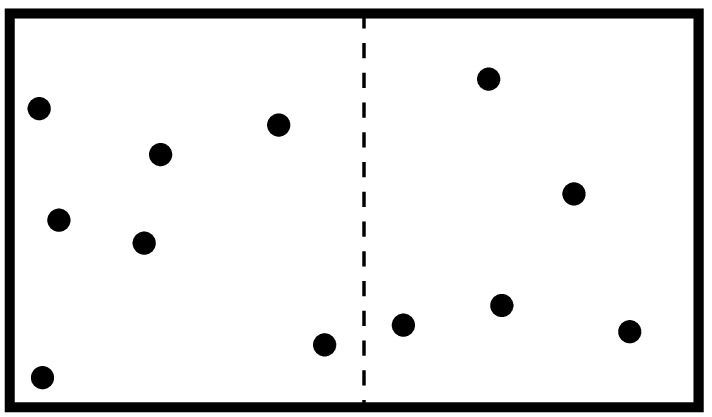}
\end{center}
De cu\'antas maneras puede ocurrir que haya $m$ part\'iculas en el
compartimiento de la izquierda y $n-m$ en el de la derecha? 

Hay $\binom{n}{m} = \frac{n!}{m! (n-m)!}$ maneras de escoger $m$
part\'iculas de entre $n$; por lo tanto, hay $\binom{n}{m}$ maneras
que haya $m$ part\'iculas en la izquierda y $n-m$ en la derecha.

Como $\log n! = \log 1 + \log 2 + \dotsc + \log n = n \log n - n + O(\log n)$,
\begin{equation}\label{eq:lmoor}\begin{aligned}
\log \binom{n}{m} &= \log n! - \log m! - \log (n-m)!\\
&= n \log n - m \log m - (n - m) \log (n-m) + O(\log n)\\
&= - n \cdot \left(\frac{m}{n} \log \frac{m}{n} + \frac{n-m}{n} \log
\frac{n-m}{n}\right) + O(\log n) .\end{aligned}\end{equation}

Qu\'e pasa si cada una de las $2^n$ maneras de colocar $n$ part\'iculas
en dos compartimientos es igualmente probable? (Hay $2^n$ maneras porque,
dada cada part\'icula, podemos elegir en cual de dos compartimientos puede
encontrarse.) Obtenemos de (\ref{eq:lmoor})
que la probabilidad que la proporci\'on de
part\'iculas en el comportamiento izquierdo sea $r = \frac{m}{n}$ es
\[\begin{aligned}
\frac{\log \binom{n}{m}}{2^n} &= e^{-n (r \log r + (1-r) \log(1-r) + \log 2)
+O(\log n)}\\ &= e^{n (H(r) + o(1))}. \end{aligned}\]
donde
\[\begin{aligned}
H(r) &= - (r \log r + (1 - r) \log (1- r) + \log 2)\\
 &= - \left(r \log \frac{r}{0.5} + (1 - r) \log \frac{1 - r}{0.5}\right) .
\end{aligned}\]
La cantidad $H(r)$ no es sino la famosa {\em entrop\'ia}.\index{entrop\'ia}
Comp\'arese con (\ref{eq:entro}).

Claro est\'a, no hay raz\'on f\'isica por la cual podamos asumir
que todas las disposiciones {\em iniciales} son igualmente probables;
muy bien podemos comenzar llenando s\'olo uno de los compartimientos
con gas. Empero, resulta poco sorprendente que, en la pr\'actica y 
con el paso del tiempo,
el sistema tienda a las proporciones $(r:1-r)$ que resultar\'ian 
m\'as probables si las $2^n$ disposiciones fueran igualmente probables.
Lo que hemos mostrado es que esto es lo mismo que decir que la entrop\'ia
$H(r)$ del sistema tiende a crecer. (Esta tendencia es muy clara, ya que
$H(r)$ est\'a en el exponente de $e^{n (H(r) + o(1))}$; por eso se habla, en la termodin\'amica,
 de una {\em ley} seg\'un la cual la entrop\'ia siempre crece.)

\begin{center}
 * * *
\end{center}

Ahora que vemos que las expresiones del tipo $r \log r$ aparecen en un modelo
 de un fen\'omeno natural, resulta sensato
esperar que el hasta ahora curioso exponente de $a \log a + 1 - a$
en (\ref{eq:zarzu}) y (\ref{eq:vclib}) no sea simplemente una consecuencia
de nuestra incapacidad. Nuestra tarea consiste ahora en dar cotas 
inferiores cercanas a las cotas superiores en (\ref{eq:zarzu}) y
(\ref{eq:vclib}), y de esta manera mostrar que el exponente $a \log a + 1 - a$
verdaderamente describe la probabilidad de las grandes desviaciones.

Comenzemos, como de costumbre, examinando variables mutuamente
independientes $X_2', X_3',\dotsc$ de distribuci\'on
\[X_p' = \begin{cases}1 &\text{con probabilidad $\frac{1}{p}$}\\
0 &\text{con probabilidad $1 - \frac{1}{p}$.} \end{cases}\]
Sea $X' = \sum_{p\leq N} X_p'$. Queremos dar cotas inferiores
para $\Prob(X > a \log \log N)$, $a>1$, y $\Prob(X < a \log \log N)$,
$a<1$.

Utilizaremos el m\'etodo del {\em ladeo exponencial}.\index{ladeo exponencial}
 Definimos
 nuevas variables $Y_2, Y_3, Y_5, \dotsc$ mutuamente independientes y de
 distribuci\'on
\begin{equation}\label{eq:garan}
Y_p = \begin{cases} 1 &\text{con probabilidad $\frac{a}{p} \cdot
\left(1 + \frac{a-1}{p}\right)^{-1}$}\\
0 &\text{con probabilidad $\left(1 - \frac{1}{p}\right) \cdot
\left(1 + \frac{a-1}{p}\right)^{-1}$.} \end{cases}\end{equation}
(Est\'a claro que estamos ``ladeando'' las variables hacia $1$, pero
donde est\'a la ``exponencial''? En $Y_p$, la probabilidad $1/p$ se ha vuelto
$a/p = a^1 \cdot (1/p) \cdot \left(1 + \frac{a-1}{p}\right)^{-1}$; 
la probabilidad $1 - 1/p$ se ha vuelto
$a^0 \cdot (1/p) \cdot \left(1 + \frac{a-1}{p}\right)^{-1}$.
Si $X_p$ tomara los valores $2, 3, \cdots$ con probabilidades no nulas,
multiplicar\'iamos dichas probabilidades por $a^2, a^3,\cdots$, 
respectivamente. El factor de $\left(1 + \frac{a-1}{p}\right)^{-1}$
est\'a all\'i simplemente para asegurar que la suma de las probabilidades
sea $1$.)

Sea $Y = \sum_{p\leq N} Y_p$. El evento $\Prob(Y> a \log \log N)$ no
es una gran desviaci\'on, sino un evento probable. Por el teorema
del l\'imite central,
\begin{equation}\label{eq:deseo}\begin{aligned}
\Prob(a \log \log N < Y \leq a \log \log N + 
\sqrt{a \log \log N}) &= \int_0^1 \frac{1}{\sqrt{2 \pi}} e^{-x^2/2} dx +
o(1)\\ &\geq \frac{1}{5} + o(1),\end{aligned}\end{equation}
(digamos). Sea $I$ el intervalo $(a \log \log N, a \log \log N +
\sqrt{a \log \log N}\rbrack$. Comparemos 
\begin{equation}\label{eq:huh}\begin{aligned} \Prob(X'> a \log \log N) 
&= \mathop{\sum_{\{x_p\}_{p\leq N}:\; x_p \in \{0,1\}}}_{
\sum_{p\leq N} x_p > a \log \log N}
\Prob(X_p' = x_p\;\; \forall p\leq N)\\
&= \mathop{\sum_{\{x_p\}_{p\leq N}:\; x_p \in \{0,1\}}}_{
\sum_{p\leq N} x_p > a \log \log N}
 \prod_{p\leq N} \begin{cases} \frac{1}{p} &\text{si $x_p = 1$}\\
1 - \frac{1}{p} &\text{si $x_p = 0$} \end{cases} 
 \end{aligned}\end{equation}
y \begin{equation}\label{eq:iskay}\begin{aligned}\Prob&(Y\in I)
= \mathop{\sum_{\{x_p\}_{p\leq N}:\; x_p \in \{0,1\}}}_{
\sum_{p\leq N} x_p\; \in I}
\Prob(Y_p = x_p\;\; \forall p\leq N)\\
&= \mathop{\sum_{\{x_p\}_{p\leq N}:\; x_p \in \{0,1\}}}_{
\sum_{p\leq N} x_p\; \in I}
\prod_{p<N} \begin{cases} \frac{a}{p} &\text{si $x_p = 1$}\\
1 - \frac{1}{p} &\text{si $x_p = 0$} \end{cases} \cdot 
\prod_{p\leq N} \left(1 + \frac{a-1}{p}\right)^{-1}\\
&= \mathop{\sum_{\{x_p\}_{p\leq N}:\; x_p \in \{0,1\}}}_{
\sum_{p\leq N} x_p\; \in I}
\Prob(X_p' = x_p\;\; \forall p\leq N) \cdot 
\prod_{p\leq N} \left(1 + \frac{a-1}{p}\right)^{-1} \cdot
\prod_{p\leq N: x_p = 1} a .\end{aligned}\end{equation}

Nuestro deseo es ir de la cota inferior (\ref{eq:deseo}) 
para (\ref{eq:iskay}) a una cota inferior para (\ref{eq:huh}).
Todo t\'ermino $\Prob(X_p' = x_p\;\;\forall p\leq N)$ 
de la suma en (\ref{eq:iskay}) aparece en (\ref{eq:huh}), puesto
que $\sum_{p\leq N} x_p \in I = (a \log \log N, a \log \log N
 + \sqrt{a \log \log N}\rbrack$
implica $\sum_{p\leq N} x_p > a \log \log N$. Claro est\'a, en
(\ref{eq:iskay}),
el t\'ermino $\Prob(X_p' = x_p\;\;\forall p\leq N)$ aparece con dos
factores; uno de ellos es
\[\prod_{p\leq N} \left(1 + \frac{a-1}{p}\right)^{-1} \ll_a
\prod_{p\leq N} \left(1 + \frac{1}{p}\right)^{-(a-1)} \ll_a
(\log N)^{-(a-1)},\] y el otro es
\[\prod_{p\leq N :\; x_p = 1} a = a^{\sum_{p\leq N} x_p}.\]
Como $\sum_{p\leq N} x_p \in I$, sabemos que
$\sum_{p\leq N} x_p \leq a \log \log N + \sqrt{a \log \log N}$. 
Por ende,
\[\prod_{p\leq N :\; x_p = 1} a \leq a^{a \log \log N + \sqrt{a \log \log N}} .\]

Concluimos que
\[\Prob(Y \in I) \ll_a \Prob(X > a \log \log N) \cdot
a^{a \log \log N + \sqrt{a \log \log N}} \cdot (\log N)^{-(a-1)} .\]
Por lo tanto
\[\Prob(X'> a \log \log N) \gg_a \Prob(Y\in I) \cdot
a^{- (a \log \log N + \sqrt{a \log \log N})} \cdot (\log N)^{a-1}.\]
Como sabemos que $\Prob(Y\in I)> \frac{1}{5} + o(1)$ (por 
(\ref{eq:deseo})), obtenemos 
\begin{equation}\label{eq:viol1}\begin{aligned}
\Prob(X'> a \log \log N) &\gg_a \frac{1}{5} a^{- (a \log \log N + \sqrt{a \log
    \log N})} \cdot (\log N)^{a-1}\\ &\gg
(\log N)^{- (a \log a - a + 1) - O((a^{1/2}  \log a) \cdot (\log \log N)^{-1/2})} \end{aligned}
\end{equation}
para $a>1$. Mediante exactamente el mismo m\'etodo, podemos obtener
\begin{equation}\label{eq:viol2}\Prob(X' < a \log \log N) \gg_a (\log N)^{
- (a \log a - a + 1) - O((a^{1/2} \log a) \cdot (\log \log N)^{-1/2})} .\end{equation}
para $a<1$. Hemos obtenido, entonces,
cotas inferiores para complementar a (\ref{eq:zarzu}).
La constante impl\'icita en $\gg$ es continua en $a$.

\begin{center}
* * *
\end{center}
A\'un nos falta derivar cotas inferiores similares para las variables
\[X_p = \begin{cases}1 &\text{si $p|n$}\\ 0 &\text{si $p\nmid n$.}\end{cases}\]
Veremos que hay un m\'etodo muy general para ir de resultados sobre $X_p'$
a resultados sobre $X_p$. Ya pudimos haberlo utilizado para traducir
las cotas superiores $\Prob(X' > a \log \log N) \ll \dots$ en cotas
superiores $\Prob(X > a \log \log N) \ll \dots$. Hemos 
esperado hasta ahora en parte
porque el m\'etodo depende de un resultado t\'ecnico com\'unmente
considerado dif\'icil: el {\em lema fundamental de la teor\'ia de cribas}.
\index{criba!lema fundamental}
(Ver el problema \ref{it:gantost} al final de esta secci\'on;
en \'este se desarrolla una prueba de una
 versi\'on d\'ebil del lema fundamental -- la \'unica versi\'on que
necesitaremos.)

Como de costumbre, comenzaremos truncando las sumas. Definimos
$S_m = \sum_{p\leq m} X_p$, y, en vez de trabajar con 
$X = S_N = \sum_{p\leq N} X_p$,
trabajaremos con $S_{g(N)}$, donde $g(N)$ es un poco menor que $N$. Acotaremos la diferencia $X - S_{g(N)}$ al final,
ya que a estas alturas esa parte es rutinaria. Escogeremos $g(N)$ en el
\'ultimo paso. 

Sea $S'_{m} = \sum_{p\leq m} X_p'$. Debemos comparar $S_{g(N)}$ con
$S_{g(N)}'$, o, m\'as precisamente, mostrar que $\Prob(S_{g(N)}> a \log \log N)$
(o $\Prob(S_{g(N)} < a \log \log N)$) es aproximadamente igual a
 $\Prob(S'_{g(N)} > a \log \log N)$
(o con $\Prob(S_{g(N)}' < a \log \log N)$). Ya tenemos cotas para
$S'_{g(N)}$, gracias a (\ref{eq:viol1}) y (\ref{eq:viol2}) (con $X=S'_{g(N)}$).

Tenemos, por una parte,
\begin{equation}\label{eq:nanette}\Prob(S_{g(N)} > a \log \log N) =
\mathop{\sum_{\{x_p\}_{p\leq g(N)}:\; x_p \in \{0,1\}}}_{\sum_{p\leq g(N)} x_p
> a \log \log N} \Prob(X_p = x_p\;\; \forall p\leq g(N)),\end{equation}
y, por otra,
\begin{equation}\label{eq:edi}\begin{aligned}
\Prob(S'_{g(N)} > a \log \log N) &=
\mathop{\sum_{\{x_p\}_{p\leq g(N)}:\; x_p \in \{0,1\}}}_{\sum_{p\leq g(N)} x_p
> a \log \log N} \Prob(X_p' = x_p\;\; \forall p\leq g(N))\\ &=
\mathop{\sum_{\{x_p\}_{p\leq g(N)}:\; x_p \in \{0,1\}}}_{\sum_{p\leq g(N)} x_p
> a \log \log N} \prod_{p\leq g(N)} 
 \begin{cases} \frac{1}{p} &\text{si $x_p = 1$}\\
1 - \frac{1}{p} &\text{si $x_p = 0$.}\end{cases}\end{aligned}\end{equation}
Consideraremos primero los t\'erminos de (\ref{eq:nanette})
con $\prod_{p: x_p = 1} p \leq N^{1-\epsilon}$.
(Aqu\'i $\epsilon>0$ es un n\'umero fijo cualquiera entre $0$ y $1$.) 

Sea $m = \prod_{p: x_p = 1} p$. Supongamos que $m\leq N^{1 - \epsilon}$.
Sea $S = \{p\leq g(N) : x_p = 0\}$. Entonces
\[\Prob(X_p = x_p\;\;\forall p \leq g(N))
= \frac{1}{m} \Prob(X_p = 0\;\;\text{$\forall p\in S$}),
\]
donde en el lado derecho de la ecuaci\'on las variables $X_p$ dependen
de un n\'umero $n$ tomado al azar entre $1$ y $N/m$, no entre $1$ y $N$.
Por el lema fundamental\index{criba!lema fundamental} de las cribas (problema \ref{it:gantost},
teorema \ref{thm:fundy}),
\[\begin{aligned}
\Prob(X_p = 0\;\;\text{$\forall p \in S$}) &=
 \prod_{p\in S} \left(1 - \frac{1}{p}\right) \cdot (1 + O_A((\log N/m)^{-A}))\\
&= \prod_{p\in S} \left(1 - \frac{1}{p}\right) \cdot
(1 + O_A((\log N)^{-A}))
\end{aligned}\]
para cualquier $A$, siempre y cuando
$\log g(N) \ll \frac{\log N/m}{\log \log (N/m)}$.
Esto \'ultimo ciertamente tiene lugar si 
 $g(N) = o\left( \frac{\log N}{\log \log N}\right)$ (donde utilizamos el hecho
que $m\leq N^{1-\epsilon}$). 

Asumamos, entonces, $g(N)\ll \frac{\log N}{\log \log N}$. Podemos
entonces concluir que
\[
\mathop{\mathop{\sum_{\{x_p\}_{p\leq g(N)}:\; x_p \in \{0,1\}}}_{\sum_{p\leq g(N)} x_p > a \log \log N}}
_{\prod_{p : x_p = 1} p \;\leq N^{1 - \epsilon}}
 \Prob(X_p = x_p\;\; \forall p\leq g(N))\]
es igual a
\[
\mathop{\mathop{\sum_{\{x_p\}_{p\leq g(N)}:\; x_p \in \{0,1\}}}_{\sum_{p\leq g(N)} x_p > a \log \log N}}_{\prod_{p : x_p = 1} p \;\leq N^{1 - \epsilon}}
\frac{1}{
\prod_{p: x_p = 1} p}
\prod_{p: x_p = 0} \left(1 - \frac{1}{p}\right)
\cdot (1 + O_A((\log N)^{-A})) ,
\]
lo cual no es sino 
\[(1 + O_A((\log N)^{-A}) \cdot
\mathop{\mathop{\sum_{\{x_p\}_{p\leq g(N)}:\; x_p \in \{0,1\}}}_{\sum_{p\leq g(N)} x_p > a \log \log N}}_{\prod_{p : x_p = 1} p \;\leq N^{1 - \epsilon}}
\Prob(X_p' = x_p \;\; \forall p\leq g(N)),\]
es decir, 
la suma de los t\'erminos de (\ref{eq:edi}) con
$\prod_{p: x_p = 1} p \leq N^{1 - \epsilon}$, multiplicada por
$(1 + O_A((\log N)^{-A}))$ . (Este es el paso crucial del m\'etodo:
hemos logrado ir de una suma que involucra a
$\Prob(X_p = x_p\;\; \forall p\leq g(N))$ a una suma que involucra a
$\Prob(X_p' = x_p\;\; \forall p\leq g(N))$, donde las variables $X_p'$
son las variables mutuamente independientes que estudiamos al principio de la
secci\'on.)

Nos queda acotar los t\'erminos de (\ref{eq:nanette}) y (\ref{eq:edi})
con $\prod_{p : x_p = 1} p > N^{1 - \epsilon}$. En un caso como el otro,
el total de tales t\'erminos es a lo m\'as
\[
\mathop{\sum_{Q \subset P}}_{\prod_{p\in Q} p > N^{1 - \epsilon}}
 \frac{1}{\prod_{p\in Q} p},
\]
donde $P$ es el conjunto de los primos $p\leq g(N)$.
Por un resultado intermedio (\ref{eq:menuh}) 
en la prueba del lema fundamental de las cribas,\index{criba!lema fundamental}
\[\mathop{\sum_{Q \subset P}}_{\prod_{p\in Q} p > N^{1 - \epsilon}}
 \frac{1}{\prod_{p\in Q} p} \ll (\log g(N)) \cdot \left(\log \left(
\frac{N^{1 - \epsilon}}{g(N)}\right)\right)
^{- \left(\log \left(
\frac{N^{1 - \epsilon}}{g(N)}\right)\right)
\cdot (1 + o(1))} \ll_A (\log N)^{-A}\]
para cualquier $A>0$, donde utilizamos el hecho que $\log(N) = o\left(
\frac{\log N}{\log \log N}\right)$. 
% Como ya sabemos
%((\ref{eq:viol1}) y (\ref{eq:viol2})) que 
%\[\begin{aligned}
%\Prob(S_{g(N)}' > a \log \log N) \gg (\log N)^{- (a \log a - a + 1) -
%O_a((\log \log N)^{-1/2})
%\\
%\end{aligned}\]

Concluimos que
\begin{equation}\label{eq:hijamia}\Prob(S_{g(N)} > a \log \log N) = \Prob(S_{g(N)}' > a \log \log N) +
O_A\left((\log N)^{-A}\right) \end{equation}
para cualquier $A>0$,
bajo la condici\'on que $g(N) = o\left(\frac{\log N}{\log \log N}\right)$.

La desigualdad (\ref{eq:viol1}) puede aplicarse directamente
a $S_{g(N)}$: como la definici\'on de las variables $X_p'$ no involucra
a $N$, podemos simplemente utilizar $g(N)$ en vez de $N$. Obtenemos
\[\Prob(S_{g(N)}' > a \log \log N) \gg (\log g(N))^{-(a \log a - a + 1) + o(1)}
.\]
Si $\log \log N - \log \log g(N) = o(\log \log N)$, entonces
$\log g(N) > (\log N)^{1 - o(1)}$. Por ende,
\[
\Prob(S_{g(N)}' > a \log \log N) \gg (\log N)^{-(a \log a - a + 1) + o(1)}\]
bajo la condici\'on que $\log \log N - \log \log g(N) = o(\log \log N)$.
La ecuaci\'on (\ref{eq:hijamia}) nos permite deducir de esto que
\begin{equation}\label{eq:inti}
\Prob(S_{g(N)} > a \log \log N) \gg (\log N)^{-(a \log a - a + 1) + o(1)}
\end{equation}

Empero,
procederemos de manera distinta. Mostraremos que la probabilidad que
$X - S_{g(N)} > \epsilon' \log \log N$ se satisfaga es muy peque\~na
-- si es que $g(N)$ satisface una cierta condici\'on f\'acil de cumplir.

Por Chebyshev-Mertens, $\sum_{g(N) <p \leq N} \frac{1}{p} = \log \log N
- \log \log g(N) + O(1)$. Procediendo exactamente como en \S \ref{sec:bambi} 
(problemas \ref{it:ino}, \ref{it:cencia} y comienzo de \ref{it:irmin}), podemos mostrar que
\[\mathbb{E}(e^{\alpha (X - S_{g(N)})}) \leq 
\left(e^{\log \log N - \log \log g(N) + O(1)}\right)^{e^{\alpha}}\]
para todo $\alpha$. (Podr\'iamos continuar como en el resto de
\S \ref{sec:bambi}, problema \ref{not:antartor}, y reemplazar
$e^{\alpha}$ por $e^{\alpha} - 1$, pero esto no ser\'a necesario.)
Por Markov,
\begin{equation}\label{eq:irlan}\Prob(X - S_{g(N)}>x) \leq \frac{
\left(e^{\log \log N - \log \log g(N) + O(1)}\right)^{e^{\alpha}}}{
e^{\alpha x}}\end{equation}
para cualquier $\alpha$. Si $\log \log N - \log \log g(N) + O(1)
= o(x)$, podemos escoger $\alpha = A+1$ para $A$ arbitariamente grande,
y entonces (\ref{eq:irlan}) nos da que
\begin{equation}\label{eq:antoni}
\Prob(X - S_{g(N)}>x) \ll_A e^{-A x}\end{equation}
para $N$ suficientemente grande. Escogemos $x = \epsilon' \log \log N$
(con $\epsilon'>0$ arbitrariamente peque\~no) y $A = A'/\epsilon'$
(para $A'$ arbitrariamente grande). Concluimos que
\begin{equation}\label{eq:killa}
\Prob(X - S_{g(N)} > \epsilon' \log \log N) \ll_{A'} (\log N)^{-A'}\end{equation}
para $N$ suficientemente grande, si se cumple la condici\'on que
$\log \log N - \log \log g(N) = o(\log \log N)$.

\begin{figure}
\centering \includegraphics[height=4in]{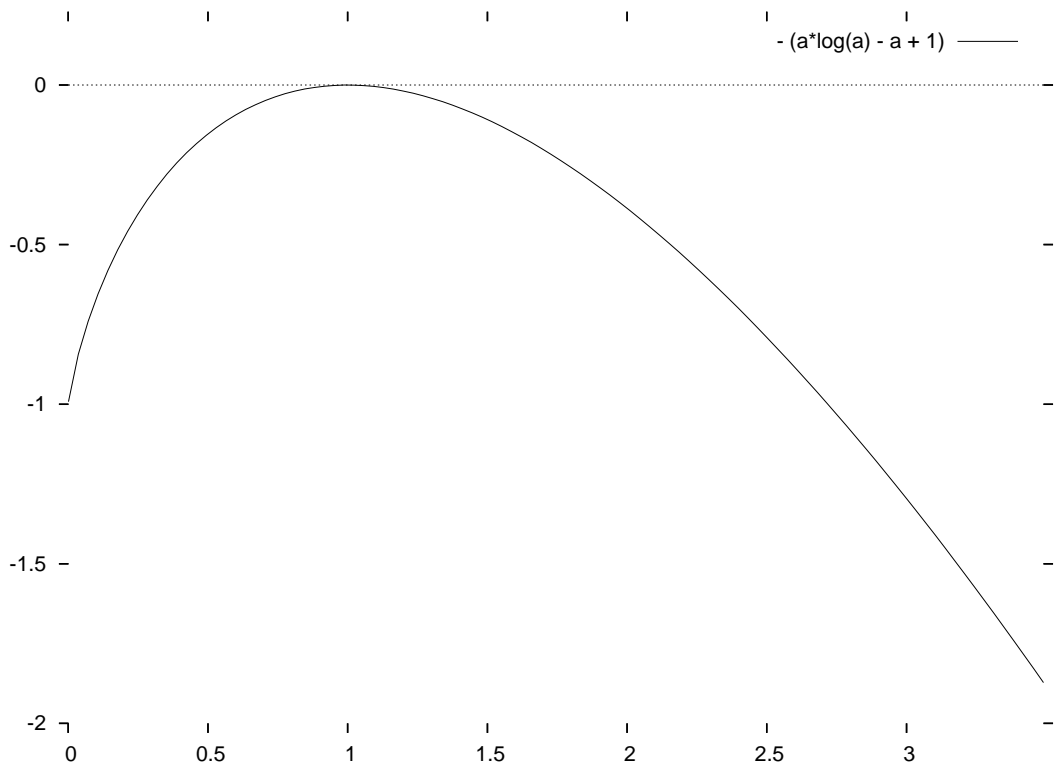}
\caption{La distribuci\'on de $\omega(n)$, vista desde lejos, en escala
logar\'itmica.
\'Este es el gr\'afico de $y = -I(x) = - (x \log x - x + 1)$. 
La probabilidad de $\omega(n) > t \log \log N$ (si $t>1$) o
$\omega(n) < t \log \log N$ (si $t<1$) es igual a
 $(\log N)^{-I(t) + o(1)}$, lo cual es lo mismo que
$(\log N)^{o(1)} \cdot \int_{t}^\infty (\log N)^{-I(x)} dx$ (si $t>1$)
o $(\log N)^{o(1)} \cdot \int_{0}^t (\log N)^{-I(x)} dx$ (si $t<1$).
}\label{fig:grande}
\end{figure}

Atemos los cabos. Por (\ref{eq:inti}),
\begin{equation}\label{eq:jlobio}\Prob(X>a \log \log N) \geq \Prob(S'_{g(N)} > a \log \log N)
\gg (\log N)^{- (a \log a - a + 1) + o(1)}\end{equation} para $a>1$. Pasemos al caso
$a<1$. Por (\ref{eq:inti}) y (\ref{eq:killa}),
\[\begin{aligned}
\Prob(X<a \log \log N) &\geq \Prob(S'_{g(N)} < (a - \epsilon') \log \log N)\\
&- \Prob(X - S_{g(N)} > \epsilon' \log \log N)\\ &\gg
(\log N)^{- ((a - \epsilon') \log (a - \epsilon') - (a - \epsilon') + 1) +
  o(1)} - 
O_{A',\epsilon'}((\log N)^{-A'})\end{aligned}\]
para $A'$ arbitrariamente grande y $\epsilon'>0$ arbitrariamente peque\~no. Como $t \to t \log t - t + 1$ es
continua, esto implica que
\begin{equation}\label{eq:jodo}
\Prob(X > a \log \log N) \gg (\log N)^{- (a \log a - a + 1) + o(1)}.\end{equation}

 Queda solamente verificar que existe un $g(N)$ que
satisface las condiciones impuestas: $\log g(N) = o\left(
 \frac{\log N}{\log \log N}\right)$
y $\log \log N - \log \log g(N) = o(\log \log N)$. La funci\'on
$g(N) = \frac{\log N}{(\log \log N)^2}$ (por ejemplo) satisface
ambas condiciones.

Hemos obtenido nuestro objetivo. Recordando
(\ref{eq:vclib}) y $X = \omega(n)$, concluimos que
\[\begin{aligned}
(\log N)^{ - (a \log a + 1 - a) + o(1)} &\ll_a 
\Prob(\omega(n) > a \log \log N) &\ll_a (\log N)^{- (a \log a + 1 - a)}
&\text{\; si $a>1$,}\\
(\log N)^{ - (a \log a + 1 - a) + o(1)} &\ll_a 
\Prob(\omega(n) <a \log \log N) &\ll_a (\log N)^{- (a \log a + 1 - a)}
&\text{\; si $a<1$,} 
\end{aligned}
\]
donde $n$ es tomado al azar entre $1$ y $N$. Las constantes dependen de
$a$ de manera continua.

\begin{center}
{\bf Notas y problemas}
\end{center}
\begin{enumerate}
\item\label{it:gantost} 
{\em Lema fundamental de las cribas (versi\'on d\'ebil).}\index{criba!lema fundamental}
\begin{enumerate}
Sea $z = N^{1/s}$, donde $s\to \infty$ cuando $N\to \infty$.
Sea $P$ un conjunto de primos $\leq z$. Queremos determinar
cu\'antos enteros $n\leq N$ son coprimos con todo $p\in P$. 

\item 

Cu\'antos enteros $n\leq N$ son impares? La respuesta es el n\'umero
de enteros $n\leq N$ (es decir, $N$) menos el n\'umero de enteros
$n\leq N$ que son pares (es decir, $\lfloor N/2\rfloor$).
Muy bien -- la respuesta es $N - \lfloor N/2\rfloor = N/2 + O(1)$. 
Cu\'antos enteros $n\leq N$ son coprimos con $2$ y $3$? Tomamos
los $N$ enteros $n\leq N$, restamos los $\lfloor N/2\rfloor$ enteros
divisibles por $2$ y los $\lfloor N/3\rfloor$ enteros divisibles por $3$,
y nos damos cuenta que hemos sustra\'ido los enteros divisibles tanto
por $2$ como por $3$ (es decir, divisibles por $6$)
por partida doble; tendremos que contarlos una vez de vuelta. Vemos, entonces,
que el n\'umero de enteros $n\leq N$ coprimos con $2$ y $3$ es
\[\begin{aligned}
N - \left\lfloor \frac{N}{2} \right\rfloor - 
\left\lfloor \frac{N}{3} \right\rfloor +
\left\lfloor \frac{N}{6}\right\rfloor &= 
N - \frac{N}{2} - \frac{N}{3} + \frac{N}{6} + 3\cdot O(1)\\
&= \left(1 - \frac{1}{2}\right) \left(1 - \frac{1}{3}\right) N + O(1)
.\end{aligned}\]

Seguimos razonando de la misma manera (enfoque que tiene el nombre de
\index{principio de inclusi\'on-exclusi\'on}
{\em principio de
  inclusi\'on-exclusi\'on}) y concluimos que la probabilidad
que un $n\leq N$ tomado al azar sea coprimo con todo $p\in P$ es
igual a 
\begin{equation}\label{eq:line}\begin{aligned}
\sum_{Q \subset P} (-1)^{|Q|} &\Prob(\text{$n$ es divisible por todo
$p\in Q$}) =
\sum_{Q \subset P} (-1)^{|Q|} \frac{1}{N} \left\lfloor \frac{N}{\prod_{p\in Q} p}
\right\rfloor\\
&= \mathop{\sum_{Q \subset P}}_{\prod_{p\in Q} p \leq N^{1 - \epsilon}}
 (-1)^{|Q|} \frac{1}{N} \left\lfloor \frac{N}{\prod_{p\in Q} p} \right\rfloor
+ \mathop{\sum_{Q \subset P}}_{\prod_{p\in Q} p > N^{1 - \epsilon}}
 (-1)^{|Q|} \frac{1}{N} \left\lfloor \frac{N}{\prod_{p\in Q} p} \right\rfloor
\end{aligned}\end{equation}
para cualquier $\epsilon>0$. Muestre que esto es
\[\prod_{p\in P} \left(1 - \frac{1}{p}\right) + O(\text{error}) +
O(N^{-\epsilon}) +
 O(\text{error}),\]
donde
\begin{equation}\label{eq:anace}
\text{error} = 
\mathop{\sum_{Q \subset P}}_{\prod_{p\in Q} p > N^{1 - \epsilon}}
 \frac{1}{\prod_{p\in Q} p} .\end{equation}
Habremos terminado una vez que acotemos este t\'ermino de error, 

\item Podemos escribir (\ref{eq:anace}) de la siguiente manera:
\begin{equation}\label{eq:vieja}
\text{error} = \sum_{k=0}^\infty\;\;
\mathop{\sideset{}{^*}\sum_{N^{1 - \epsilon} 2^k < n \leq N^{1 - \epsilon} 2^{k+1}}}_{p|n \Rightarrow
p\leq z} \frac{1}{n}
\leq
\sum_{k=0}^\infty \frac{1}{N^{1 - \epsilon} 2^k}
\mathop{\sideset{}{^*}\sum_{N^{1 - \epsilon} 2^k < n \leq N^{1 - \epsilon} 2^{k+1}}}_{p|n \Rightarrow
p\leq z} 1
,\end{equation}
donde $\sum^*$ denota que la suma recorre
s\'olo aquellos enteros $n$ que no tienen divisores cuadrados.
Lo que es verdaderamente crucial es que
todo n\'umero $n$ en la sumas $\sum_n$ en (\ref{eq:vieja}) es {\em desmenuzable},\index{n\'umeros desmenuzables} i.e., no tiene
factores primos grandes. 

Podemos, entonces, utilizar (\ref{eq:astrid}), que nos
da la (m\'as bien peque\~na) probabilidad que un n\'umero sea desmenuzable (y libre
de divisores cuadrados). Obtenemos 
\begin{equation}\label{eq:menuh}
\text{error} \leq \sum_{k=0}^{\infty} u_k^{-u_k (1 + o(1))},\end{equation}
donde $u_k = \frac{\log N^{1-\epsilon} 2^k}{\log z}$. 

\item Como
$u_k = a + k b$, donde $a = \frac{\log N^{1 - \epsilon}}{\log z} = (1 -
\epsilon) s$ y $b = \frac{\log 2}{\log z}$, tenemos
\[\begin{aligned}u_k^{-u_k (1 + o(1))} &= 
(a + k b)^{- (a + k b) (1 + o(1))}\\ &\leq 
 a^{-a (1 + o(1))}  e^{- k b (1+ o(1)) \log a} .\end{aligned}\]
Por lo tanto
\[\text{error} \leq a^{-a (1 + o(1))} \sum_{k=0}^{\infty} e^{- k \delta},\]
donde $\delta = b (1 + o(1)) \log a$. 
Muestre que 
$\sum_{k=0}^{\infty} e^{-k \delta} \ll \frac{1}{\delta}$ para $\delta>0$.
As\'i,
\[\begin{aligned}
\text{error}&\ll \frac{a^{-a (1 + o(1))}}{\delta} \ll \frac{\log z}{\log s}
\cdot ((1-\epsilon) s)^{- (1 - \epsilon) s \cdot (1 + o(1))}\\ &\ll
(\log z) s^{- (1-2\epsilon) s (1+o(1))}.\end{aligned}\]
El error $(\log z) \cdot s^{- (1 - 2 \epsilon) s \cdot (1 + o(1))}$
 ser\'a $O((\log N)^{-A})$ cuando
$s> 3 A \frac{\log \log N}{\log \log \log N}$.
Si $s\gg \log \log N$, el error
ser\'a $\ll_A (\log N)^{-A}$ para todo $A>0$.

%, y, en particular,
%ser\'a mucho m\'as peque\~no que el t\'ermino principal.
\item Hemos probado
\begin{thm}[Lema fundamental de las cribas, version d\'ebil]\label{thm:fundy}
\index{criba!lema fundamental}
Sea $z = N^{1/s}$, donde $s\gg \log \log N$.
Sea $P$ un subconjunto de 
$\{p\leq z: \text{$p$ primo}\}$. Entonces, para $n$ tomado
al azar entre $1$ y $N$,
\begin{equation}\label{eq:lemun}
\Prob(\text{$n$ es coprimo con todo $p\in P$}) =
(1 + o(1))  \cdot \prod_{p\in P} \left(1 - \frac{1}{p}\right) 
.\end{equation}
Expl\'icitamente,
\begin{equation}\label{eq:hawaii}\begin{aligned}
\Prob&(\text{$n$ es coprimo con todo $p\in P$})\\
 &=
\prod_{p\in P} \left(1 - \frac{1}{p}\right) 
+ 
O\left((\log z) \cdot s^{- (1 - 2 \epsilon) s \cdot (1 + o(1))}\right) + 
O(N^{1-\epsilon})\\
&= (1 + O_A((\log N)^{-A})) \cdot 
\prod_{p\in P} \left(1 - \frac{1}{p}\right)
\end{aligned}\end{equation}
para cualquier $\epsilon>0$ y cualquier $A>0$.
\end{thm}

{\em Nota.} 
La versi\'on fuerte del lema fundamental de las cribas 
consiste en la aseveraci\'on que (\ref{eq:lemun})
rige no s\'olo para $s\gg \log \log N$, sino 
para todo $s \to \infty$. No probaremos la versi\'on fuerte aqu\'i.\\

\end{enumerate}

\item El lema fundamental es 
muy razonable: como la probabilidad que $n$ sea coprimo con un $p\in P$
dado es $1 - 1/p$, la ecuaci\'on (\ref{eq:lemun}) nos dice (como muchas
otras cosas que hemos probado) que
los eventos $p|n$ se comportan hasta cierto punto como variables mutuamente independientes.
Empero, 
si $z = N^{1/s}$ y $s$ no va a $\infty$ cuando $N\to \infty$,  
la ley (\ref{eq:lemun}) no rige.
%(Cuando $s\to \infty$ se cumple pero $s = o(\log \log N)$, entonces,
%como dijimos, (\ref{eq:lemun}) es a\'un cierta; simplemente no la hemos
%llegado a probar.)

Mostremos esto en el caso m\'as simple.
Si $z = N^{1/2}$ y $P$ es el conjunto de todos los primos
$p\leq P$, entonces $n$ es coprimo con todo $p\in P$ s\'i y s\'olo s\'i
$n$ es primo. La probabilidad que $n$ sea primo es
\begin{equation}\label{eq:realidad}(1 + o(1)) \cdot \frac{1}{\log N}\end{equation}
(este es el {\em teorema de los n\'umeros primos},
\index{teorema!de los n\'umeros primos}
 el cual no probaremos).
La aseveraci\'on (\ref{eq:lemun}) nos dar\'ia, de otro lado,
$(1 + o(1)) \cdot \prod_{p\leq N^{1/2}} \left(1 - \frac{1}{p}\right)$. Se puede mostrar
que
\begin{equation}\label{eq:prediccion}
\prod_{p\leq N^{1/2}} \left(1 - \frac{1}{p}\right) = 
\frac{2 e^{-\gamma}}{\log N} (1 + o(1)),\end{equation}
donde $\gamma$ es la constante de Euler $\gamma = 0.577\dotsc$. Lo
importante aqu\'i es que $2 e^{-\gamma} \ne 1$, por lo cual la
predicci\'on natural (\ref{eq:prediccion}) no es compatible con la
realidad (\ref{eq:realidad}).\\

\item {\em Desviaciones moderadas.}\index{desviaciones!moderadas}

\begin{enumerate} \item
Si $X - \mathbb{E}(X)$ est\'a
en la escala de $\sqrt{\Var(X)}$ (es decir, entre
$\frac{1}{1000} \sqrt{\Var(X)}$ y $1000 \sqrt{\Var(X)}$, por ejemplo),
hablamos del {\em l\'imite central},\index{l\'imite central}
 o de {\em peque\~nas desviaciones}.
Si $X - \mathbb{E}(X)$ est\'a en la escala de $\Var(X)$ o $\mathbb{E}(X)$,
hablamos de {\em grandes desviaciones}.\index{desviaciones!grandes} A\'un no hemos examinado el caso
en el cual $X - \mathbb{E}(X)$ es bastante m\'as grande que $\sqrt{\Var(X)}$
y bastante m\'as peque\~no que $\Var(X)$; como cabr\'ia esperar, esto se
llama una {\em desviaci\'on moderada}.

Lo que queremos examinar es
\[\Prob(X > (1 + \Delta(N)) \log \log N) \text{\; y }
\Prob(X < (1 - \Delta(N)) \log \log N),\]
donde $\Delta(N)>0$, $\Delta(N) = o(1)$ y $(\log \log N)^{-1/2} = o(\Delta(N))$.
(Si $\Delta(N) = o(1)$ no se cumpliera, tendr\'iamos una gran desviaci\'on;
si $(\log \log N)^{-1/2} = o(\Delta(N))$ no se cumpliera, tendr\'iamos una
peque\~na desviaci\'on.)

\item Las cotas superiores de grandes desviaciones son a\'un v\'alidas, como
podemos ver repasando sus pruebas.
Por lo tanto:
\[\Prob(X > (1 + \Delta(N)) \log \log N) \leq (\log N)^{-I(1 + \Delta(N))},\]
donde $I(a) = a \log a - a + 1$. Mediante una serie de Taylor, verifique que
\[I(1 + \Delta(N)) = \frac{1}{2} \Delta^2(N) + O(\Delta^3(N)) =
\frac{1}{2} \Delta^2(N) \cdot (1 + o(1)).
\]
Concluimos que
\[\Prob(X > (1 + \Delta(N)) \log \log N) \leq (\log N)^{-\frac{1}{2} 
\Delta^2(N) \cdot (1 + o(1))},\]
y, similarmente,
\[\Prob(X < (1 - \Delta(N)) \log \log N) \leq (\log N)^{-\frac{1}{2} 
\Delta^2(N) \cdot (1 + o(1))}.\]
\item Las cotas inferiores tendr\'an que ser rehechas con m\'as cuidado:
el sumando de $o(1)$ en el exponente de (\ref{eq:jlobio}) y (\ref{eq:jodo})
es ahora demasiado burdo, ya que $- (a \log a + 1 - a)$ ser\'a
bastante m\'as peque\~no que $o(1)$.

Las cotas (\ref{eq:viol1}) y (\ref{eq:viol2}) son a\'un v\'alidas. 
Muestre que, como estamos asumiendo $(\log \log N)^{-1/2} = o(\Delta(N))$,
las cotas (\ref{eq:viol1}) y (\ref{eq:viol2}) toman la forma
\begin{equation}\label{eq:rino}\begin{aligned}
\Prob(X' > (1 + \Delta(N)) \log \log N) &\gg (\log N)^{-\frac{1}{2} 
\Delta^2(N) \cdot (1 + o(1))}\\
\Prob(X' < (1 - \Delta(N)) \log \log N) &\gg (\log N)^{-\frac{1}{2} 
\Delta^2(N) \cdot (1 + o(1))}.\end{aligned}\end{equation}
\item La transici\'on de $X'$ a $X$ se hace como antes. La \'unica
dificultad reside en el hecho que ya no basta probar que
$\Prob(X - S_{g(N)} > \epsilon' \log \log N)$ es peque\~na;
debemos asegurarnos que $\Prob(X - S_{g(N)} > \epsilon' \Delta(N) \log \log
N)$
sea peque\~na. Para esto tendremos que modificar el valor de $g(N)$.

Sustituya $x = \epsilon' \Delta(N) \log \log N$ y $A = 1/\epsilon$
 en (\ref{eq:antoni}) y 
obtenga
\begin{equation}\label{eq:respig}\Prob(X - S_{g(N)} > \epsilon' \Delta(N) \log \log N) \leq (\log N)^{- \Delta(N)},
\end{equation}
bajo la suposici\'on que $\log \log N - \log \log g(N) + O(1) = o(\Delta(N)
\log \log N)$. 
Como $\frac{1}{2} \Delta^2(N) = o(\Delta(N))$, el t\'ermino de error
(\ref{eq:respig}) es mucho m\'as chico que el t\'ermino principal
$(\log N)^{- \frac{1}{2} \Delta^2(N) (1 + o(1))}$.

\item Muestre que se puede definir $g(N)$ de tal manera que 
\[\log g(N) = o\left(\frac{\log N}{\log \log N}\right)\;\;\;\;\;\;\;\text{y}
\;\;\;\;\;\;\;
\log \log N - \log \log g(N) = o(\Delta(N) \log \log N).\] Recuerde que
$\Delta(N)\gg \sqrt{\log \log N}$.) Concluya que
\begin{equation}\label{eq:otto}\begin{aligned}
\Prob(X' > (1 + \Delta(N)) \log \log N) &\ll (\log N)^{-\frac{1}{2} 
\Delta^2(N) \cdot (1 + o(1))}\\
\Prob(X' < (1 - \Delta(N)) \log \log N) &\ll (\log N)^{-\frac{1}{2} 
\Delta^2(N) \cdot (1 + o(1))}.\end{aligned}\end{equation}

\item
Podemos escribir $t = \Delta(N) \cdot (\log \log N)^{1/2}$. (La
desviaci\'on es entonces $\Delta(N) \cdot (\log \log N) = 
t \sqrt{\log \log N}$.) Est\'a claro que
\[(\log N)^{-\frac{1}{2} 
\Delta^2(N) \cdot (1 + o(1))} = e^{-\frac{1}{2} t^2 (1 + o(1))}\]
As\'i, (\ref{eq:rino}) y (\ref{eq:otto}) nos dicen que
la normal nos da, por lo menos, una idea aproximada
de la escala de la probabilidad de las desviaciones moderadas.\\
\end{enumerate}

\item Podemos utilizar el lema fundamental (a\'un en su versi\'on debil)
para dar una prueba alternativa del teorema de Erd\H{o}s-Kac (teorema
\ref{thm:erdkac}). (Lo siguiente est\'a muy cercano del camino
 seguido originalmente por Erd\H{o}s y Kac mismos.) 
\index{teorema!de Erd\H{o}s-Kac}
Sea $g(x)$ una funci\'on tal que
$g(x) = O\left(x^{1/\log \log x}\right)$ y 
$\log \log x - \log \log g(x) = o(\sqrt{\log \log N})$; podemos
tomar $g(x) = x^{1/\log \log x}$, como en la primera prueba que dimos
del teorema. 

Proc\'edase como en esa prueba hasta (\ref{eq:moesg}), i.e.,  muestre
que podemos trabajar con la suma truncada 
\[S_{g(N)} = \frac{1}{\sqrt{\log \log g(N)}} \cdot \sum_{p\leq g(N)}
 (X_p - 1/p).\]
Podemos luego utilizar
el lema fundamental en la misma manera que lo usamos en esta secci\'on,
y as\'i mostrar que
\[\Prob(S_{g(N)}\leq t) = \Prob(S'_{g(N)}\leq t) + o(1)\]
para todo $t$, donde $S_{g(N)}'$ es la suma de variables
mutuamente independientes $X'_p$:
\[S_{g(N)}' = \frac{1}{\sqrt{\log \log g(N)}} \cdot 
 \sum_{p\leq g(N)} (X_p' - 1/p).\] 
En otras palabras, $S_{g(N)}$
tiende a la misma distribuci\'on que $S'_{g(N)}$ -- es decir, tiende a la
normal.

\item {\em Entrop\'ia relativa.}\index{entrop\'ia!relativa}
Consideremos tres compartimientos separados por membranas porosas.
Llen\'emoslos de gas. Sea $p_j$ la probabilidad que una part\'icula
dada (todas son intercambiables) est\'e en el  compartimiento $j$:
\begin{center}
\includegraphics[height=2in]{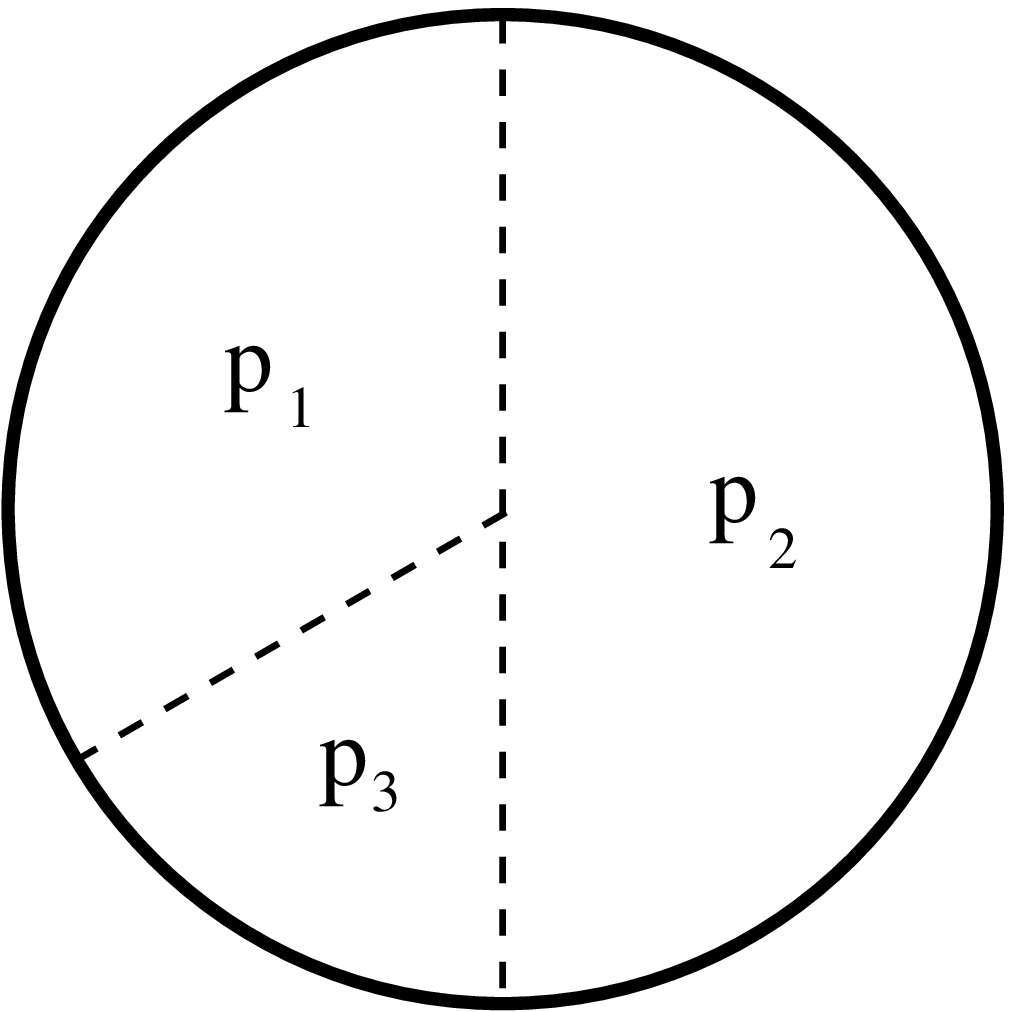}
\end{center}
Cu\'al es la probabilidad que haya $n_1 = r_1\cdot n$ part\'iculas
en el primer
compartimiento, $n_2 = r_2 \cdot n$ en el segundo, y $n_3 = r_3\cdot n$
en el tercero? 

La probabilidad de una configuraci\'on espec\'ifica -- es decir, la
probabilidad que $n_j$ part\'iculas espec\'ificas est\'en en
el compartimiento $j$ -- es
\[p_1^{n_1} \cdot p_2^{n_2} \cdot p_3^{n_3}.\]
El n\'umero de configuraciones con $n_j$ part\'iculas en la c\'amara
$j$ es
\[\frac{n!}{n_1! \cdot n_2! \cdot n_3!} .\]
Luego, la probabilidad que haya $n_j = r_j \cdot n$ part\'iculas en
el compartimiento $j$ es
\[P = p_1^{n_1} p_2^{n_2} p_3^{n_3} \frac{n!}{n_1! \cdot n_2! \cdot n_3!}.\]

Cu\'anto es esto, aproximadamente? Como en el caso que vimos al principio
de la secci\'on, extraemos el logaritmo:
\[\begin{aligned}
\log P &= n_1 \log p_1 + n_2 \log p_2 + n_3 \log p_3
+ n \log n - n + O(\log n)\\
&- ((n_1 \log n_1 - n_1) + (n_2 \log n_2 - n_2) + (n_3 \log n_3 - n_3) +
O(\log n))\\
&= - \left( n_1 \log \frac{n_1/n}{p_1} + n_2 \log \frac{n_2/n}{p_2} +
n_3 \log \frac{n_3/n}{p_3}\right)+ O(\log n)\\
&= - n \left(r_1 \log \frac{r_1}{p_1} + r_2 \log \frac{r_2}{p_2} +
r_3 \log \frac{r_3}{p_3}\right) + O(\log n).
\end{aligned}\]
Por lo tanto,
\[P = e^{- n (H + o(1))}\]
donde $H = r_1 \log \frac{r_1}{p_1} + r_2 \log \frac{r_2}{p_2} +
r_3 \log \frac{r_3}{p_3}$.
En general, definimos
\[H = \sum r_j \log \frac{r_j}{p_j} \text{\;\;\;\;\; (entrop\'ia relativa)}\]
y tenemos que la probabilidad que una proporci\'on $r_j$ de las part\'iculas
est\'en en el comportamiento $j$ es
\[P = e^{- n (H + o(1))} .\]

\item {\em La entrop\'ia relativa y los primos.} Sea una partici\'on
\[\{\text{los primos}\} = P_1 \cup P_2 \cup P_3
\text{\;\;\;\;\; $(P_1 \cap P_2 = P_1 \cap P_3 = P_2 \cap P_3 = \emptyset)$}\]
tal que, para cada $j=1,2,3$,
\[\left(\frac{\log z_1}{\log z_0}\right)^{p_j} \ll
\mathop{\prod_{p\in P_j}}_{z_0 < p \leq z_1} \left(1 + \frac{1}{p}\right) \ll
\left(\frac{\log z_1}{\log z_0}\right)^{p_j}\]
Sea
\[\begin{aligned}
\omega_j(n) &= \text{n\'umero de divisores primos de $n$ en $P_j$}\\
 &= |\{p\in P_j : p|n\}|\end{aligned}\]
Queremos estimar
\[\Prob((a_j - \epsilon) \log \log N < \omega_j(n) \leq
(a_j + \epsilon) \log \log N\;\;\;\forall j)
\]
Por ejemplo: si $p_1 = \frac{1}{6}$, $p_2 = \frac{1}{2}$,
$p_3 = \frac{1}{3}$, cu\'al es la probabilidad que $N$ tenga
$\sim \frac{1}{2} \log \log N$, $\sim 5 \log \log N$ y
$\sim \frac{1}{10} \log \log N$ divisores de cada tipo, respectivamente?
\begin{enumerate}
\item Como en \S \ref{sec:bambi}, problema \ref{not:antartor}, pruebe que
\[
\sum_{n\leq N} e^{\alpha_1 \omega_1(n) + \alpha_2 \omega_2(n) +
\alpha_3 \omega_3(n)}/n \ll (\log N)^{p_1 e^{\alpha_1} + p_2 e^{\alpha_2}
+ p_3 e^{\alpha_3}} \]
para $\alpha_1$, $\alpha_2$, $\alpha_3$ cualesquiera, y, prosiguiendo
como en el mismo problema,
\begin{equation}\label{eq:algem}
\mathbb{E}\left(e^{\sum_j \alpha_j \omega_j(n)}\right) \ll
(\log N)^{\sum_j p_j e^{\alpha_j} - 1} .\end{equation}
\item Por Markov,
\[\Prob((a_j - \epsilon) \log \log N < \omega_j(n) \leq
(a_j + \epsilon) \log \log N\;\;\;\forall j)\]
es a lo m\'as
\begin{equation}\label{eq:peluc}\frac{\mathbb{E}\left(e^{\sum_j a_j \omega_j(n)}\right)}{
e^{\sum \alpha_j a_j \log \log N + O_{\sum_j a_j}(\epsilon)}}\end{equation}
para $\alpha_j$ arbitrarios. Utilice (\ref{eq:algem}) y encuentre
los $\alpha_j$ para los cuales (\ref{eq:peluc}) es m\'inimo. 
Muestre, en conclusi\'on, que 
\[\Prob((a_j - \epsilon) \log \log N < \omega_j(n) \leq
(a_j + \epsilon) \log \log N\;\;\;\forall j)
\ll (\log N)^{- I_{\vec{p}}(\vec{a}) + O(\epsilon)},
\]
donde \[I_{\vec{p}}(\vec{a}) = \left(\sum a_j \log \frac{a_j}{p_j}\right)
+ 1 - \sum_j a_j .\]
\item
Proceda como en la secci\'on presente para mostrar que
\[\Prob((a_j - \epsilon) \log \log N < \omega_j(n) \leq
(a_j + \epsilon) \log \log N\;\;\;\forall j)
\gg (\log N)^{- I_{\vec{p}}(\vec{a}) + O(\epsilon) + o(1)} .\]

La cantidad $I_{\vec{p}}(\vec{a})$ es llamada la {\em entrop\'ia relativa}
de $\vec{a}$ con respecto a $\vec{p}$.
\end{enumerate}
\end{enumerate}
%\section{Caminatas aleatorias. Ley del arcoseno.}

%Cual es la palabra tecnica? Tambien -- confirmar ``esperanza''.
%Gauss-Dirichlet
%poner en claro que no est\'as asumiendo movimiento Browniano (pero
%mencionarlo)!
%\section{Ley del logaritmo iterado.}

%\section{Los grandes divisores: Gauss-Dirichlet.}
%o: convergencia de Poisson?
%Stein-Chen...
%tanto numero de grandes divisores como $\omega(n)$ (pero:Sathe-Selberg)
%connexion con el hecho que (1-rho_p) no converge exactamente a Poisson?

%\include{aplic}
\appendix
%    Include appendix "chapters" here.
\chapter{Rudimentos de probabilidades}\label{ap:proba}

Un {\em evento aleatorio} $E$ es algo que puede ya sea ocurrir o no: digamos,
la lluvia de ma\~nana. Los ejemplos extremos son los eventos con probabilidad
$0$ -- es decir, los que se sabe con certeza que no ocurrir\'an -- y los eventos
con probabilidad $1$ -- es decir, los que se sabe con certeza que ocurrir\'an.
Todo otro evento tiene probabilidad entre $0$ y $1$. Denotamos la probabilidad
del evento $E$ mediante $\Prob(E)$.

Una {\em variable aleatoria} $X$ puede tomar cualquier valor dentro de un
conjunto. Los casos m\'as comunes son las variables que toman valores
dentro de un conjunto finito o infinito de enteros (``variables discretas'')
 y las variables que toman
 valores dentro de los reales o alg\'un
otro espacio vectorial dado (``variables continuas'').
La cantidad de lluvia que caer\'a ma\~nana es un ejemplo del segundo
tipo de variable; el n\'umero de d\'ias de lluvia del a\~no pr\'oximo
es un ejemplo del primero. 
Los eventos aleatorios son, claro est\'a, un caso particular del primer tipo:
pueden verse como las variables que toman los valores $0$ y $1$, o ``no'' y 
``s\'i''. Tales variables son llamadas {\em variables de Bernoulli}:\index{variables de Bernoulli}
\[X = \begin{cases} 1 &\text{con probabilidad $p$}\\
0 &\text{con probabilidad $1-p$.}
\end{cases}\]

La {\em funci\'on de cuant\'ia} o {\em funci\'on de probabilidad}\index{funci\'on de probabilidad}
de una variable discreta $X$ 
es la funci\'on $f$ que asigna a cada valor posible $x$ su probabilidad
$f(x) = \Prob(X=x)$. 
La suma $\sum f(x)$ siempre es $1$, ya que la probabilidad que la variable
tome alguno de los posibles valores es $1$.
La funci\'on de cuant\'ia de una variable de Bernoulli es
\[f(x) = \begin{cases} p &\text{para $x=1$}\\
1-p &\text{para $x=0$}\\
0 &\text{para todo otro $x$.}\end{cases}\]
La funci\'on de cuant\'ia de la variable ``cara de un dado'' ser\'ia
\[f(x) = \begin{cases} 1/6 &\text{si $x\in \{1,2,3,4,5,6\}$}\\
0 &\text{de otra manera,} \end{cases}\]
a menos, por supuesto, que el dado est\'e trucado.

Una variable continua generalmente toma cada uno
de sus valores posibles con probabilidad $0$: la probabilidad que caigan
exactamente $\pi$ cent\'imetros de lluvia ma\~nana es cero, o infinitesimal. 
Sin embargo, una tal variable a\'un puede ser descrita por una funci\'on
de probabilidad, llamada, en este caso, 
{\em funci\'on de densidad}.\index{funci\'on de densidad}
 Digamos que la variable en cuesti\'on toma valores
en $\mathbb{R}$. La funci\'on de densidad de una tal variable
es una funci\'on $f:\mathbb{R}\to \mathbb{R}_0^+$ cuya integral es $1$.
La probabilidad que la variable tome su valor entre $a$ y $b$
est\'a dada por la integral
\[\int_a^b f(x) dx .\]

\begin{figure}
\centering \includegraphics[height=1.5in]{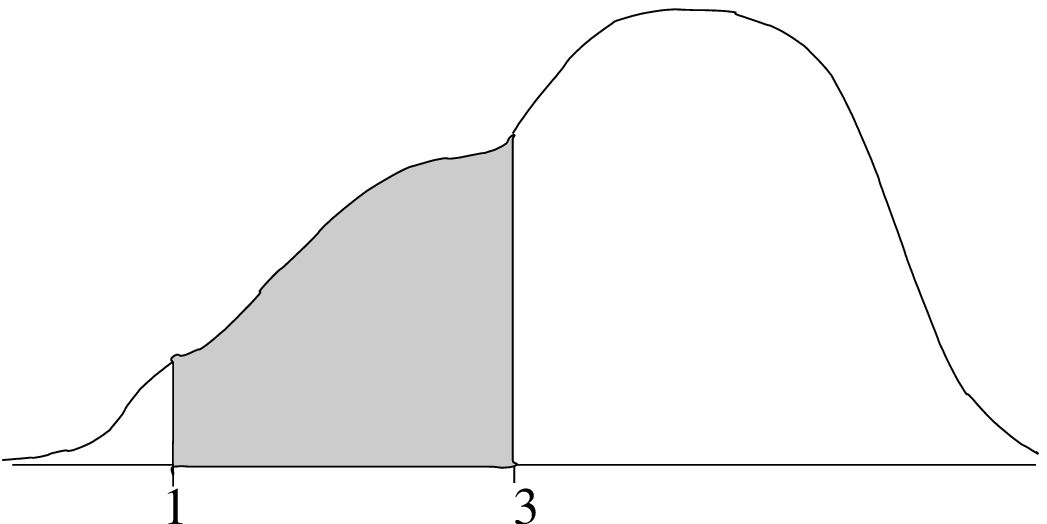}
\caption{Distribuci\'on de una variable continua. La probabilidad que
la variable tome un valor entre $1$ y $3$ es igual al \'area
marcada.}
\end{figure}

\begin{figure}
\centering \includegraphics[height=1.115in]{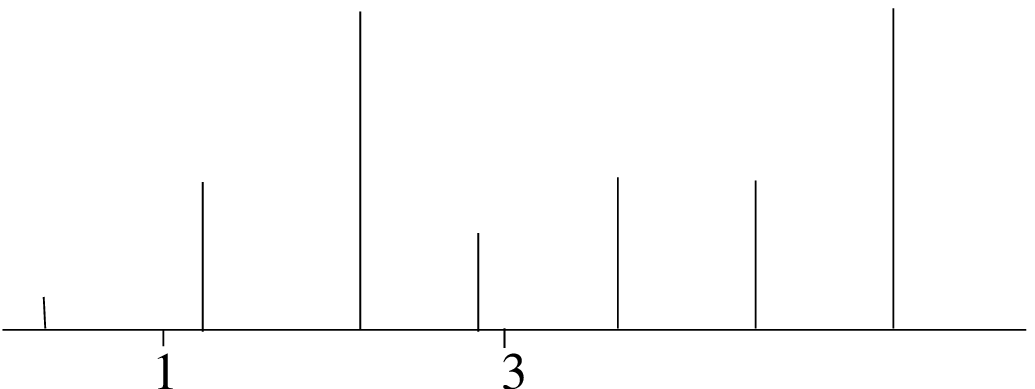}
\caption{Distribuci\'on de una variable discreta. La probabilidad que
la variable tome un valor entre $1$ y $3$ es igual a la
suma de los valores $f(x)$ de la funci\'on de distribuci\'on
para $x$ entre $1$ y $3$, o, lo que es lo mismo, a la suma
de las alturas de las barras entre $1$ y $3$.}
\end{figure}

Decimos que una variable $X$ tiene la {\em distribuci\'on uniforme}\index{distribuci\'on uniforme}
si todos sus valores son igualmente probables. Tanto en el caso continuo
como en el caso discreto, $X$ tiene la distribuci\'on uniforme si y s\'olo
si su funci\'on $f$ es una funci\'on constante.
Tanto un dado justo como una moneda justa tienen la distribuci\'on uniforme --
si se define el dominio como $\{1,2,\dotsc,6\}$ y
$\{\text{cara}, \text{sello}\}$ en el otro, claro est\'a.

\begin{small}
Muy a menudo, los mismos enunciados y las mismas pruebas valen para
las variables discretas y continuas si se utilizan sumas
en un caso e integrales en el otro. Tambi\'en puede haber variables
con un rango en parte continuo y en parte discreto. Por ello, lo correcto
es tener un s\'olo marco para todas las distribuciones, de tal manera
que la distinci\'on entre las variables discretas y las variables continuas 
desaparezca en el plano formal. 
El lector puede adivinar que tal marco nos es dado por la 
integraci\'on de Lebesgue; en dicha perspectiva, las sumas son un
caso particular de las integrales. 

Tal es el formalismo aceptado en estos d\'ias, por excelentes razones. 
Empero, no nos preocuparemos, y hablaremos como si nuestras
variables fueran discretas o continuas dependiendo de lo que haga que nuestra
notaci\'on sea m\'as conveniente. 

Una alternativa elemental consiste en definir la {\em funci\'on de
distribuci\'on acumulada} $P_X(x)$ de la siguiente manera:
\begin{equation}\label{eq:cafe}
P_X(x) = \Prob(X\leq x).
\end{equation}
En el caso continuo, $P(x) = \int_{-\infty}^x f(t) dt$; en el caso 
discreto, $P(x) = \sum_{t\leq x} f(t)$. Como la definici\'on (\ref{eq:cafe})
es v\'alida en ambos casos, uno puede utilizar $P_X(x)$ en vez de
$f(t)$ y de esta manera hablar de ambos tipos de
distribuci\'on a la vez.

No utilizaremos la funci\'on de distribuci\'on acumulada. Cuando decimos
``la distribuci\'on $f(x) = e^{-x}$'' o ``la distribuci\'on
$f(n) = \frac{1}{e \cdot n!}$'', queremos decir ``la distribuci\'on
continua con funci\'on de 
densidad $f(x) = e^{-x}$'' o ``la distribuci\'on discreta dada por la 
funci\'on de cuant\'ia $f(n) = \frac{1}{e \cdot n!}$'', respectivamente;
el caso que se tiene en mente estar\'a claro en el contexto.
\end{small}

\begin{center}
* * *
\end{center}

La {\em probabilidad condicional} $\Prob(E_1 | E_2)$ es la probabilidad 
de un evento $E_1$ dado que el evento $E_2$ ocurre. Tenemos
\[\Prob(E_1|E_2) = \frac{\Prob(E_1\wedge E_2)}{\Prob(E_2)},\]
donde $E_1\wedge E_2$ se define como el evento que tanto $E_1$ como $E_2$
ocurran. (Si la probabilidad que llueva ma\~nana es $0.1$, y la probabilidad
que llueva y enfr\'ie es $0.07$, la probabilidad que enfr\'ie, dado que
llover\'a, es $0.7$.)\index{probabilidad condicional}

Decimos que dos eventos $E_1$, $E_2$ son {\em independientes} si\index{independencia!de eventos}
\[\Prob(E_1|E_2) = \Prob(E_1) \text{\;\; y\;\;} \Prob(E_2|E_1) = \Prob(E_2).\]
Decimos que dos variables $X$, $Y$ son independientes si
\[\Prob(X=x|Y=y) = \Prob(X=x) \text{\;\; y\;\;} \Prob(Y=y|X=x) = \Prob(Y=y) .
\]\index{independencia! de variables}
para $x$, $y$ cualesquiera. En otras palabras, $X$ y $Y$ son independientes
si el valor tomado por una no nos dice nada acerca del valor de la otra.
(Digamos: saber que ma\~nana llover\'a en Alaska no nos ayuda a saber si
es que ma\~nana llover\'a en Iquitos.)

Si tenemos variables $X_1, X_2, \dotsc, X_n$, decimos
que son {\em independientes en pares} si $X_i$ y $X_j$ son independientes
\index{independencia!en pares}
para $1\leq i,j,\leq n$, $i\ne j$ cualesquiera. Decimos que 
$X_1, X_2, \dotsc, X_n$ son {\em mutuamente independientes} si
\index{independencia!m\'utua}
\[\Prob(X_i = x_i | X_1 = x_1, \dotsc, X_{i-1} = x_{i-1}, X_{i+1} = x_{i+1},
\dotsc, X_n = x_n)\]
para $x_1, x_2,\dotsc, x_n$ e $i$ cualesquiera. 

\begin{xca}
Muestre que, si $X_1, X_2,\dotsc, X_n$ son mutuamente independientes,
entonces son independientes en pares.
\end{xca}

Sin embargo, si $X_1, X_2,\dotsc, X_n$ son independientes en pares, ello no
es de ninguna manera suficiente para que sean mutuamente independientes.

\begin{center}
* * *
\end{center}

Sea $X$ una variable que toma valores dentro de los reales (o los complejos).
La {\em esperanza} $\mathbb{E}(X)$ de $X$ es\index{esperanza}
%VER TERMINOLOGIA!
\[\mathbb{E}(X) = \sum_x \Prob(X = x) \cdot x \]
(o $\int_{-\infty}^{\infty} f(x) x dx$ en el caso continuo, donde $f$
es la funci\'on de densidad).
En otras palabras, se trata de una especie de promedio. 
%Por ejemplo, si $X$ es
%un dado (sin truco), entonces $\mathbb{E}(X) = 3.5$. Si $X$ estuviera
%trucado para darnos $6$ la mitad del tiempo, y tomara los otros valores
%con igual probabilidad, tendr\'iamos $\mathbb{E}(X) = 4.5$.

%(Como de costumbre, estamos trabajando en el caso discreto, en parte para
%que nuestra notaci\'on sea sencilla. En el caso continuo, tendr\'iamos
%\[\mathbb{E}(X) = \int f(x) x dx,\]
%donde $f$ es la funci\'on de distribuci\'on de $X$.)

\begin{xca}
Muestre que, si $X$ e $Y$ son independientes, entonces
\[\mathbb{E}(X Y) = \mathbb{E}(X) \mathbb{E}(Y).\] 
\end{xca}
\'Esta es una condici\'on necesaria, pero no suficiente, para que $X$ e $Y$
sean independientes.

Sean $X$ una variable y $E$ un evento. 
Definimos la {\em esperanza condicional} \index{esperanza condicional}
$\mathbb{E}(X|E)$ de $X$ dado $E$ de la siguiente manera:
\[\mathbb{E}(X|E) = \sum_x \Prob(X = x|E) \cdot x .\]

\begin{center}
* * *
\end{center}

Existen diversas maneras de describir una variable, m\'as all\'a de su
distribuci\'on (que nos da una descripci\'on completa) y su esperanza.
La m\'as com\'un es la {\em varianza}:\index{varianza}
\[\Var(X) = \mathbb{E}((X - \mathbb{E}(X))^2) .\]
La varianza ser\'a grande si los valores de $X$ tienden a alejarse mucho
de la esperanza de $X$, y peque\~na si esto no sucede. 

La {\em desviaci\'on est\'andar} $\sigma(X)$ no es sino la ra\'iz cuadrada
de la varianza:\index{desviaci\'on estandar}
\[\sigma(X) = \sqrt{\Var(X)} .\]
La desviaci\'on est\'andar nos da una buena idea de la escala de las
desviaciones -- esto es, las distancias de los valores de $X$ de la esperanza
de $X$. Puede suceder que $X$ est\'e a una o dos desviaciones est\'andar de
su esperanza gran parte del tiempo, pero podr\'a estar a m\'as de diez
desviaciones est\'andar de distancia de su esperanza a lo m\'as una de cada
$100$ veces ({\em desigualdad de Chebyshev}).\index{desigualdad!de Chebyshev}

La varianza y la desviaci\'on est\'andar no distinguen entre las veces
en que $X$ toma valores m\'as grandes y m\'as peque\~nos que su esperanza.
Hablamos de la {\em cola superior} de la distribuci\'on cuando queremos referirnos a aquellos
valores posibles de $X$ que son mucho m\'as grandes que $\mathbb{E}(X)$;
decimos {\em cola inferior} para referirnos a los valores que son mucho
m\'as peque\~nos que $\mathbb{E}(X)$. 

%Decimos, coloquialmente, que una distribuci\'on es de {\em cola ligera}
%si $\Prob(|X - \mathbb{E}(X)|>t)$) decrece muy r\'apidamente a medida que
%$t\to \infty$; de lo contrario, decimos que se trata de una distribuci\'on
%de {\em cola pesada}. El ejemplo arquet\'ipico de una distribuci\'on
%de cola ligera es la normal $f(x) = \frac{1}{\sqrt{2\pi}} e^{-x^2/2}$;
%dos distribuciones de cola pesada muy comunes son la exponencial
%($f(x) = \alpha e^{-\alpha x}$, $x\geq 0$, para alg\'un $\alpha>0$ fijo;
% distribuci\'on continua)
%y la distribuci\'on de Poisson ($f(n) = e^{-\lambda} \frac{\lambda^n}{n!}$,
%$n=0,1,2,\dotsc$, para alg\'un $\lambda>0$ fijo; distribuci\'on discreta).
%Podemos ver que ``cola pesada'' es un t\'ermino relativo; fuera de la
%teor\'ia de las probabilidades,
%dir\'iamos que $f(x) = e^{-x}$ decrece r\'apidamente -- es s\'olo que
%decrece lentamente comparada con $e^{-x^2/2}$.
\begin{center}
* * *
\end{center}

Supongamos que tenemos dos variables $X$ e $Y$ con la misma funci\'on
de densidad, excepto por un argumento; digamos, por ejemplo,
que la funci\'on de densidad de $X$ es 
$g(x) = \begin{cases} 1 &\text{si $0\leq x\leq 1$}\\0 
&\text{de otra manera,}
\end{cases}$ y la funci\'on de densidad de $Y$ es
$h(x) = \begin{cases} 1 &\text{si $0<x\leq 1$}\\$0$ &
\text{de otra manera}
\end{cases}$. Entonces tiene sentido decir que $X$ e $Y$ poseen la misma
distribuci\'on: la probabilidad que $X$ sea exactamente $1$ es 
infinitesimal, de todas maneras, y la probabilidad que $X$ est\'e entre
$a$ y $b$ (para $a<b$ cualesquiera) es igual a la probabilidad que
$Y$ est\'e entre $a$ y $b$.

Esto sugiere la definici\'on siguiente.
Sean dadas una variable $Z$ y una sucesi\'on de variables
$Z_1, Z_2, Z_3,\dotsc$.
Decimos que las variables $Z_1, Z_2, Z_3,\dotsc$ 
{\em convergen en distribuci\'on a $Z$} si, para $a<b$ cualesquiera,
\begin{equation}\label{eq:ahorita}
\lim_{n\to \infty} \Prob(a< Z_n< b) = \Prob(a<Z<b).\end{equation}\index{convergencia en distribuci\'on}

\begin{small}
Desde el punto de vista de la integraci\'on de Lebesgue, esta no es
sino la ``convergencia d\'ebil''. Por qu\'e? Ahora bien, si una sucesi\'on
de funciones de densidad (o cuant\'ia) $f_n$ converge en el sentido de la
convergencia d\'ebil a una funci\'on de densidad (o cuant\'ia) $f$, la
sucesi\'on de funciones de distribuci\'on acumulada
  $P_n(x) = \int_{-\infty}^x f_n(t) dt$ 
converge para cada $x$ a la funci\'on de distribuci\'on acumulada
$P(x) = \int_{-\infty}^x f_n(t) dt$. (Esto no es sino
(\ref{eq:ahorita}).) No es dif\'icil probar que la convergencia
de $P_n(x)$ a $P(x)$ es, incluso, uniforme en $x$. (Utilize
$f_n(t), f(t)\geq 0$ y $\int_{-\infty}^\infty f_n(t) dt =
\int_{-\infty}^{\infty} f(t) dt = 1$.)
\end{small}

Muy bien puede suceder que las variables $Z_n$ sean discretas, y que su
l\'imite $Z$
una variable continua. 

\chapter{Comentarios finales}

En esta introducci\'on a la teor\'ia probabil\'istica de n\'umeros, nos
centramos en el estudio de los divisores primos de un n\'umero aleatorio,
y no en el estudio de un primo aleatorio, o de los primos cercanos
a un entero aleatorio. La raz\'on principal es que a\'un se sabe poco
con certeza acerca de estas otras preguntas. Por ejemplo,
se conjetura, pero no se sabe, la probabilidad que un entero aleatorio $n$
sea primo y que $n+2$ tambi\'en sea primo (\S \ref{sec:varianza},
nota \ref{it:gotor}). Las m\'as de las veces, lo \'unico que se posee es
cotas superiores dadas por la teor\'ia de cribas. Existen modelos tanto
fruct\'iferos como imperfectos - por ejemplo, el modelo de Cram\'er\index{modelo
de Cram\'er} (ver, por ejemplo, \cite{G}), que
dice que el evento que un n\'umero $n$ sea primo y el evento que un n\'umero
$m$ distinto sea primo 
se comportan muchas veces como si fueran eventos independientes. Est\'a claro
que esto no debe ser cre\'ido completamente: por ejemplo, si $m=n+1$ y $n>2$,
los n\'umeros $n$ y $m$ no pueden ser ambos primos, lo cual ser\'ia 
una posibilidad si estuvieramos hablando de variables independientes.

En todo el texto, obedecimos 
a una restricci\'on autoimpuesta: nos abstuvimos de usar 
teor\'ia de la medida y an\'alisis complejo. Para proseguir en el
estudio de la teor\'ia probabil\'istica de n\'umeros, es necesario usar los dos.
El an\'alisis complejo es sumamente \'util para el estudio de los primos. 
Toda la teor\'ia anal\'itica de n\'umeros depende del an\'alisis; se trata de
un caso cl\'asico de c\'omo el estudio de lo continuo puede ayudar en el 
estudio de lo discreto. La idea principal es que, para estudiar sumas finitas,
como $\sum_{n\leq N} \Lambda(n)$, $N$ variable, debemos estudiar sumas infinitas
$\sum_n \Lambda(n) n^{-s}$, $s$ variable. Estas sumas infinitas se tratan
como funciones complejas de $s$, generalmente anal\'iticas o merom\'orficas.

La teor\'ia de la medida se considera hoy en d\'ia como necesaria para
desarrollar la teor\'ia de probabilidades sobre una base rigurosa. Si se
profundiza en el estudio de la teor\'ia probabil\'istica de n\'umeros
sobre una base puramente intuitiva, se llega f\'acilmente al punto donde
el lenguaje mismo falta. Veamos, por ejemplo, el caso de 
las caminatas aleatorias.
Tomemos un n\'umero $n$ al azar. Consideremos los primos $p=2,3,5,\dotsc$
en orden. En cada paso, si $p$ divide $n$, damos un paso a la derecha;
si $p$ no divide $n$, damos un paso mucho m\'as corto a la izquierda. 
Billingsley \cite{Bi1} prob\'o que la caminata que resulta tiende en 
distribuci\'on\index{convergencia en distribuci\'on}
al mismo l\'imite que una caminata aleatoria, es decir, el movimiento
Browniano. Ahora bien, qu\'e quiere decir que una caminata ``tiende'' a un
cierto tipo de crecimiento en distribuci\'on? No se trata simplemente
de una sucesi\'on de n\'umeros que convergen a otro n\'umero. Resulta dif\'icil
formular el resultado -- ni que decir de su prueba -- sin la teor\'ia
de la medida.

Hay dos desarrollos recientes notables. En primer lugar, la teor\'ia de
cribas ha probado ser m\'as flexible y potente de lo que se cre\'ia hasta
ahora, si se suplementa con otros m\'etodos; est\'an all\'i resultados
inesperados en la teor\'ia anal\'itica de n\'umeros obtenidos
por Friedlander-Iwaniec, Goldston-Pintz-Yildirim, y otros. En segundo lugar,
la teor\'ia erg\'odica no s\'olo esta iluminando la teor\'ia de n\'umeros
-- inclu\'ido el estudio de los primos -- sino que esta haciendo posibles
pruebas de resultados que no se cre\'ian anteriormente accesibles. 
Un ejemplo muy reciente e impresionante es el teorema de Green y Tao sobre
los n\'umeros primos en progresiones aritm\'eticas; no es que su prueba
establezca propiedades
sumamente delicadas de los n\'umeros primos, sino m\'as
bien que muestra que algunas leyes muy precisas no son delicadas, al 
punto que deben regir para todo subconjunto de los enteros con ciertas propiedades
generales -- propiedades que los primos satisfacen. Toda la teor\'ia erg\'odica
esta basada sobre la teor\'ia de la medida, y ser\'ia imposible comenzar
el estudio de la primera sin utilizar la segunda.

La bibliograf\'ia tiene como fin ser \'util antes que completa. El libro
de Tenenbaum \cite{T} es una introducci\'on est\'andar y detallada
 al tema, con mucho m\'as an\'alisis que la presente monograf\'ia. El
libro de Iwaniec y Kowalski \cite{IK} se ha vuelto la obra can\'onica de la
teor\'ia anal\'itica de n\'umeros para la \'epoca actual. Feller \cite{Fe}
es un texto cl\'asico de probabilidades del cual generaciones se han beneficiado. Todos los art\'iculos citados est\'an entre los esenciales sobre el tema;
la mayor parte de ellos son de lectura razonablemente accesible.

%tambien: conceptos basicos de probabilidades
\backmatter
%    Bibliography styles amsplain or harvard are also acceptable.

\printindex
\end{document}